\documentclass[12pt,leqno]{article} 
\newcommand{\mylabel}[1]{\label{#1}}

\newcommand{\eop}{\hfill $|||$ \vspace{.1in}}

\newtheorem{lemma}{Lemma}
\newtheorem{corollary}[lemma]{Corollary}
\newtheorem{proposition}[lemma]{Proposition}

\newtheorem{theorem}{Theorem}
\newtheorem{conjecture}{Conjecture}

\newcommand{\cc}{{\bf C}}
\newcommand{\rr}{{\bf R}}
\newcommand{\qq}{{\bf Q}}

\newcommand{\zz}{{\bf Z}}
\newcommand{\pp}{{\bf P}}
\newcommand{\AAA}{{\bf A}}

\newcommand{\Cc}{{\cal C}}
\newcommand{\Ee}{{\cal E}}
\newcommand{\Ff}{{\cal F}}
\newcommand{\Gg}{{\cal G}}

\newcommand{\Tt}{{\cal T}}
\newcommand{\Oo}{{\cal O}}
\newcommand{\Dd}{{\cal D}}

\newcommand{\Bb}{{\cal B}}
\newcommand{\Mm}{{\cal M}}

\newcommand{\Kk}{{\cal K}}
\newcommand{\Vv}{{\cal V}}
\newcommand{\Uu}{{\cal U}}
\newcommand{\Ll}{{\cal L}}
\newcommand{\Hh}{{\cal H}}

\newcommand{\Nn}{{\cal N}}
\newcommand{\Pp}{{\cal P}}

\newcommand{\Yy}{{\cal Y}}
\newcommand{\Zz}{{\cal Z}}

\newcommand{\Gm}{{\bf G}_m}
\newcommand{\Aone}{{\bf A}^1}
\newcommand{\Aquot}{{\cal A}}

\begin{document}

\section*{Nonabelian mixed Hodge structures}

\medskip

\noindent
Ludmil Katzarkov
\footnote{University of California at Irvine, Irvine, CA 92697, USA. 
Partially supported by NSF Career Award  DMS-9875383 and A.P. 
Sloan research fellowship.},
Tony Pantev
\footnote{University of Pennsylvania, 
209 South 33rd Street
Philadelphia, PA 19104-6395, USA.
Partially supported by NSF Grant DMS-9800790 and A.P. Sloan Research Fellowship.},
Carlos Simpson
\footnote{CNRS, Universit\'e de Nice-Sophia Antipolis, Parc Valrose, 
06108 Nice Cedex
2, France.}

{\scriptsize
\noindent
{\bf Abstract}---
We propose a definition of ``nonabelian mixed Hodge structure'' together
with a construction associating to a smooth projective variety $X$ and
to a nonabelian mixed Hodge structure $V$, the ``nonabelian cohomology
of $X$ with coefficients in $V$'' which is a (pre-)nonabelian mixed
Hodge structure denoted $H=Hom(X_M, V)$. We describe the basic
definitions and then give some conjectures saying what is supposed to
happen. At the end we compute an example: the case where $V$ has 
underlying homotopy type the complexified $2$-sphere, and mixed Hodge 
structure coming from its identification with $\pp ^1$. For this example
we show that $Hom (X_M,V)$ is a namhs for any smooth projective variety $X$. 
}

\medskip

\noindent
{\bf Introduction}---p. \pageref{intropage}
\newline
{\bf Part I: Nonabelian weight filtrations}
\newline
{\sc Conventions}---p. \pageref{conventionspage}
\newline
{\sc Nonabelian filtrations}---p. \pageref{nonabfiltpage}
\newline
{\sc  Perfect complexes and Dold-Puppe linearization}---p. \pageref{doldpuppepage}
\newline
{\sc Further study of filtered $n$-stacks}---p. \pageref{furtherstudypage}
\newline
{\sc Filtered and $\Gm$-equivariant perfect complexes}---p. \pageref{gmequivariantpage}
\newline
{\sc Weight-filtered $n$-stacks}---p. \pageref{weightfilteredpage}
\newline
{\sc Analytic and real structures}---p. \pageref{analyticrealpage}
\newline
{\bf Part II: The main definitions and conjectures}
\newline
{\sc Pre-nonabelian mixed Hodge structures}---p. \pageref{pnamhspage}
\newline
{\sc Mixed Hodge complexes and linearization}---p. \pageref{mhcpage}
\newline
{\sc Nonabelian mixed Hodge structures}---p. \pageref{namhspage}
\newline
{\sc Homotopy group sheaves}---p. \pageref{htygrpshpage}
\newline
{\sc The basic construction}---p. \pageref{basicconstructionpage}
\newline
{\sc The basic conjectures}---p. \pageref{basicconjecturespage}
\newline
{\sc Variations of nonabelian mixed Hodge structure}---p. 
\pageref{fvnamhspage}
\newline
{\bf Part III: Computations}
\newline
{\sc Some morphisms between Eilenberg-MacLane pre-namhs}---p. 
\pageref{basicpage}
\newline
{\sc Construction of a namhs $\Vv$ of homotopy type $S^2_{\cc}$}---p. 
\pageref{aconstructionpage}
\newline
{\sc The namhs on cohomology of a smooth projective variety with
coefficients in $\Vv$}
---p. \pageref{thenamhspage}

\newpage

\begin{center}
{\Large \bf  Introduction}
\end{center}
\mylabel{intropage}

The nonabelian aspect of homotopy theory has been present in mixed Hodge
theory from very early on. This started with the work of
Deligne-Griffiths-Morgan-Sullivan \cite{dgms} on rational homotopy
theory, and then Morgan \cite{Morgan} put a mixed Hodge structure on the
differential graded algebra which calculates the homotopy groups of a
simply connected complex variety (or the nilpotent completion of the
fundamental group of a
non-simply connected one). This construction was later improved by Hain
\cite{Hain} \cite{Hain2} who in particular obtained a much better
functoriality. 

Many applications of the Morgan-Hain construction have been given, for example in
Torelli-type theorems \cite{ccm}, \cite{Carlson}, the use
of the strictness property \cite{Lasell}, in the study of algebraic
cycles \cite{Green}, \cite{Voisin}, for the Shafarevich conjecture for
nilpotent fundamental groups \cite{Katzarkov} \cite{Leroy}, etc.
There have recently been new attempts to better understand what is
going on, see \cite{Arapura2}, \cite{Pearlstein}, and to simplify the relatively abstract
nature of Hain's construction, see \cite{Leroy}, \cite{Kaenders}.

In an important development, Hain obtained a mixed Hodge
structure on the {\em relative Malcev completions} of the fundamental
group at any representation underlying a variation of Hodge structure \cite{Hain3}.
This significantly extends the ``amount'' of the fundamental group that
is seen by mixed Hodge theory. 

The variation of the mixed Hodge structures on homotopy groups, when the
underlying variety varies in a family, has been studied by Navarro-Aznar
\cite{NavarroAznar} \cite{NavarroAznar2}
and Wojtkowiac \cite{Wojtkowiac}.

There are two basic gaps in the classical mixed Hodge homotopy theory.
The first is that one should really treat a ``higher functoriality'':
the set of homotopy $n$-types forms an $n+1$-category, so one would like
to have an $n+1$-category of ``homotopy $n$-types with mixed Hodge
structure''. The fact that we will then have an $n+1$-functor from the
$1$-category of varieties, to this $n+1$-category, contains subtle
additional ``homotopy coherence'' data. 

The second gap is that up until fairly recently, 
very little has been done on the problem of the
higher homotopy theory of non-simply connected varieties. A first look
at this problem may be considered to come from the work of
Green-Lazarsfeld \cite{GreenLazarsfeld} on cohomology jump loci, 
although they don't talk about
mixed Hodge structures. Arapura has looked at the Green-Lazarsfeld
results using mixed 
Hodge structures \cite{Arapura}. Kapranov \cite{Kapranov} and Toen
\cite{Toen}
also have methods which go in the direction of looking at the higher 
homotopy theory of non-simply connected varieties, but these haven't
been
fully developed yet.

The third author of the present paper has proposed, as a remedy to the
above problems, to look at ``nonabelian cohomology''. One chooses a
``coefficient  object'' which is an $n$-stack $T$, and looks at the {\em
nonabelian cohomology} $H=\underline{Hom}(X_B, T)$. The coefficient object
encodes data which is much the same as that of a complex differential
graded algebra, however on the level of the fundamental group of $T$ we 
have $\pi _1(T)=G$ an affine group scheme. The cohomology stack $H$ now
lies over the moduli $1$-stack ${\cal M}(X_B, G)
=\underline{Hom}(X_B, BG)$ of representations of $\pi
_1(X^{\rm top})$ in $G$, and $H$ encodes higher homotopy data relative
to the universal family of representations of $\pi _1 (X^{\rm top})$.
In this way we get an approach towards the problem of higher homotopy in
the presence of a fundamental group. The construction $X \mapsto
\underline{Hom}(X_B,T)$ is automatically $n+1$-functorial in $X$ (in $T$
too), so this addresses the issue of ``higher functoriality''.
Furthermore, this higher functoriality shows up in concrete ways, for
example in the case of a family of varieties one obtains a secondary
Kodaira-Spencer map \cite{secondaryKS} which detects a ``higher
variation'' of Hodge structure even when all of the classical variations
of (mixed) Hodge structure associated to the family are constant. 

One of the  major problems in this nonabelian cohomology approach has been
the lack of a notion of weight filtration, and consequently the lack of
a notion of mixed Hodge structure. The present paper aims to start to
close this gap. Briefly, what we do is to propose a definition of
nonabelian mixed Hodge structure and in particular, the notion of ``weight
filtration'' that goes into this. We state a number of conjectures about
the notion that we define, but only go a short way toward justifying
these conjectures by explicitly treating a particular example in Part
III. 

Our search for a notion of nonabelian mixed Hodge structure is 
motivated by the possibility of obtaining applications similar to 
what has already been done with mixed Hodge structures on homotopy groups.
One of the main things would be to obtain restrictions on the homotopy
types of smooth projective varieties (or even general quasiprojective varieties).
This might include generalized ``quadratic singularity'' statements, for 
nonabelian cohomology stacks.
Other possibilities might include Torelli-type theorems, strictness results,
and restrictions
on the monodromy actions (or generalized monodromy actions) for holomorphic
families of smooth projective varieties. 
One might get a nonabelian version of the theory of realizations of
motives \cite{Huber}.
Finally, we think that the techniques we
use here might also be of interest in generalizing the geometric Langlands
correspondence to higher nonabelian cohomology.

Before saying more about the ideas that we present here, we
would like to take note of a few other directions that one could
possibly take. It should be possible to remedy the problem of ``higher
functoriality''
by looking directly at Morgan's or Hain's construction and consistently
keeping
track of what is going on in a suitable $n+1$-categorical or
$\infty$-categorical way. The work of Hinich \cite{Hinich} would
probably be relevant here. 

Arapura's treatment of the Green-Lazarsfeld results on jump loci
\cite{Arapura} makes explicit use of mixed Hodge
structures on dgas and discusses higher homotopy in the sense that he treats the jump loci
for higher cohomology.  He also has some results for support loci for higher rank
local systems \cite{Arapura3}. These results are a step in the direction
indicated in the previous paragraph.

It might be possible to construct a ``big''
differential graded algebra which encloses mixed-Hodge theoretic
information about higher homotopy, relative to all the
representations of $\pi _1$. The first and second authors have been
thinking about this for some time now and hope to treat
it in a future paper. 

There is a relationship between simply  connected differential
graded algebras and simply connected geometric $n$-stacks, and we have
to some extent exploited this in our research on the questions treated
in the present paper; however, our calculations in this direction remain
largely heuristic. It would be good to understand rigorously the foundations
of these calculations; this would help to bridge the gap with the
earlier results of Morgan and Hain.

The big differential graded algebra approach might also be related to the notion of higher
Malcev completion with its conjectural mixed Hodge structure 
mentionned at the end of \cite{limits}. On the downside,
these objects seem to be too ``big'' (see the example given at the beginning of
\cite{aaspects}) and in a certain sense, one could think of the nonabelian
cohomology approach as being a way of extracting finite dimensional
information from these big objects.

Another approach which uses the notion of differential graded algebra
but in a slightly different manner, is Kapranov's ``extended moduli
space of representations'' \cite{Kapranov}, part of the idea for which
apparently comes from Kontsevitch's folkloric ``derived deformation
theory''. Kapranov's extended moduli space visibly contains higher
homotopy data, because the tangent complex at a point includes all
higher cohomology groups. We have not been able to identify any further
the exact nature of this homotopy-theoretic information.
It seems quite likely that there should be some type of ``mixed Hodge
structure''
on Kapranov's space, and that this would include a  mixed Hodge structure on
whatever higher homotopy information is contained therein. 

For us, an interesting application of Kapranov's idea would
be if one could define a notion of ``extended $n$-stack'' in order to be
able to have an ``extended nonabelian cohomology stack'' which would
extend $\underline{Hom}(X_B, T)$ in the same way that Kapranov's object
extends $\underline{Hom}(X_B, BG)$. The reason why this would be helpful
is that the extended  object could be expected to be formally smooth
(this hope seems to be due to Kontsevitch, Deligne, Drinfeld). This
smoothness would serve to avoid a large number of the problems which we
have with the annihilator ideals in the present paper. 
Thus, while representing a further abstraction, such a combination of
ideas could end up eventually simplifying things.

Yet another direction is that taken recently by Toen \cite{Toen}. He
proposes to do an $\infty$-categorical version of the idea of Tannaka
duality, to view a homotopy type as the Tannaka dual of its $\infty$-category
of local systems of perfect complexes with monoidal structure (tensor
product). In this point of view, Toen expects to be able to define
Hodge and weight filtrations and to obtain some type of mixed Hodge
structure.  

To close this part, we note that 
one of the major advances in mixed Hodge theory,
Saito's theory of mixed Hodge modules \cite{Saito} \cite{Saito2} \cite{Saito3}, treats
mostly the abelian cohomology theory. It would be good to have a
``nonabelian''
version of Saito's theory directed toward a notion of Hodge theory for
nonabelian cohomology with local coefficients varying constructibly
with respect to a stratification.
This seems to be a long way off and that
could be considered as a long-range goal of any of the above approaches.

Now we get to the ideas that are contained in the present paper. 
The domain stack $X_M$ which we shall use is 
an object enclosing $X_B$, $X_{DR}$ and $X_{Dol}$, see p.
\pageref{basicconstructionpage} for more
details. 
In the abelian theory the idea of combining together several different
types
of cohomology theories goes back a long way, see Huber \cite{Huber} for
a compendium. 
Several people including B. Toen and M. Rosellen had asked us
whether there could be any type of ``nonabelian motive'' taking into
account
the de Rham, Dolbeault and Betti cohomologies. While there doesn't seem
to be anything deep going on here, the notation $X_M$ is a useful
bookkeeping device. 

Our first  preliminary calculations showed that in order to have any
hope of having some kind of ``mixed Hodge structure'' on a nonabelian
cohomology stack $\underline{Hom}(X_M,T)$, the coefficient stack $T$ had
itself to have a ``mixed Hodge structure''. 
Thus our search for the
notion of nonabelian mixed Hodge structure concerned both the
coefficient stack and the answer stack. The idea is to obtain a
collection of conditions for the definition, such that if $T$ satisfies
these conditions then so does $\underline{Hom}(X_M, T)$. Actually the
set of conditions which we propose here is just a first attempt; it
could in the future turn out to be wrong either by having too many
conditions, or by having not enough.

The first and
most obvious element 
is that, in order to integrate the ``weight filtration''
into the notion of $n$-stack, we will use the Rees space approach to the
notion of filtration: a filtered object is an object over $\Aone$ together
with an action of $\Gm$ covering the standard action of $\Gm$ on 
$\Aone$. If the objects in question are vector spaces, then to 
a filtered vector space one obtains a corresponding equivariant bundle
over $\Aone$ which is essentially its {\em Rees module}.
The origins of the notion of Rees module are actually closely related to the case
of ``filtered schemes'' (cf Fulton \cite{Fulton} pp 90-91, who cites Gerstenhaber
\cite{Gerstenhaber} and Rees \cite{Rees}, treating filtered rings; 
these are all related to the 
``deformation to the normal cone'' which we discuss on p. \pageref{defts} below).

The idea of using the Rees correspondence for the Hodge filtration, and
particularly for the Hodge filtration on nonabelian cohomology, is 
discussed in \cite{naht}, \cite{santacruz}. The idea goes back to Deligne's
approach to the twistor space \cite{delignetwistor}, 
Deninger's formal version of the Rees correspondence for the Hodge filtration
\cite{Deninger}, and the Neisendorfer-Taylor 
deformation between de Rham and Dolbeault
cohomology \cite{NiesendorferTaylor}. 
Perhaps the first place where the notion of Rees bundle
appears in connection with the weight filtration, is in Sabbah \cite{Sabbah}.
Recently (at Irvine, June 1998) 
Voevodsky made the intriguing remark that the tensor category
of mixed motives had two fiber functors, one for the underlying de Rham
cohomology and one for the associated-graded of the weight filtration. These
clearly correspond to the two points of the stack $\Aquot$ defined in the
following paragraph, so this remark can also be seen as a precedent for using the
Rees space approach to look at the weight filtration.

Thus we fix the notion of weight filtration as being an equivariant
object over $\Aone$, or equivalently an object over the quotient stack
$$
\Aquot := \Aone / \Gm .
$$
It follows relatively logically that the notion of a {\em pair} of
filtrations, the weight filtration and the Hodge filtration, is to be
interpreted as being an object over the product
$$
\Aquot _{wt} \times \Aquot _{hod}
$$
where the superscripts indicate separate copies of $\Aquot$, the first one for
the weight filtration and the second one for the Hodge filtration. 
We denote a geometric $n$-stack provided with two filtrations, by 
$(V,W,F)$; this object really means a geometric $n$-stack 
$$
Tot ^{W,F}(V)\rightarrow \Aquot _{wt} \times \Aquot _{hod}.
$$

An important element in Hodge theory is the interplay between the Hodge
filtration and its complex conjugate. This reappears in the notion of
``twistor space'' if one wishes to avoid the Hodge filtration altogether.
It would certainly be possible to do what we are doing here in the 
``twistor'' context, to obtain a notion of {\em nonabelian mixed twistor
structure}. However, this would add on yet another layer of abstraction
which we don't feel is justified at the current time. Thus, for the present 
purposes, we introduce the notion of {\em real structure}, which allows  us
literally to take the complex conjugate of the Hodge filtration
for the purposes of measuring the ``purity'' of the Hodge filtration.

The {\em associated-graded} of the weight filtration is an object $Gr^W(V)$
obtained by restricting to $[0]\in \Aquot _{wt}$. The ``point'' represented
by the origin is actually the quotient stack $0/\Gm $ with $\Gm$ acting
trivially on the point. This quotient
stack is the classifying stack denoted $B\Gm$. Thus the origin of $\Aquot$
is in fact the classifying stack $B\Gm$ and the restriction $Gr^W(V)$ is an
object over $B\Gm$, which really means an $n$-stack with action of $\Gm$.
Another way to understand this is to think of the total stack as being a
$\Gm$-equivariant stack over $\Aone$; since the origin is a fixed point,
the fiber over the origin is a $\Gm$-equivariant stack.

The homotopy groups of the associated-graded $Gr^W(V)$ are vector spaces
with $\Gm$ action, so they decompose into pieces which one could denote by
$\pi _iGr^W_k(V)$. These are of course the analogues of the classical
graded pieces of the weight filtration. They are provided with Hodge
filtrations, and real structures, so it makes sense to look at the
complex conjugate of the Hodge filtration and ask, for example, that
the Hodge filtration be $h$-opposed with its complex conjugate on 
$\pi _iGr^W_k(V)$. Here $h$ is to be chosen according to Deligne's
arithmetic of \cite{hodge3} which combines the degree of the weight filtration
$k$ and  the homotopical degree $i$: we ask that $h =
k-i$ (cf p. \pageref{footnotepage}).

The above ideas are enough to try for a sketch of definition of the notion
of nonabelian mixed Hodge structure. We define an object which we call a {\em
pre-namhs} (see p. \pageref{pnamhspage}) 
which basically consists of a geometric $n$-stack with a bifiltration
$(V_{DR},F,W)$ (a bifiltration is determined by its total stack $Tot^{W,F}(V)$
over $\Aquot \times \Aquot$ cf p. \pageref{bifiltpage}), plus a real 
stack $V_{B,\rr }$
with a filtration $W$ and an analytic equivalence between $(V_{DR}, W)$ and
the complexified $(V_{B,\cc }, W)$.  The associated graded $Gr^W(V_{DR})$
is a $\Gm$-equivariant geometric $n$-stack, so (after linearization) 
it breaks up into components
which we denote $Gr^W_k(V_{DR})$. These components have Hodge filtrations $F$,
and analytic real structures coming from their analytic identifications with
the complexification of $Gr^W_k(V_{B,\rr})$.  In  particular, the homotopy
group sheaves $\pi _iGr^W_k(V_{DR})$ are actually just complex vector spaces
(for $i\geq 2$), and they 
have Hodge filtrations and real structures. Thus it makes sense to ask for a
purity condition relating the Hodge filtration and its complex conjugate (this
condition also involves $k$ and $i$, see below). In fact, we don't put on any
such condition in the notion of pre-namhs, because it appears useful to
manipulate
a notion of bare object not necessarily satisfying any conditions.

This definition of pre-namhs is modelled in an obvious way on Deligne's
definitions of
mixed Hodge structure and mixed Hodge complex. More generally, if $OBJ$ is 
any $n+1$-stack (with real structure and analytic extension) then one can define
in a similar way a notion of ``pre-mixed Hodge object of type $OBJ$'' 
(cf page \pageref{premhobject}) as
being a collection consisting of a bifiltered object $(V_{DR}, W,F)$ of $OBJ$,
a real filtered object $(V_{B,\rr}, W)$, and an analytic equivalence $\zeta$.
For $OBJ=VECT$ we recover the notion of pre-mixed Hodge structure, and for
$OBJ=PERF$ we recover the notion of pre-mixed Hodge complex. The case of a
pre-namhs is obtained from $OBJ = nGEOM$.

The above notion of pre-namhs has an obvious problem, 
namely that $\pi _1$ is only
an algebraic group and $\pi _0$ is a sheaf of sets. Thus, it isn't clear
what type of purity condition to put here, even though the $\pi _0$ and 
$\pi _1$ are key parts of what we would like to look at here. Nonetheless,
we made some calculations and looked at some examples (in fact, preliminary
versions of the computations we do in Part III) and came to the conclusion
that there was also a more subtle problem even with respect to the $\pi _i$ 
for $i\geq 2$, in the above way of proceeding. In particular, if one tries
to define a nonabelian mixed Hodge structure as being a pre-namhs which
satisfies some purity conditions on the Hodge filtrations on the homotopy 
group sheaves, then there seems to be missing  the ``strictness'' which
one has come to expect in the mixed Hodge situation.

Our next idea was to look at a specific example where we know that there
should be a ``nonabelian mixed Hodge structure''. The best case we could
find was that of $S^2$. On the one hand, this seems to be an interesting
coefficient stack for nonabelian cohomology (cf \cite{secondaryKS}); the
``schematic sphere''  had previously been studied in the stable range in
\cite{BreenEkedahl}. On
the other hand, it has a motivic structure being the homotopy type of
$\pp ^1$ and in fact this ``motive'' plays an important role for example
in Morel's and Voevodsky's theory \cite{Morel}. In particular,
it is clear that whatever the definition of ``nonabelian mixed Hodge
structure''
is, the homotopy type of $\pp ^1$ should carry one of these.  

In doing the calculations, the main thing to recall is the 
``shift'' of the weight filtration by the homotopical degree 
already mentionned above. This shift
first appeared in Deligne's ``Hodge III'' \cite{hodge3}. Basically there
are three numbers involved. The {\em Hodge degree} which is the integer
$h$ such that the Hodge filtration and its complex conjugate are
$h$-opposed (in a more general ``mixed twistor'' context, $h$ would be
the slope of a semistable bundle over the twistor line $\pp ^1$). 
The {\em weight} $w$ is the degree in the weight
filtration, equal to $k$ on the graded piece $Gr^W_k(V)$. 
And the {\em homotopical degree} $i$ indicates that we are
looking at $\pi _i$. In a mixed Hodge complex, the homotopical degree is
minus the cohomological degree, i.e. one should consider $\pi _i =
H^{-i}$. The mixed Hodge condition ${\bf MHC}$ is a linear relationship between
these parameters:
$$
h = w -i.
$$
See our discussion on page \pageref{footnotepage} for a justification of
the sign in this formula.

Now we get back to our calculation for $\pp ^1$. We have 
(in complexified homotopy theory)
$$
\pi _2 (\pp ^1)= H_2(\pp ^1)  = \cc ,
$$
which is of Hodge type $(-1, -1)$. This gives $h(\pi _2)=-2$. The
homotopical degree is of course $i(\pi _2)= 2$, and the formula ${\bf
MHC}$ determines $w(\pi _2)=0$. Note of course that the degree $w$ in the
weight filtration is the degree in the actual nonabelian filtration we
will put on, which is different from the level in the end-result mixed
Hodge structure on $\pi _2$ exactly because of the formula ${\bf MHC}$. 

On the other hand,
$$
\pi _3 (\pp ^1)= \pi _2 (\pp ^1) \otimes _{\cc}\pi _2 (\pp ^1),
$$
the isomorphism in question being given by the Hopf map or Whitehead
product.
This has Hodge type $(-2,-2)$ so $h(\pi _3)=-4$. Again of course $i(\pi
_3)=3$ and our formula gives $w(\pi _3)= -1$. Now, in terms of the
actual unshifted weight filtration, the Whitehead product map 
$$
\pi _2  \otimes \pi _2 \rightarrow \pi _3
$$
goes from a $1$-dimensional space of weight $w=0$ to a space of weight
$w=-1$. In particular, and here is the {\sc key observation}:
\newline
---the Whitehead products vanish on the associated-graded of the actual
(unshifted) weight filtration.

This observation was confirmed by heuristic calculations with
differential graded algebras. 

It led to the {\sc main idea}: that the associated-graded
for the weight filtration in a nonabelian mixed Hodge structure
should be an abelian object. 

This idea fits in nicely with the whole
``quantization'' industry (cf {\tt q-alg} or {\tt math.QA}), 
where the idea of ``deformation quantization'' is to deform an abelian or
classical situation, to a  nonabelian one. The weight filtration
considered as a $\Gm$-equivariant family of objects parametrized by
$\Aone$, may be considered as a deformation of the associated-graded
$Gr^W(V)$ which is the fiber over $0\in \Aone$. It should  be possible
to explore more deeply the relationship between our notion of weight
filtration and the notion of quantization, but we haven't done that.

This parallel raises the interesting question if any of the recent
considerable progress in understanding deformation quantization, as
exemplified by Kontsevich \cite{Kontsevich} for example, might be useful
in understanding the weight filtration in higher homotopy.

The notion of linearization plays an important role in Toen's higher
Tannaka duality \cite{Toen}.

On a technical level, another motivation for our main idea can be found in the
construction of the Whitehead product on page \pageref{whitehead}. 
The shift by the homotopical degree described above conflicts with the fact that
the Whitehead product goes from $\pi _i \otimes \pi _j$ to $\pi _{i+j-1}$.
Without a supplementary condition, this product would be going between 
shifted mixed Hodge structures with shifts differing by $1$. The condition that
the Whitehead products vanish on the associated-graded of the weight filtration
means that they descend the weight filtration by one step; this extra step is
just what is needed to counteract the extra $1$ in the homotopical degree
occuring in the Whitehead product (cf p. \pageref{whitemotivatepage}). 
With this condition, when the mixed Hodge
structures on the $\pi _i$ are shifted back to being actual mixed Hodge
structures, the Whitehead product becomes a morphism of mixed Hodge  structures.
This situation provides a strong 
motivation for asking that the associated-graded object
$Gr ^W(\Vv )$ be an abelian object.

To put our main idea into practice, we augment the definition of ``pre-namhs''
described above, by adding a notion of ``linearization''. If $\Vv$ is 
a pre-namhs, then one obtains a pre-namhs denoted $Gr^W(\Vv )$,
which is split in the sense that the weight filtration comes directly from
an action of $\Gm$. A {\em linearization} of $\Vv$ is the data of a split
pre-mixed Hodge complex $\Cc$ with an equivalence between $Gr^W(\Vv )$
and the Dold-Puppe pre-namhs of $\Cc$. (A technical point is that we require
$\Cc$ to be split by a $\Gm$-action splitting the weight filtration, and we
require the equivalence to respect splittings.) 
Without going into all of the definitions
(for which we refer to pages \pageref{perfectpage}, 
\pageref{hemiperfectpage}, \pageref{dplinpage}, \pageref{weightfilteredpage},
\pageref{mhcpage},
\pageref{lpnamhspage}) we just say here that 
giving $Gr^W(\Vv )$ the structure of being the Dold-Puppe of a complex is
basically the same thing as giving it an ``infinite loop-space structure''
which means a ``linear structure'' in an appropriate sense. One effect of this
structure is that $\pi _iGr_k^W(\Vv )$ becomes a complex vector space
(with Hodge filtration, real structure) for all $i\geq 0$. The other effect is
that the Whitehead products for $Gr^W(\Vv )$ vanish, putting into effect our
main idea. We denote the object
$\Cc $ by $LGr^W(\Vv )$, and if such a structure is provided then we say that
$\Vv$ is a {\em linearized pre-namhs}.

Now the purity condition on a linearized pre-namhs is easy to state: it is
just the condition that the  pre-mixed Hodge complex $LGr^W(\Vv )$ should be an
actual mixed Hodge complex. 

This condition is still not enough. It is insufficient on the level of $\pi _0$.
To understand this part, it will be useful to suppose that the objects we are
talking about are schemes rather than $n$-stacks. The weight filtration is a
$\Gm$-equivariant deformation from $V$ to its associated-graded $Gr^W(V)$; on
the other hand, the linearization amounts to an equivalence between $Gr^W(V)$
and a complex vector space. In particular, if the deformation is flat then $V$
is smooth (at least locally near where the deformation is taking place).
At non-smooth points, we could deform $V$ to, say, its normal cone, but that
would then have to be embedded in a vector space. Thus, in general we expect
that the weight filtration is not flat over $\Aquot$. 

In the case where the weight filtration {\em is} flat, which means that we are
looking locally near some smooth points of $V$, it seems likely that the 
``mixed Hodge complex'' condition described above, is in fact sufficient
to define a good notion of nonabelian mixed Hodge structure. The purity
condition on the level of the associated-graded, extends outward in the
deformation in the same way as it does for mixed Hodge complexes. 

In the non-flat case, we look more carefully at the algebraic geometry of the
situation and end up defining the notion of {\em annihilator ideal}. Since the
present introduction is already too long, we won't go into any further details
here but refer the reader directly to the text. We just mention that we isolate
three conditions to put on the annihilator ideals denoted ${\bf A1}$,
${\bf A2}$, and ${\bf A3}$. The first two concern only the weight filtration
and the third is a compatibility with the split mixed Hodge structure on the
associated-graded part. These conditions amount to limitations on the way the
non-flat weight filtration can constitute a ``jump'' when going from $Gr^W(V)$
to $V$. Condition ${\bf A3}$ complements the purity condition on the
associated-graded $Gr ^W$ and requires that this condition persist after
the jump from $Gr ^W(V)$ down to $V$.

A note of caution about the conditions on the annihilator ideals is in order.
We have come up with these conditions only fairly recently in our research on
this problem. We feel that the three conditions represent the right collection
of conditions, but this is more of an opinion than an established fact. 
We support this opinion by showing how the whole theory works in a specific
computation in Part III, and also give some evidence showing how these
conditions already arise in the context of mixed Hodge complexes (cf p.
\pageref{truncatingmhc}).
However, it remains possible that in the future it might turn out to be better
to modify these conditions in some way (in particular, in order to be able to
obtain the various conjectures that we make). 

After presenting the definition of nonabelian mixed Hodge structure, we turn
to the basic construction where this notion appears, namely the construction of
a nonabelian mixed
Hodge structure $\Hh = \underline{Hom}(X_M, \Vv )$ when the coefficient object 
$\Vv$ is itself a nonabelian mixed Hodge structure. We construct 
$\Hh$ as a linearized pre-namhs, when $\Vv$ is a linearized pre-namhs
(modulo a technical hypothesis that the fundamental group object of
$\Vv$ should be represented by a flat linear group scheme).
Then we
state as Conjecture \ref{nonabcoh}, that if $\Vv$ is a nonabelian mixed Hodge structure
then $\Hh$ should be one too.  We don't prove this conjecture in general;
however in Part III we construct a specific pre-namhs $\Vv$ whose underlying
homotopy type is the complexified $2$-sphere, and we show that for this specific
$\Vv$, the nonabelian cohomology stack $\Hh$ is a nonabelian mixed Hodge
structure. This example is fairly instructive in showing what is going
on. 

In the section presenting our conjectures, we also state several other
desireable things
that should be possible to prove given enough time. We refer the reader there
for the statements.

Given: (1) the length necessary
just to explain the definition of nonabelian mixed Hodge structure; (2)
that any proof of these conjectures would probably take a lot more time to do,
a lot more space to write, and be even more highly
technical; and (3) the illustrative nature of the example which we can 
already give in Part III, we felt that
it would be reasonable to leave these main statements as conjectures,
and to post what we have done up until now. 

So, to conclude this introduction, the current paper is basically just a 
presentation of the idea of the definition of nonabelian mixed Hodge structure,
with a computation of a specific example plus 
some indications of how the future developpement of this idea might go.
Unfortunately, even with this relatively limited goal, the paper is too long. 

{\em Acknowledgements:} We would like to thank A. Beilinson, A. Hirschowitz, 
M. Kontsevich, C. Teleman, and
B. Toen for helpful conversations on topics closely connected to 
what is being done here. The second and third authors would like to thank 
UC Irvine for their hospitality during the preliminary period of research
on this question.

\newpage

\begin{center}
{\Large \bf Part I: Nonabelian weight filtrations}
\end{center}

\begin{center}
{\large \bf Conventions}
\end{center}
\mylabel{conventionspage}

We work over the site $\Zz$ of schemes of finite type over $\cc$ with
the etale topology.
We look at $n$-stacks for $0\leq n < \infty$. 

Throughout the paper, all geometric $n$-stacks (which are $n$-stacks of
$n$-groupoids over the site $\Zz$)---see \cite{geometricN}---
will be assumed to be
{\em very presentable} (see \cite{kobe} \cite{relativeLie});
the notation ``geometric $n$-stack''
is taken to mean ``very presentable geometric $n$-stack''. 
We also recall that the definition of ``geometric'' in \cite{geometricN}
includes the condition of being ``of finite type'' in an appropriate
sense.

Denote by $nGEOM$ the $n+1$-stack of very
presentable geometric $n$-stacks. 
This $n+1$-stack includes all not necessarily invertible morphisms, so
it is a stack of $n+1$-categories but not of $n+1$-groupoids. One can
think of it as being a stack of simplicial categories or ``Segal
categories'' cf \cite{ahcs}.  
For a scheme $Z$,
the $n+1$-category $nGEOM (Z)$ is by definition the $n+1$-category
of very presentable geometric $n$-stacks $T\rightarrow Z$. 

This description also works when the base is a stack. In fact,
throughout the paper we shall use the equivalence between two different
ways of thinking about families of $n$-stacks: if $S$ is an $n$-stack
of groupoids
then the $n+1$-stack of $n$-stacks over $S$, i.e. 
$$
\underline{Hom}(I, nSTACK) \times _{nSTACK}\{ S\}
$$
(the stack of arrows in $nSTACK$ whose target is $S$),
is equivalent to 
\linebreak
$\underline{Hom}(S, nSTACK)$. This is discussed in
\cite{aaspects}; admittedly that discussion is incomplete but 
we don't currently wish to add anything to it. Furthermore, we state
without proof that if $f:S\rightarrow nSTACK$ corresponds to an
$n$-stack $T\rightarrow S$, and if $S$ is geometric, then $T$ is
geometric if and only if the classifying morphism $f$ lands in $nGEOM$.
Thus if $S$ is a geometric $n$-stack then a morphism $S\rightarrow nGEOM$ is the
same thing as a geometric $n$-stack $T$ over $S$, i.e. a morphism of
geometric $n$-stacks $T\rightarrow S$. We shall use this equivalence
tacitly throughout the paper. 

In the case $n=0$, a geometric $0$-stack is just an Artin algebraic space.
If $n=1$ then a geometric $1$-stack is an Artin algebraic $1$-stack.
The ``very presentable'' condition corresponds to asking that the 
stabilizer groups for  points of the $1$-stack be affine algebraic groups.

\begin{center}
{\large \bf  Nonabelian filtrations}
\end{center}
\mylabel{nonabfiltpage}

Fix an $n+1$-stack $OBJ$ which we think of as parametrizing a ``type of
object''. The typical example is when $OBJ = VECT$ is the $1$-stack
of vector bundles. This example leads back to the classical notion of filtered
vector space. The example which we shall use is when $OBJ = nGEOM$, 
the $n+1$-stack of
geometric
$n$-stacks. Another useful example is $OBJ = nPERF$, the $n+1$-stack of
perfect complexes of length $n$ (see below).

Recall that a {\em filtered object of $OBJ$} is a
$\Gm$-equivariant morphism from $\Aone$ to $OBJ$, in other words a
morphism of $n+1$-stacks
from (the $1$-stack quotient thought of as an $n+1$-stack) $\Aone /\Gm$, to
$OBJ$.

For the rest of the paper we establish the notation
$$
\Aquot := \Aone /\Gm 
$$
(the quotient $1$-stack),
and we denote the substacks:
$$
\Gm /\Gm \cong \ast \;\; \mbox{by}\;\; [1] \subset \Aquot ,
$$
and
$$
\{ 0\} /\Gm \cong B\Gm \;\; \mbox{by}\;\; [0] \subset \Aquot .
$$

We obtain the $n+1$-stack of filtered objects in $OBJ$, denoted
$$
F.OBJ := \underline{Hom}(\Aquot , OBJ).
$$
Note here and throughout, that whenever we write a formula of this kind
using $Hom$ or internal $\underline{Hom}$ we make the tacit assumption
that
the range object has been replaced with a fibrant replacement, cf 
\cite{vk} \cite{ahcs}.

The {\em underlying object} of a filtered object, is the preimage of $[1]\in \Aquot$, call
it
$V$. It is a point of $OBJ$, i.e. a morphism $V: \ast \rightarrow OBJ$ or
equivalently
a global section $V\in OBJ (Spec (\cc ))$. If $V$ is a global section
of $OBJ$ then a ``filtration of $V$'' is the specification of a
filtered object in $OBJ$ together with an equivalence between the
underlying object and $V$; we denote such a specification by a letter
such as $W$, and the morphism $\Aquot \rightarrow OBJ$ corresponding
to $W$  is then denoted $Tot ^W(V)$. Introduce also the notation
$$
T^ W(V) := Tot ^W(V) |_{\Aone};
$$
in terms of $n$-stacks over a base this could also be written
$$
T^ W(V) := Tot ^W(V) \times _{\Aquot} \Aone .
$$

We denote by $Gr^W(V)$
the {\em associated-graded object}, which is defined to be the restriction to
$[0]\subset \Aquot$. This is an object with an action of $\Gm$, 
because it is really a morphism $B\Gm \rightarrow OBJ$. It is occasionally more
convenient to think of $Gr^W(V)$ as being the object itself (i.e.  the section $\ast
\rightarrow OBJ$ obtained by restriction to the basepoint $\ast
\rightarrow
B\Gm$), provided
with an action of $\Gm$.  With all of these notations we denote a
filtered object in the usual way by a pair of letters $(V,W)$.

\begin{center}
{\bf Linearization}
\end{center}

\mylabel{linearizationpage}

Suppose now that $\Phi : OBJ'\rightarrow OBJ$ is a functor of
$n+1$-stacks. A {\em $\Phi$-linearized} (or just {\em linearized}) filtered object
is a filtered object $(V,W)$ together with
an object $LGr^W(V): B\Gm \rightarrow OBJ'$ and an equivalence
$$
Gr ^W(V) \cong \Phi (LGr^W(V)).
$$
We denote by $F^{\Phi }. OBJ$ the $n+1$-stack of $\Phi$-linearized filtered objects.
The terminology ``linearization'' is due to the fact that we will use
this with $OBJ'$ being a stack of ``linear objects'' (perfect complexes).

\begin{center}
{\bf Filtered vector spaces}
\end{center}

\mylabel{vectorpage}

Before getting to the nonabelian examples which are basic to this paper,
we start by examining the case of filtered vector spaces. 
The references for this are \cite{Rees}, \cite{naht}, \cite{delignetwistor},
\cite{santacruz}, \cite{Sabbah}.
As well as motivating the nonabelian definitions, this discussion is
necessary for fixing our conventions regarding the indexing of
filtrations. 

Apply the
above general discussion to the $1$-stack $VECT$ of vector bundles: for
a scheme $X$ we put $VECT(X)$ equal to the category of vector bundles
(of finite rank) on
$X$. A filtered object of $VECT$ is thus a morphism $\Aquot \rightarrow
VECT$, i.e. a $\Gm$-equivariant vector bundle over $\Aone$. The fiber
over $[1]$ is a vector space $V$, and the total vector bundle is denoted
$T^W(V)=Tot^W(V) |_{\Aone}$. An element $v\in V$ gives rise to a unique $\Gm$-invariant
section $\tilde{v}=\Gm v$ over $\Gm \subset \Aone$. We say that $v\in W_kV$ if 
the section $t^k\tilde{v}$ extends to a section of $T^W(V)$ over $\Aone$;
i.e. if $\Gm v$ vanishes to order  at least $-k$ as a section of
$T^W(V)$. Here $t$ is the coordinate vanishing at the origin of $\Aone$. 
\mylabel{filtdefpage}
We obtain
sub-vector spaces $W_kV\subset V$ which form an {\em increasing}
filtration of $V$:
$$
\ldots W_k V \subset W_{k+1}V \subset \ldots .
$$
The fact that $T^W(V)$ is a vector bundle over $\Aone$ means that this
filtration is exhaustive (and since $V$ has finite rank, it has only a
finite number of nonzero steps).
Furthermore, the fiber of $T^W(V)$ over $0\in \Aone$ is naturally
identified with the associated-graded of the actual filtration defined above:
$$
T^W(V)|_0 \cong \bigoplus Gr^W_k(V) = \bigoplus W_kV / W_{k-1}V.
$$
The identification between $W_kV/W_{k-1}V$ and a subspace of $T^W(V)|_0$
is provided by the map
$$
v\mapsto t^k \tilde{v}(0).
$$
Now note that $T^W(V)$ is a $\Gm$-equivariant bundle, in particular the
fiber over the origin has an action of $\Gm$. To fix conventions for
this action (which has an impact on the sign of the indexing integers below),
we say that the action is considered as an action on the
total space covering the action $\lambda \cdot t = \lambda t$ on
$\Aone$.
If $u$ is a section of a $\Gm$-equivariant bundle then $tu$ is another
section and we have $\lambda \cdot (tu)= \lambda ^{-1}t (\lambda \cdot u)$
(the sign can be seen by looking at the trivial action on the trivial
bundle with the unit section $u$). 
If $\lambda \in \Gm$ then applying the above to the invariant section 
$\tilde{v}$ we get
$$
\lambda \cdot (t^k\tilde{v})=\lambda ^{-k} t^k \tilde{v}.
$$
In particular, the image of $W_kV/W_{k-1}V$ is the subspace of 
$T^W(V)|_0$ on which $\lambda \in \Gm$ acts by the character $\lambda
^{-k}$. \mylabel{signmotpage}
Note that this will motivate our choice of sign in the notations of Theorem \ref{decomp}
and Corollary \ref{decompfilt} below. 

The above construction gives an equivalence between
$F.VECT (Spec (\cc ))$ and the category of filtered vector spaces.
The functor going in the other direction is the ``Rees
module''
construction, see \cite{Rees}, \cite{naht}, \cite{delignetwistor},
\cite{santacruz}, \cite{Sabbah}.

\begin{center}
{\bf Filtered vector sheaves}
\end{center}

The reader is urged to skip this subsection.
\mylabel{filtvectsheafpage}
Let $VSCH$ denote the $1$-stack of vector schemes, and $VSH$ the
$1$-stack of vector sheaves. Recall that a vector scheme over a base $Y$
is a group scheme which locally looks like the kernel of a map of vector
bundles,
and a vector sheaf  (see \cite{Auslander} \cite{Hirschowitz}
\cite{relativeLie}) is a sheaf of groups on the big site (Zariski or
etale)
which locally looks like the kernel of a map of coherent sheaves.
We have inclusions 
$$
VECT \subset VESCH \subset VESH.
$$
We consequently obtain notions of ``filtered vector schemes'' and
``filtered
vector sheaves''. These notions are different but the distinction is
somewhat subtle, because the underlying objects are the same: a vector
sheaf or vector scheme over $Spec (\cc )$ is always just a finite
dimensional vector space. However, a vector scheme on
$\Aone$ is not necessarily a vector bundle, and in turn a vector sheaf
is not necessarily a vector scheme. Thus for a vector space $V$ we  have
three distinct notions of filtration on $V$, namely a vector bundle
filtration; a vector scheme filtration; and a vector sheaf filtration. 

In all three cases, the definition of the actual filtration of $V$ by
subspaces, using the extension property of sections (i.e. that $v\in W_kV
\Leftrightarrow t^k\tilde{v}$ extends to a section of $T^W(V)$), works 
equally well.
Thus, if $(V,W)$ is a filtered vector sheaf then it induces a filtration
of $V$ by sub-vector spaces $W_kV$. This construction is not an equivalence of
categories, however, and so should be used with care. 

We will further discuss these notions when they are necessary; usually
they are not. Unless
otherwise indicated, a ``filtered vector space'' will always mean an
object of $F.VECT$, i.e. a vector space provided with a filtration in
the stack of vector bundles.

\begin{center}
{\bf Filtered geometric $n$-stacks}
\end{center}

The basic example  
\mylabel{filtgeompage} 
which we use all the time is when $OBJ := nGEOM$ is the $n+1$-stack of
geometric
$n$-stacks of groupoids (for finite $n$). In this case, a filtered
object $(V,W)$ consists of the specification of a total geometric
$n$-stack $Tot^W(V)\rightarrow \Aquot$ together with an equivalence
between the fiber over $[1]$ and $V$; the associated-graded $n$-stack 
$Gr^W(V)$ is
the fiber over $0$ with its $\Gm$-action. 

We will study this situation more closely further on below.

\begin{center}
{\bf Split filtrations}
\end{center}

\mylabel{gradetofilt}

Suppose $U\rightarrow B\Gm$ is a $\Gm$-equivariant geometric $n$-stack
(let $V$ be the fiber $V= U\times _{B\Gm}\ast $; then to be precise $V$
has an action of $\Gm$ and $U$ is the quotient). 
Then we can define the {\em filtration of $V$ associated to this 
action} in the following way. Note that we have a morphism $\Aquot
\rightarrow B\Gm$ whose fiber is $\Aone$. Set
$$
Tot ^W(V):= U \times _{B\Gm}\Aquot .
$$
The fiber of this over $[1]$ is just the fiber of $U$ over $\ast$,
namely it is $V$; and the associated-graded or fiber over the origin, is
$Gr^W(V)=U\rightarrow B\Gm$. 

A filtered object obtained from this construction is called a {\em
split} filtered object. However, when we speak of the $n+1$-category of
split filtered objects, we usually mean that the splitting is included
in the structure, in other words this is just the $n+1$-category of
$\Gm$-equivariant objects, and one is conscious that it has a functor
towards the category of filtered objects.

If we apply this discussion to the case of
filtered vector spaces, we get the standard construction which takes a vector
space with $\Gm$-action to the associated filtered vector space, turning
the grading into a filtration.

\begin{center}
{\bf Bifiltered objects}
\end{center}

\mylabel{bifiltpage}

The next example shows why we wrote up the above definition of filtered
object in a totally general way: we can apply it to $OBJ := F.nGEOM$,
the $n+1$-stack of filtered geometric $n$-stacks. We obtain an
$n+1$-stack
$F.F.nGEOM$ of ``geometric $n$-stacks having two filtrations''. This
phrase is in quotation marks because, as we shall see, the datum of an
object in $F.F.nGEOM$ is different from the data of an object $V$
provided with two different filtrations $W_1$ and $W_2$; the
``compatibility between the filtrations'' amounts to extra structure. 

In fact, let's back up and analyse what it means to have an object in
$F.F.OBJ$.  We have
\begin{eqnarray*}
F.F.OBJ & := & \underline{Hom}(\Aquot , F.OBJ)\\
& = &\underline{Hom}(\Aquot ,\underline{Hom}(\Aquot ,OBJ))\\
&=& \underline{Hom}(\Aquot \times \Aquot , OBJ).
\end{eqnarray*}
Put another way, an object of $F.F.OBJ$ is a $\Gm \times
\Gm$-equivariant morphism $\AAA ^2 \rightarrow OBJ$. In passing note
that there is an action of the symmetric group $S_2$ on $\Aquot \times \Aquot$,
and hence on $F.F.OBJ$; the transposition acts as ``interchanging the
two filtrations''. More generally there is an action of $S_m$ on 
$F.\ldots F.OBJ$. 

Now getting back to the case where $OBJ = nGEOM$, an object of $F.F.nGEOM$ is the
same thing as a geometric $n$-stack $T\rightarrow \Aquot \times \Aquot$.
As before we shall adopt the notation $(V,W,F)$ for an object of
$F.F.nGEOM$,
where $V$ is the underlying object and $W$ and $F$ denote the two
filtrations. The total stack which ``is'' the pair of filtrations,  
shall now be denoted by 
$$
Tot^{W,F}(V) \rightarrow \Aquot \times \Aquot .
$$
Its fiber over $([1], [1])$ is identified with $V$. Note that with our
notations we have
$$
Tot^{W,F}(V) = Tot ^F(Tot^W(V)).
$$
We denote the fiber of $Tot^{W,F}(V)$
over
$0$ in the first factor, by $(Gr ^W(V), F)$ and the fiber over $0$ in
the second factor, by $(Gr^F(V), W)$; these are both filtered objects
with $\Gm$-action, i.e. morphisms $B\Gm \rightarrow F.nGEOM$. The fiber
over $(0,0)$ is denoted $Gr^{W,F}(V)$, and we have
$$
Gr^{W,F}(V) \cong Gr^F(Gr^W(V)) \cong Gr^W(Gr^F(V)).
$$
It is an object with $\Gm \times \Gm$-action, i.e. a morphism
$B\Gm \times B\Gm \rightarrow nGEOM$.
We shall sometimes call an object $(V,W,F)$ a {\em bifiltered geometric $n$-stack}.

We can see the difference between the above notion of
bifiltered object, and the notion of the same object having two
different filtrations. If $(V,G)$ and $(V,H)$ are two filtered objects
having the same underlying object, then we have 
$$
Tot^G(V)\rightarrow \Aquot ,
$$
and
$$
Tot^H(V)\rightarrow \Aquot .
$$
These agree over $[1]\in \Aquot$ so we can glue them together to get
a geometric $n$-stack 
$$
T\rightarrow \Aquot  \cup ^{[1]} \Aquot .
$$
Note  that
$$
\Aquot  \cup ^{[1]} \Aquot  = (\AAA ^2 - \{ 0\} ) / \Gm \times
\Gm .
$$
Thus $T$  is a $\Gm \times \Gm$-equivariant object over $(\AAA ^2 - \{ 0\}
)$.
On the other hand we have seen that a bifiltered object means a 
$\Gm \times \Gm$-equivariant object over $\AAA ^2$. Thus, the extra
information needed to go from an object with two separate filtrations,
to a bifiltered object, is an extension of the $n$-stack $T \rightarrow 
(\AAA ^2 - \{ 0\})$ to an
$n$-stack
$T'$ over $\AAA ^2$. 

In the case of bifiltered vector spaces, which is the case
$OBJ = VECT$, the extension in question is automatic; this corresponds to the 
``lemme de Zassenhaus'' of (Deligne \cite{hodge2}, 1.2.1),
which says that the double associated-graded taken in one order is equal to
the double associated-graded taken in the other order.
In our language, an object in $OBJ = VECT$ with two
different filtrations is a $\Gm \times \Gm$-equivariant vector bundle
over $(\AAA ^2 - \{ 0\})$, whereas an object with two filtrations which
are compatible in the sense that they satisfy the ``lemme de Zassenhaus'' 
is a $\Gm \times \Gm$-equivariant vector bundle
over $\AAA ^2$. Recall that 
every algebraic vector bundle over $(\AAA ^2 - \{ 0\})$
extends uniquely to a vector bundle over $\AAA ^2$ (this
is proved by extending the algebraic vector bundle to a reflexive
sheaf by taking the direct image from the open set to all of $\AAA ^2$, and
then using the theorem that reflexive sheaves on smooth surfaces are bundles).
The common double associated-graded is just the fiber over $(0,0)$ of this 
extension. The fact that
one doesn't have a ``lemme de Zassenhaus'' for three filtrations 
translates the
fact that the extension property for vector bundles doesn't  
hold for $(\AAA ^3 - \{ 0 \})$.

It would be interesting to be able to say something about the problem of
existence and/or uniqueness
extensions of families of other types of objects
(schemes, geometric $n$-stacks, perfect complexes) from $(\AAA ^2 - \{ 0\})$
to $\AAA ^2$. While the ``lemme de Zassenhaus'' probably fails in all of
these contexts, it would be good to understand just to what extent it fails.

\begin{center}
{\bf Bifiltrations and linearization}
\end{center}

We finish with a word about bifiltrations where one of the filtrations
is linearized.
\mylabel{bifiltstruct}
In practice below we will work
with ``weight filtrations'' which are filtered objects $(V,W)$ provided
with linearization via a certain functor
$\Phi$ as envisaged above. We will then want to look at filtered objects
in $F^{\Phi}.nGEOM$ (the second filtration being the ``Hodge filtration'').

We will still denote objects of $F.F^{\Phi}.nGEOM$
as
triples $(V,W,F)$, where $(V,W)$ is an linearized filtered object
and $F$ is a filtration of $(V,W)$ considered as an object of
$F^{\Phi}.nGEOM$. Concretely (and denoting for the time being the 
linearization functor by $\Phi : OBJ' \rightarrow nGEOM$)
this means that the graded object $(Gr^W(V),
F)$ is linearized along $\Phi$, i.e. it is the  image
of a filtered object $(LGr^W(V),
F)$
of $OBJ'$ (and also it has the $\Gm$-action because of
being the associated-graded for $W$) or in other words,
$$
(LGr^W(V),
F) : B\Gm \rightarrow F.OBJ' (Spec (\cc )),
$$
with given an equivalence
$$
\Phi (LGr^W(V) , F) \cong (Gr^W(V), F)
$$
of $\Gm$-equivariant objects in $F.nGEOM$.

\begin{center}
{\large \bf  (Hemi-)perfect complexes and Dold-Puppe linearization}
\end{center}
\mylabel{doldpuppepage}

In this section we look at the specific linearization functor $\Phi$ which will
interest us. It is the forgetful
functor from geometric $n$-stacks of spectra, to geometric $n$-stacks.
Equivalently, it is the Dold-Puppe functor from perfect $n$-complexes,
to geometric $n$-stacks.  Let $nPERF$ denote the $n+1$-stack of perfect complexes of
length $n$ (supported in the interval $[-n,0]$). 
Recall from \cite{ahcs} that this may be defined in the
following way. 
\mylabel{perfectpage} 
For any scheme $Z\in \Zz$ we put $nPerf(Z)$ equal to the category of
complexes of sheaves 
\footnote{By 
\cite{relativeLie} it doesn't matter whether we work with complexes of
sheaves of $\Oo$-modules, or complexes of sheaves of abelian groups,
since we are working on the big site.}
on $Z$ which are locally quasiisomorphic to
complexes of vector bundles of length $n$ (supported in $[-n,0]$, say). 
Then set
$nPERF (Z)$ equal to the
simplicial (or Segal) category obtained by localizing
$nPerf(Z)$ by dividing out by quasiisomorphisms, in the sense of
Dwyer-Kan \cite{DwyerKan}. This simplicial category is in fact $n+1$-truncated, i.e. it
may be considered as an $n+1$-category. It is proved in \cite{ahcs} that
the resulting $n+1$-prestack is an $n+1$-stack.

The functor 
$$
DP: nPERF \rightarrow nGEOM
$$
is just the Dold-Puppe functor \cite{Illusie} \cite{geometricN}.

Another way of
defining the $n+1$-stack $nPERF$ is as follows. Fix $N>n$. Then 
the $n+N+1$-substack 
$$
(n+N)GEOM^{ptd,N-conn} \subset (n+N)GEOM
$$ 
consisting of pointed, $N$-connected 
geometric (and by our convention, very presentable) $n+N$-stacks, is in
fact $n+1$-truncated. (Note that an object of this substack, over a
scheme $Z$, is an $n+N$-stack $T\rightarrow Z$ provided with a section
$t: Z\rightarrow T$ and such that $T$ is relatively $N$-connected over $Z$.)
This substack is equivalent to $nPERF$. The
equivalence in one direction is obtained by taking a perfect  complex
supported in $[-n,0]$ and shifting it to $[-n-N, -N]$ then applying
Dold-Puppe.
In the other direction, take the singular cochain complex for example. The
fact that these two constructions are inverses is the fact that rational
homotopy theory in the stable range is equivalent to the theory of
complexes of rational vector spaces
\cite[Ch.8\S1,
Theorem~7]{Margolis}. (We leave it to the reader to fill
in the details of the proof of this equivalence).

In this second point of view, the functor  
$$
DP : nPERF = (n+N)GEOM^{ptd,N-conn} \rightarrow nGEOM
$$
is the functor which to an object $(T,t)$ associates the $N$-th iterated
loop stack of $T$ at the basepoint $t$. 

{\bf  Example:} \mylabel{perfectpointpage}
A perfect complex over $Spec (\cc )$ is just a complex of $\cc$-vector
spaces. In this case, it is non-canonically quasiisomorphic to a complex
with zero differential (the direct sum of its cohomology spaces) so its
Dold-Puppe
decomposes as a product of things of the form $K(\cc ^{a_i}, i)$.

\begin{center}
{\bf Hemiperfect complexes}
\end{center}

Somewhat unfortunately, we will sometimes have to deal with a notion slightly
more general than that of ``perfect complex'' which we call 
``hemiperfect complex''. This notion appears, for complexes of
length $1$, in \cite{LMB} for example. The origin of the problem is that
if
$K^{\cdot}$ is a perfect complex (over a scheme $Z$, say) 
supported in the interval $[-n, \infty
)$
then its truncation $C^{\cdot} = \tau _{\leq 0}(K^{\cdot})$ is no longer
a perfect complex. Rather, the $C^i= K^i$ are vector bundles for $i<0$
but
$$
C^0 = \ker ( K^0 \stackrel{d}{\rightarrow} K^1)
$$
is only a vector scheme. 

Thus we introduce the following definition: \mylabel{hemiperfectpage}
a {\em hemiperfect complex of length $n$} over a scheme $Z$ is a  
complex of sheaves of groups
on the big etale site over $Z$, which is locally quasiisomorphic to a
complex of the form $C^{\cdot}$ where $C^i$ is zero outside the interval
$[-n, 0]$, where the $C^{-i}$ are vector bundles for $-n\leq -i \leq -1$,
and where $C^0$ is a vector scheme. Let $nHPerf(Z)$ denote the
$1$-category of such objects, and let $nHPERF (Z)$ denote the Dwyer-Kan
localization dividing out by quasiisomorphisms. We state without proof
that:
\newline
---the Segal category $nHPERF(Z)$ is $n+1$-truncated, so it may be
considered as an $n+1$-category;
\newline
---if $K$ and $L$ are hemiperfect complexes, then 
the morphism $n$-groupoid between $K$ and $L$ in $nHPERF(Z)$ is
calculated by the formula
$$
nHPERF(Z)_{1/}(K,L) = DP (\tau _{\leq 0}({\bf R}Hom (K,L)))
$$
where the ${\bf R}Hom$ is taken in the derived category of complexes of
sheaves of abelian groups or sheaves of $\Oo$-modules on the big Zariski site (we state
that the two answers are the same); and
\newline
---the $n+1$-prestack $Z\mapsto nHPERF(Z)$ is actually an $n+1$-stack on
the site $\Zz$.

These statements are the analogues of statements proved for perfect
complexes in \cite{ahcs}, see also \cite{relativeLie} for the fact that
the ${\bf R}Hom$  may be taken in complexes of sheaves of abelian groups
or sheaves of $\Oo$-modules. We assume that the same proofs work for
hemiperfect complexes.

If $B$ is an $n$-stack then a {\em hemiperfect complex over $Z$} is
defined as being a morphism $B\rightarrow nHPERF$.

It is sometimes useful to have a second point of view
\mylabel{secondpage} on perfect
complexes, for hemiperfect complexes. This is somewhat problematic,
because a simple transcription of that discussion cannot work. Indeed,
an $n+N$-stack which is geometric and also $N$-connected relative to the
base scheme $Z$ is automatically smooth and corresponds to a perfect complex.
In particular, it cannot come
from a hemiperfect complex which is not perfect. Thus we have to 
use the notion of very presentable $n$-stack. 

A hemiperfect
complex $\Cc$ may be viewed as a relatively pointed very presentable $n+N$-stack
$(\tilde{\Cc},p)$ over $Z$ (i.e. $p$ is a section of $\tilde{\Cc} /Z$) such that
$\tilde{\Cc}$ is
relatively
$N$-connected, and such that the $N$-th looping $\Omega ^N(\tilde{\Cc},p)$ is
geometric. If furthermore  $\Omega ^N(\tilde{\Cc},p)$ is {\em smooth} relative
to $Z$, then $\tilde{\Cc}$ itself will be geometric and this will correspond to a
perfect complex. 

Our statement of the equivalence between the notion of truncated perfect
complex and the notion of $(\tilde{\Cc},p)$ as above, is given with
only a sketch of proof: starting with $\Cc$ we set $\tilde{\Cc}$ to be
the
Dold-Puppe of the shift of $\Cc$ by $N$ places to the left. On the other
hand
given $\tilde{\Cc}$ we construct $\Cc$ as the relative homology complex
of $\tilde{\Cc}$ over $Z$, which is then
shifted by $N$ places to the right. The proof that if $\Cc$ is hemiperfect
then $\Omega ^N(\tilde{\Cc},p)$ is geometric, is the same as the proof
for perfect complexes given in \cite{geometricN}. The proof in the other
direction,
that if $\Omega ^N(\tilde{\Cc},p)$ is geometric then $\Cc$ is
hemiperfect, is obtained by looking at $\Cc$ as the relative tangent
complex
of $\Omega ^N(\tilde{\Cc},p)$ over $Z$.

We still obtain a Dold-Puppe functor which we again denote by
$$
DP: nHPERF \rightarrow nGEOM.
$$
It sends $\Cc$ to the geometric $n$-stack $\Omega ^N(\tilde{\Cc},p)$, or
equivalently it may be viewed directly as the Dold-Puppe of $\Cc$ see
\cite{Illusie}.

\begin{center}
{\bf Cohomology and homotopy sheaves}
\end{center}

If $C^{\cdot}$ is a (hemi)perfect complex of length $n$ on a scheme $Z$ then its {\em $-i$-th
cohomology} is a vector sheaf denoted $H^{-i}(C^{\cdot})$ over $Z$, nonzero for
$0\leq i\leq n$ (see \cite{Hirschowitz}). We have the following compatibility with the relative
homotopy sheaves of the Dold-Puppe:
$$
\pi _i(DP(C^{\cdot})/Z) = H^{-i}(C^{\cdot}).
$$
Thus, the $H^{-i}(C^{\cdot})$ plays the role of the homotopy
sheaf $\pi _i$ of the perfect complex. 

In terms of the second point of view given above, if we think of $C^{\cdot}$ as
corresponding to a relatively $N$-connected very presentable $n+N$-stack $\tilde{C}$
over $Z$, then $H^{-i}(C^{\cdot})$ may be defined as being 
equal to the relative homotopy group sheaf
$\pi _{N+i}(\tilde{C}/Z)$. In these terms, the above compatibility becomes
tautological.

\begin{center}
{\bf Dold-Puppe linearization}
\end{center}

\mylabel{dplinpage}

We now say
that a {\em Dold-Puppe-linearized filtered geometric $n$-stack} is an object
$$
(V,W) \in F^{DP}.nGEOM,
$$
i.e. it is a filtered geometric $n$-stack with extra structure on the
associated-graded
via the functor 
$$
DP : nHPERF \rightarrow nGEOM.
$$
Concretely, this means that $(V,W)$ is a quintuple
$$
(V,W) = (V, \epsilon , Tot ^W(V) \rightarrow \Aquot , \eta , 
LGr^W(V))
$$
where $V$ is a geometric $n$-stack; $Tot ^W(V) $
is a geometric $n$-stack over $\Aquot $; $\epsilon$ is an equivalence
$$
\epsilon : V \cong Tot ^W(V) \times _{\Aquot} [1] ;
$$
where 
$$
LGr ^W(V) : B\Gm \rightarrow nHPERF
$$ 
is a hemiperfect complex with
$\Gm$-action;
and where $\eta$ is an equivalence between $DP$ of this latter, and the
fiber of $Tot^W(V)$ over $0$:
$$
\eta :DP (LGr^W(V))\cong Tot ^W(V) \times _{\Aquot} B\Gm .
$$

We denote by $LF.nGEOM:= F^{DP}.nGEOM$, the $n+1$-stack of
Dold-Puppe-linearized
filtered geometric $n$-stacks.

If the object $LGr^W(V)$ is a perfect complex then we say 
\mylabel{perfectlinpage} that the
linearization is a {\em perfect linearization}. In the general case we
may sometimes call it a {\em hemiperfect linearization}  in order to avoid
confusion,
but if no specification is made then the linearization is {\em a priori}
only hemiperfect.

The compatibility between cohomology groups of a complex 
and homotopy groups of the Dold-Puppe, gives in the case of a linearized
filtration
the formulae
$$
\pi _iGr ^W(V) = H^{-i}(LGr^W(V)).
$$
In particular, even $\pi _0Gr^W(V)$ and $\pi _1Gr^W(V)$ are vector
sheaves
(i.e. vector spaces when the base is a point).

\begin{center}
{\large \bf  Further study of filtered $n$-stacks}
\end{center}
\mylabel{furtherstudypage}

If $(V,W)$ is a filtered geometric $n$-stack then 
$$
T^W(V):= Tot^W(V) \times _{\Aquot} \Aone
$$
is a geometric $n$-stack mapping to $\Aone$. This situation leads to a
certain amount of algebraic geometry which is basically the same as the
geometry of a scheme mapping to $\Aone$. In studying this, it is important
to recall \cite{geometricN} that there is a scheme of finite type $Z$ with a surjective
smooth map 
$$
Z\rightarrow T^W(V).
$$

\begin{center}
{\bf Coherent sheaves on geometric $n$-stacks}
\end{center}

Before getting to our geometric study of $T^W(V)$ or $Tot^W(V)$, we
digress for a moment to define the abelian category $Coh(Tot^W(V))$ of
coherent sheaves. Let $COH$ be the $1$-stack (in additive categories)
which to each scheme $Y$ 
assigns the additive category of coherent sheaves on $Y$. If $X$ is a
geometric $n$-stack, put $Coh(X):= \underline{Hom}(X, COH)$. 

\begin{lemma}
\mylabel{coh}
If $X$ is a geometric $n$-stack, then $Coh(X)$ is an abelian category
containing a canonical object $\Oo$. If $f:Z\rightarrow X$ is a flat
morphism of geometric $n$-stacks, then the induced morphism 
$Coh(X)\rightarrow Coh(Z)$ is exact.  If $f$ is flat and surjective then this
morphism
is an exact inclusion of abelian categories.
\end{lemma}
{\em Proof:}
The proof is by induction on $n$. Note that if $X$ is a scheme then 
$Coh(X)=COH(X)$ is the usual abelian category of coherent sheaves, and
smooth  morphisms between schemes induce exact pullback functors. 
The case of schemes takes care of the case $n=-1$  in a certain sense in
the induction, so we don't bother starting with $n=0$. Suppose the
statement is known for $n-1$. Suppose $X$ is a geometric $n$-stack and
choose a smooth surjection from a scheme $Z\rightarrow X$. 
By induction, we know that 
$$
Coh(Z), \;\;\; Coh (Z\times _XZ)\;\; \mbox{and} \;\; Coh(Z\times _XZ\times _XZ)
$$
are abelian categories, and that the pullback functors $p_i$ and
$p_{ij}$
between them are exact embeddings of abelian categories. Now an object
$\Ff$ of $Coh(X)$ is the same thing as a pair $(\Gg , g)$ where 
$\Gg \in Coh(Z)$ and $g: p_1^{\ast}(\Gg )\cong  p_2^{\ast}(\Gg )$
is an isomorphism such that the composition of $p_{12}^{\ast}(g)$ and 
$p_{23}^{\ast}(g)$ is equal to $p_{13}^{\ast}(g)$. From this description
and the fact that the functors $p_i$ and $p_{ij}$ involved are exact, it
is clear that $Coh(X)$ admits kernels, cokernels, and that the coimages
are equal to the images; i.e. $Coh(X)$ is abelian. 
Furthermore it is clear that the pullback
functor
$Coh(X)\rightarrow Coh(Z)$ is a faithful and exact. 

The pair $(\Oo , 1)$ corresponds to the object $\Oo$ of $Coh(X)$.

Suppose now that $f: X\rightarrow X'$ is a flat morphism 
(resp. flat surjection). Let
$Z'\rightarrow X'$ be a smooth surjection from a scheme. Then $X\times
_{X'}Z'\rightarrow X$ is a smooth surjection. Let $Z\rightarrow X\times
_{X'}Z'$ be a smooth surjection from a scheme, so $Z\rightarrow X$ is a
smooth surjection from a scheme. Note that $Z\rightarrow Z'$ is a flat
morphism
(resp. flat surjection) of schemes.
In the diagram 
$$
\begin{array}{ccc}
Coh (X') & \rightarrow & Coh (X) \\
\downarrow && \downarrow \\
Coh (Z') & \rightarrow & Coh (Z)
\end{array}
$$
the vertical arrows are exact and faithful, and the bottom arrow is 
exact (resp. exact and faithful). Therefore the top arrow is exact
(resp. exact and faithful).
\eop

This lemma implies, among other things, that the objects of $Coh (X)$ are
noetherian. 

For any $n$-stack $X$ one can define an additive category 
$$
Coh(X):= \underline{Hom}(X, COH). 
$$
It doesn't look likely that $Coh(X)$ will be
abelian for an arbitrary $n$-stack $X$, but we don't have a
counterexample. Note however that $Coh(X)$ depends only on the
$1$-truncation $\tau _{\leq 1}(X)$ so any $n$-stack whose $1$-truncation
is equivalent to the $1$-truncation of a geometric $n$-stack, will yield
an abelian category $Coh(X)$.

\begin{center}
{\bf Geometric study of $Tot^W(V)$}
\end{center}

We can now continue with our study of $T^W(V)$ and $Tot^W(V)$.
To start with, we say that $(V,W)$ is {\em flat} if $T^W(V)\rightarrow
\Aone$ is a flat  morphism. This means that the morphism of
schemes $Z\rightarrow \Aone$ is flat. 

We obtain from the previous lemma 
an abelian category $Coh(T^W(V))$ of coherent sheaves on the
geometric $n$-stack $T^W(V)$. There is a canonical
object, the structure sheaf $\Oo \in Coh (T^W(V))$. On the other hand,
let $t$ denote the coordinate on $\Aone$. 
This pulls back to a section of
the structure sheaf, which we can view as a morphism in the category
$Coh (T^W(V))$, eventually raised to a power:
$$
t^m: \Oo \rightarrow \Oo .
$$
Let $\ker (t^m)$ denote the subobject of $\Oo$ which is the kernel of $t^m$.
It is a coherent sheaf of ideals (i.e. a subobject of the structure
sheaf) on $T^W(V)$. We have
$$
\ker (t^m) \subset \ker (t^{m+1}).
$$
Since $\Oo$ is noetherian (look over any scheme in $\Zz$ with a smooth surjection
to $T^W(V)$), this increasing sequence of ideals eventually becomes
stationary.

We now note that these ideals are ``$\Gm$-equivariant'' i.e. they 
actually come from ideals (for which we
use the same notation) $\ker (t^m)$ on $Tot ^W(V) = T^W(V) /\Gm$.
To see this, let $\Oo \{ 1\}$ be the structure sheaf $\Oo$ but with
action of $\Gm$ twisted by the standard character of $\Gm$; it is an
object in $Coh (Tot^W(V))$. The coordinate $t$ is a morphism 
$$
t: \Oo \rightarrow \Oo \{ 1\}
$$
and its powers are  morphisms 
$$
t^m : \Oo \rightarrow \Oo \{ m\} .
$$
Thus $\ker (t^m)$ has a meaning on $Tot ^W(V)$.

These ideals are related to the canonical flat
sub-filtration of a filtered object, defined in the following lemma:

\begin{lemma}
\mylabel{flatsub}
Suppose $(V,W)$ is a filtered geometric $n$-stack. Then the
$n+1$-category of objects $(V,W')\rightarrow (V,W)$ (filtrations on $V$ 
mapping to $W$) such that $(V,W')$ is flat, has a final object $(V,
W^{\rm fl})$, which we call the {\em canonical flat subfiltration}.
Furthermore, $Tot^{W^{\rm fl}}(V)$ is the closed geometric substack of
$Tot^W(V)$
defined by the coherent sheaf of ideals $\ker (t^m)$ (for $m\gg 0$).
\end{lemma}
{\em Proof:}
Define $(V, W^{\rm fl})$ by the condition that 
 $Tot^{W^{\rm fl}}(V)$ be the closed geometric substack of
$Tot^W(V)$
defined by the coherent sheaf of ideals $\ker (t^m)$ (for $m\gg 0$).
Note that for $m\gg 0$ this sequence of ideals becomes stationary, as
can be verified by going to a scheme $Z$ with smooth surjection to $Tot^W(V)$. 
Note also that if $Z^{\rm fl}$ is the closed subscheme of $Z$ defined by
the sheaf of ideals $\ker (t^m)$ on $Z$ (which is also the pullback of
$\ker (t^m)$ from $Tot^W(V)$) then $Z^{\rm fl}\rightarrow 
Tot^{W^{\rm fl}}(V)$ is a smooth surjection. 

The coordinate $t$ is not well-defined on $\Aquot$ but on any
scheme mapping to $\Aquot$ we can locally chose a definition of $t$.
On overlaps of neighborhoods of definition, the coordinates $t$ differ
by invertible elements. 
Note that $\ker (t^m)$ on $Z$, for $m\gg 0$, is the ideal of $t$-torsion
in $\Oo (Z)$. In particular, $\Oo (Z^{\rm fl})$ has no $t$-torsion so
the map $Z^{\rm fl}\rightarrow \Aquot$ is flat. 

Suppose that $(V,W')$ is a flat filtration mapping to $(V,W)$, and 
let $Z'$ be a scheme
with a smooth surjection to $Tot^{W'}(V)$, so that $Z' \rightarrow \Aone$ is flat. 
We may assume that there is a lifting to a map $Z'\rightarrow Z$. 
Then $\Oo (Z')$ has no
$t$-torsion, so the $t$-torsion in $\Oo (Z)$ goes to zero in $\Oo (Z')$.
Thus, the map factors through $Z'\rightarrow Z^{\rm fl}$. This implies
that the map 
$$
Tot ^{W'}(V) \rightarrow Tot^W(V)
$$
factors through a map 
$$
Tot ^{W'}(V) \rightarrow Tot^{W^{\rm fl}}(V).
$$
The factorization is unique since 
$$
Tot^{W^{\rm fl}}(V)\subset Tot^W(V)
$$
is a closed substack defined by a sheaf of ideals. This completes the proof.
\eop

Next we note that $Coh (Gr^W(V))$ is equivalent to the category of 
objects in $Coh (Tot ^W(V))$ on which $t$ acts as zero. The operation of
restriction of a coherent sheaf from $Tot ^W(V)$ to $Gr ^W(V)$
corresponds to tensoring by $\Oo _{Gr}:= \Oo /t\Oo \{ -1\}$. If $\Ff$ is
a coherent sheaf on $Tot^W(V)$ then
$$
\Ff |_{Gr^W(V)} = \Ff \otimes _{\Oo} \Oo _{Gr}.
$$
All of these statements can be seen by using Lemma \ref{coh} to 
reduce to the case of a scheme
$Z$ surjecting smoothly to $Tot^W(V)$. 

Finally we
define the {\em annihilator ideals} of $(V,W)$ to be the images
$$
Ann ^W(t^m; V):= im (\ker (t^m) \otimes _{\Oo}\Oo _{Gr} \rightarrow \Oo
_{Gr}).
$$
These are ideals on $Gr^W(V)$, in an increasing sequence,
$$
Ann ^W(t^m; V) \subset Ann ^W(t^{m+1}; V)
$$
and this sequence eventually becomes stationary (because $\Oo
_{Gr}$ is noetherian).

\begin{lemma}
Suppose in the above situation, that for all $m\geq 1$ we have
$$
Ann ^W(t; V) = Ann ^W(t^{m}; V).
$$
Then the $\ker (t^m)$ are all equal, and they are already supported on
$Gr^W(V)$, i.e. we have
$$
\ker (t^m) = Ann ^W(t; V).
$$
\end{lemma}
{\em Proof:}
It suffices to prove this for a morphism of schemes $Z\rightarrow \Aone$
(for the ideals $\ker ( \ldots )$ and $Ann (\ldots )$ defined in the
same way on $Z$). To see that this suffices, let $Z$ be a scheme with a
smooth surjection
$$
Z\rightarrow Tot ^W(V) \times _{\Aquot} \Aone ,
$$
and note also that $Z$ has a smooth surjection to $Tot^W(V)$.
The ideals defined relative to the map $Z\rightarrow \Aone$ are
the pullbacks of the corresponding ideals on $Tot^W(V)$, and since 
$Z\rightarrow Tot^W(V)$ is a smooth surjection, the formulae in question
may be proven after pullback to $Z$.

Now the question for $Z\rightarrow \Aone$ is just a statement in usual
algebraic geometry. To justify it, note that we have an exact sequence
$$
0 \rightarrow \Tt \rightarrow \Oo (Z) \rightarrow \Ff \rightarrow 0
$$
where $\Tt$ is the submodule of $t$-torsion in $\Oo (Z)$, and $\Ff$ is 
flat over a neighborhood of $0\in \Aone$. Note that the torsion is of
finite
length, in other words $\Tt$ is really a $\cc [t]/t^M$-module for some $M$.
Then $\ker (t^m)$ and $Ann (t^m, Z)$ depend only on $\Tt$
(these
objects may be defined for any $\cc [t]$-module, and the statement is that
the
objects for $\Tt$ coincide with those for $\Oo (Z)$).  On the other
hand,
the standard proof of the Jordan normal form works even in this context
where the vector spaces are of infinite dimension but the order of
nilpotency is bounded, so we may write
$$
\Tt = \bigoplus M_i, \;\;\;\; M_i = \cc [t]/t^{e_i}.
$$

On a piece $M_i$, direct calculation shows that 
$Ann (t^m; M_i)$ is zero if $m<e_i$ and is $\cc$ if $m\geq e_i$. Thus the
condition that $Ann (t^m; M_i)=Ann (t; M_i)$ implies that $e_i=1$. 
On the pieces $M_i$ with $e_i=1$ (which are just of the form $M_i=\cc$)
one immediately gets that $\ker (t^m)$ is independent of $m$ and that
these are equal to $Ann (t, M_i)$. 

In this proof, the condition of the lemma 
(which we shall call ${\bf A1}$ below)
 is seen as equivalent to the statement that
the $t$-torsion in $\Oo (Z)$ is actually a  $\cc =\cc [t]/(t)$-module. 
\eop

\begin{center}
{\bf Conditions ${\bf A1}$ and ${\bf A2}$}
\end{center}

We now state the {\em first and second conditions on the annihilator ideals} that
will be used (for the weight filtration) 
in the definition of nonabelian mixed Hodge structure:
\mylabel{a1a2}
\newline
${\bf A1}(V,W)$\, This condition says that for all $m\geq 1$ we have
$$
Ann ^W(t; V) = Ann ^W(t^{m}; V).
$$
This says that the conclusion of the previous lemma applies. 

If this condition holds, then $Ann ^W(t;V)$ is an ideal on $Tot^W(V)$, and the quotient by
this ideal is the structure sheaf for the canonical flat extension
$Tot ^{W^{\rm fl}}(V)$. 
\newline
${\bf A2}(V,W)$\, This condition says that $Ann ^W(t,V)\subset
\Oo _{Gr^W(V)}$ is contained in the ideal of
$0\in Gr^W(V)$ (the point $0$ being defined because of the structure of
complex,
giving a structure of vector space on $\pi _0Gr^W(V)$).

\begin{corollary}
If $(V,W)$ satisfies conditions ${\bf A1}(V,W)$ and ${\bf A2}(V,W)$ then
there is a quasi-finite flat cover $Y\rightarrow \Aone$, containing a lift
of the origin, such that
the morphism $Tot^W(V)\times _{\Aquot}Y \rightarrow Y$ admits a section.
\end{corollary}
{\em Proof:}
Let $(V, W^{\rm fl})$ denote the canonical flat subfiltration. Our
conditions ${\bf A1}$ and ${\bf A2}$  imply that for $m\gg 0$ the ideals
$\ker (t^m)$ (which are all equal) restricted to the fiber of the total space over the origin,
are not $(1)$. This ideal defines $Tot^{W^{\rm fl}}(V)$, in particular
$Gr^{W^{\rm fl}}(V)$ is nonempty. To finish the proof of the corollary,
we may now assume that $(V,W)$ is a filtered object with $Tot^W(V)$ flat
over
$\Aquot$ and such that $Gr^W(V)$ is nonempty. Thus 
$$
T^W(V)\rightarrow \Aone
$$
is a flat morphism whose fiber over $0$ is nonempty.
Choose a smooth surjection from a scheme $Z\rightarrow T^W(V)$
so that $Z\rightarrow \Aone$ is a flat morphism of schemes of
finite type.  In this situation
it is automatic that there is $Y\rightarrow \Aone$ as in the statement, such that there
exists a section $Y\rightarrow Z$.
\eop

This corollary says that, if conditions ${\bf A1}$ and ${\bf A2}$ 
hold, then $Gr^W(V)$ is not ``totally detached''
from $V$. 

{\em Counterexample (1):} 
\mylabel{counter1page}
A simple construction gives an example why we need
some condition of type ${\bf A2}$ to insure that $Gr^W(V)$ isn't
detached from $V$. Let $V$ be a scheme, and let $U\rightarrow B\Gm$ be
any scheme with $\Gm$ action. Then set 
$$
Tot ^W(V):= V \coprod U
$$
with its map to $\Aquot$ (sending $V$ to $[1]$ and $U$ to the $B\Gm$
which is the origin). This defines a perfectly good filtration $(V,W)$
in which $Gr ^W(V)=U$. However, $U$ doesn't really have anything to do
with $V$. In this construction, the ideals $\ker (t^m)$ and $Ann ^W(t^m; V)$
are all equal to
the unit ideal $(1)$ of $\Oo _U$. The canonical flat extension $(V,
W^{\rm fl})$ is the filtration whose total space is just $V$ and whose
associated-graded is empty.

{\em Counterexample (2):}
\mylabel{counter2page}
 One cannot conclude, in the above corollary,
existence of a section over $\Aone$ (and certainly not existence of an
equivariant section either) as shown by the following example. Let
$f: Y\rightarrow \Aone$ be a double cover ramified at the origin, and let 
$T:= Y\times \Gm$. Define a morphism $p:T\rightarrow \Aone$ by 
$p(y,t):= tf(y)$. This morphism is $\Gm$-equivariant (by construction)
so it defines a morphism $Y= T/\Gm \rightarrow \Aquot$, in  other words
we get a filtered object $(V,W)$ with $V$ equal to the fiber of $T$ over
$1\in \Aone$ (which is also equal to the inverse image of $\Gm \subset
\Aone$ in $Y$) and with $Tot^W(V):= T/\Gm =Y$. We have
$$
Tot^W(V) \times _{\Aquot}\Aone = T\stackrel{p}{\rightarrow} \Aone .
$$
The fiber of $T$ over the origin has multiplicity two, i.e. the morphism
$p$ is ramified along the fiber; this implies that there doesn't exist a
section over $\Aone$ (or over any etale neighborhood of the origin).

\begin{center}
{\large \bf  Filtered and $\Gm$-equivariant perfect complexes}
\end{center}
\label{gmequivariantpage}

We can apply our discussion of filtered objects, to the stack $nPERF$
of length-$n$ perfect complexes, to obtain an $n+1$-stack $F.nPERF$ of
filtered perfect complexes. The same holds for hemiperfect
complexes,
and we get an $n+1$-stack $F.nHPERF$ of
filtered hemiperfect complexes. 

We leave for another place the task of establishing the relationship
between filtered complexes in the usual sense (cf Illusie \cite{Illusie}) for
example, and objects of $F.nPERF$ or $F.nHPERF$.

\begin{center}
{\bf Delooping and truncation}
\end{center}

Suppose $\Cc$ is a complex of sheaves supported in $[-n, 0]$. Let
$\Omega ^k\Cc$ denote the complex obtained by shifting $\Cc$ to the
right by $k$ places, and doing the ``intelligent'' truncation in degrees
$\leq 0$ which consists of replacing the term of degree zero by the
kernel of its differential to the term of degree one. Thus $\Omega
^k\Cc$ is supported in the interval $[k-n, 0]$. If $k>n$ this produces
the zero complex. We call this complex the {\em $k$-th loop complex of
$\Cc$}. This terminology is due to the following compatibility with the
Dold-Puppe operation:
$$
DP (\Omega ^k \Cc ) \cong \Omega ^k(DP \Cc , 0)
$$
where the terminology on the right means taking $k$ times the loop space
based at the zero section of $DP (\Cc )$. The equivalence is a canonical
weak homotopy equivalence of $n-k$-stacks.

If $\Cc$ is a  hemiperfect complex then $\Omega ^k\Cc$ is again
hemiperfect.
The operation $\Omega$ combines a
shift and a truncation, so 
if $\Cc$ is perfect, the loop
complex might be hemiperfect without being perfect. In fact this is the
main motivation for introducing hemiperfect complexes: they arise as the
truncations of perfect complexes.

\begin{lemma}
\mylabel{strdef}
Suppose $(\Cc , F)$ is a filtered hemiperfect complex. Suppose that the
cohomology sheaves $H^{-i}(Tot ^F(\Cc )/\Aquot )$ are vector bundles
over
$\Aquot$. Then there exists a splitting, namely an equivalence in
$F.nHPERF$
between $(\Cc , F)$ and $\bigoplus _i (H^{-i}(\Cc ), F)$. In particular
$(\Cc , F)$ is a filtered perfect complex, as are all of the loop
complexes
$(\Omega ^k\Cc , F)$. For filtered perfect
complexes the condition (which we denote by ${\bf Str}$) 
that the $H^{-i}$ be vector bundles
over $\Aquot$,  is equivalent to strictness of the
differentials or equivalently degeneration of the spectral sequence at
$E_1$. Finally, strictness of the differentials implies that the cohomology 
vector sheaves are vector bundles.
\end{lemma}
{\em Proof:}
Suppose the cohomology vector sheaves are vector bundles.
The total complex $Tot ^F(\Cc )$ is classified by a collection of extension 
classes in higher $Ext$ groups between the cohomology vector sheaves.
If the vector sheaves are vector bundles then these $Ext$ groups coincide
with the usual ones, and in fact are calculated cohomologically:
$Ext^i(U,V)= H^i (\Aquot , U^{\ast}\otimes V)$. But $\Aquot$ is cohomologically
trivial for vector bundle coefficients (use the principal $\Gm$-bundle
$\Aone\rightarrow \Aquot$ to see this). Thus the higher $Ext$ groups vanish which
implies that the complex splits. 

If the complex splits, it is clear that the spectral sequence degenerates.
This latter condition is equivalent by Deligne \cite{hodge2} Proposition 1.3.2,
to strictness of the differentials; and finally, strictness of the differentials
means that the differentials of the total complex are strict morphisms of vector
bundles, which in turn implies that the cohomology vector 
sheaves are vector bundles.
\eop

If $(\Cc , W)$ is a filtered hemiperfect complex then we can replace
the vector scheme $(Tot^W(\Cc ))^0$ by its flat subscheme 
(cf Lemma \ref{flatsub}
p. \pageref{flatsub})
and obtain the total complex for a filtered perfect complex. We call this
the {\em perfect subfiltration} denoted $W^{\rm fl}\subset W$. 

\begin{center}
{\bf Annihilator ideals for filtered hemiperfect complexes}
\end{center}

In the spirit of the previous remark,
we can also look at the annihilator ideals of a hemiperfect complex.
This discussion constitutes, of course, an abelian version of the
previous
discussion of annihilator ideals for $n$-stacks, and could have been
done before as a motivation for the nonabelian version of the
discussion.
As the notation is actually simpler in the general nonabelian case, we
have not chosen that order of things.

Suppose the base $B$ is $\Aone$ or $\Aquot$. If $C^{\cdot}$ is a
hemiperfect complex, then we define as before the {\em annihilator
ideals} $Ann (t^m; C^{\cdot})$ which are ideals on $H^0(C^{\cdot})$. In
this case they are ideals of linear subspaces (in particular they
are contained in the maximal ideal defining the origin, in other words
condition 
${\bf A2}$ is always satisfied).  To define the annihilator ideals,
$C^0\rightarrow C^{\cdot}$ serves as a smooth surjection from a scheme,
and the annihilator ideals are defined as
descending from the corresponding ideals defined
on the vector scheme $C^0$.
In particular, if
$C^{\cdot}$
were actually a perfect complex (i.e. if $C^0$ were a vector bundle)
then the annihilator ideals would vanish.

It also follows that the annihilator ideals are compatible with the
Dold-Puppe
construction.

If the base is $\Aquot$
so in fact $C^{\cdot}$ corresponds to a filtered complex $(\Cc , W)$
(which we call a ``truncated filtered complex'') then we denote the
annihilator ideals by $Ann ^W(t^m; \Cc )$.

We have the following
relationship between the annihilator ideals and the spectral sequence of a
filtered complex. It generalizes Lemma \ref{strdef} one step further,
to spectral sequences degenerating at $E_2$.

\begin{lemma}
\mylabel{e2}
Suppose $(\Dd , W)$ is a filtered perfect complex of length $n+k$, and let $\Cc = \Omega
^k\Dd $ be the loop complex which is a filtered hemiperfect complex of
length $n$. Suppose that the spectral sequence 
$E_{\cdot}(\Dd )$ degenerates at $E_2$.  Then we have condition ${\bf A1}$
saying that
$$
\forall m \geq 1,\;\;\; Ann ^W(t^m; \Cc )=  Ann ^W(t; \Cc ).
$$
\end{lemma}
{\em Proof:}
Write $\Dd = \{ D^i, \,\, d(i): D^i \rightarrow D^{i+1}\}$ with $D^i$ being
filtered vector spaces, which we can think of as being vector bundles
over $\Aquot$. We obtain complexes $\{ E^i_r, \,\, d_r(i) : E^i_r
\rightarrow E^{i+1}_r\}$ forming the spectral sequence. 

Fix $i$; then the behavior of
the differentials $d_r(i)$ for a given $i$ (i.e. their vanishing or not)
depends only on the original differential between filtered vector spaces
$d(i)$. To see this, note that if we consider the complex $\tilde{D}^{\cdot}$
consisting only of the terms $D^i$ and $D^{i+1}$ then this gives rise to
a spectral sequence $\tilde{E}^{\cdot}_r$ and 
$$
\tilde{E}^i_r \rightarrow E^i_r
$$
is surjective, and 
$$
E^{i+1}_r \rightarrow \tilde{E}^{i+1}_r
$$
is injective, and the differential $\tilde{d}_r(i)$ factors through
$d_r(i)$; thus $d_r(i)=0$ if and only if $\tilde{d}_r(i)=0$. 

We claim that if $d_r(0)=0$ for all $r\geq 2$, then the condition of the
lemma holds for the truncated complex $\Cc$. Note that the condition of
the lemma, for the truncated complex $\Cc$, again depends only on the
differential $d(0): D^0 \rightarrow D^1$ of the original complex. In
view of all of this, we can reduce to consideration of a complex $\Dd$
concentrated only in degrees $0$ and $1$; for the rest of the proof we
make this  hypothesis.

Now, our complex $\Dd$ breaks up into a direct sum of complexes where
either one of the terms is zero, plus a sum of complexes of the form 
$$
d: L\rightarrow M
$$
where $L$ and $M$ are one-dimensional filtered vector spaces and 
$d$ is a nonzero map between them.  To see this decomposition, work with
$$
d(0): D^0 \rightarrow D^1
$$
as a morphism of vector bundles over $\Aquot$ and use the Gauss method
to reduce the matrix of $d(0)$ to a diagonal form (taking note of the
$\Gm$-equivariance, which is preserved if the Gauss operations are done
correctly).

Now in view of this decomposition (and noting that the cases where
$D^0=0$
or $D^1=0$ are trivial), it suffices to treat the case of a complex of
the form $d: L \rightarrow M$ with $L$ and $M$ being one-dimensional.
Let $l$ and $m$ denote the respective weights of $L$ and $M$ (i.e. the
places where the single nonzero steps occur in the filtrations).
Note that $m\leq l$.
In this case, it is easy to see that the condition of degeneration at
$E_2$ is equivalent to the condition $m\geq l-1$. From this we deduce
that the differential, considered as a map between trivial bundles over
$\Aone$, is multiplication either by $1$ or by $t$ (here, pull back from $\Aquot$ to
$\Aone$). From this, we get the desired statement about the annihilator ideals.
\eop

\begin{center}
{\bf A decomposition for $\Gm$-equivariant perfect complexes}
\end{center}

In a Dold-Puppe linearized filtration, the object $LGr^W(V)$ is a
$\Gm$-equivariant hemiperfect complex. In this subsection we point out that
such a thing decomposes into terms $LGr^W_k(V)$ according to the
character of the action
of $\Gm$.

The choice of sign is motivated by the discussion of filtered vector
spaces on p. \pageref{signmotpage}.

\begin{proposition}
\mylabel{decomp}
Let $X$ be any $n$-stack. 
A hemiperfect  complex $C$ over $X\times B\Gm$,
breaks up canonically as a direct sum
$$
C = \bigoplus _k C_k
$$
where the $C_k$ are hemiperfect complexes on $X\times B\Gm$ such that the 
action of $\Gm$ on the cohomology vector sheaves of $C_k$ is via the
character $t\mapsto t^{-k}$.
\end{proposition}
{\em Sketch of proof:}
The proof uses two techniques which we don't develop here (but which
will, one hopes, be developped elsewhere). These techniques are:
tensoring hemiperfect complexes by line bundles; 
and direct images of hemiperfect complexes.
Let $\omega ^k$ be the rank one line bundle on $B\Gm$ corresponding
to the character $t\mapsto t^k$ (and use the same notation for its
pullback to $X\times B\Gm$). Denoting by $p: X\times B\Gm \rightarrow X$
the projection, put
$$
C_k := \omega ^{-k}\otimes p^{\ast} (p_{\ast}(C\otimes \omega ^{k}))
$$
where $p_{\ast}$ is the (full derived) direct image functor
for hemiperfect complexes, which in this case
gives a hemiperfect complex as an answer.
We have a morphism of hemiperfect complexes
$$
i: C \rightarrow \bigoplus  C_k.
$$
Noting that the higher cohomology of $B\Gm$ with vector sheaf
coefficients vanishes, and that the global sections functor just takes
the invariants, we obtain that $i$ is an isomorphism and that the $C_k$
have the desired property.
\eop

If the complex $C$ is perfect then the components $C_k$ are also perfect.
This is because 
the fact that the  higher cohomology of $B\Gm$ vanishes means that 
no truncation operation is necessary to give back complexes 
supported in $[-n,0]$, and
the direct image of a perfect complex remains perfect. 

{\em Remark:}
The construction of the $C_k$ and the isomorphism $i$, are
$n+1$-functorial in $C$ and in $X$. In particular we obtain a morphism
of $n+1$-stacks 
$$
\underline{Hom}(B\Gm , nHPERF) \rightarrow \prod _k nHPERF ,
$$
$$
C\mapsto \left( \ldots , p_{\ast}(C\otimes \omega ^{-k}), \ldots 
\right) 
$$
which is fully faithful with essential image the ``restricted product''
(consisting of those elements which are zero except for a finite number
of places).

Applying the theorem to the case $X=\Aquot$,  we obtain:

\begin{corollary}
\mylabel{decompfilt}
If $(C, F)$ is a $\Gm$-equivariant filtered perfect or 
hemiperfect complex then it 
decomposes canonically as a direct
sum of terms $(C_k, F)$ such that 
the action of $\Gm$  on the $H^i(Tot^F(C_k))$
is via the character $t\mapsto t^{-k}$.  
\end{corollary}

{\em Notation:} If $(V,W,LGr^W(V))$ is a Dold-Puppe linearized
filtration
then we denote by $LGr^W_k(V)$ the component of degree $k$ in the
decomposition given by the previous corollary (recall that 
$LGr^W(V)$ is a $\Gm$-equivariant perfect complex). Similarly, if 
$(V,W)$ is a filtered (hemi)perfect complex then $Gr^W(V)$ is again a 
$\Gm$-equivariant (hemi)perfect complex, and again the previous corollary
gives a decomposition into components denoted $Gr^W_k(V)$.
We call these components (in both cases) the ``$k$-th graded pieces''
of the associated-gradeds.

{\bf Exercise:}
\mylabel{gmdecomppage}
The non-canonical decomposition of p. \pageref{perfectpointpage}
also holds for  $\Gm$-equivariant perfect complexes; in other words,
if $C: B\Gm \rightarrow nPERF$ is a 
 $\Gm$-equivariant perfect complex then it decomposes 
as 
$$
C = \bigoplus _i H^{-i}(C) [i],
$$
non-canonically but still into $\Gm$-equivariant pieces.

\begin{center}
{\large \bf  Weight-filtered $n$-stacks}
\end{center}
\mylabel{weightfilteredpage}

Before getting to the full notion of nonabelian mixed Hodge structure,
we will first discuss just that part which concerns the weight filtration.
Note that this is what is new in the present paper: the Hodge filtration has been
the subject of several previous papers, and its integration into the picture
is essentially straigthforward. On the other hand, for the weight filtration,
we have introduced an essentially new idea which is that of ``linearization''
cf page \pageref{dplinpage}. To explain
this alone we introduce the definition of a ``weight-filtered $n$-stack''.

{\bf Definition:} A {\em pre-weight-filtered $n$-stack} is a quadruple
$$
\{ (V, W) , \;\; LGr ^W(V), \;\; \zeta \} 
$$
where $(V,W)$ is a filtered geometric $n$-stack (cf page
\pageref{filtgeompage} and page \pageref{furtherstudypage}), where
$LGr^W(V)$ is a 
hemiperfect complex of length $n$ (cf page \pageref{hemiperfectpage}), 
and 
$$
\zeta : DP(LGr^ W(V))\cong Gr^W(V)
$$
is a Dold-Puppe linearization (cf page \pageref{dplinpage}). 

In other words, the $n+1$-stack of pre-weight filtered $n$-stacks is just
$F^{DP}.nGEOM$. 

Recall that we have considered two  conditions
which one can put on a filtration:
\newline
---${\bf A1}(V,W)$:\, that 
the annihilator ideals $Ann^W(t^m; V)$ are all equal for $m\geq 1$; qnd
\newline
---${\bf A2}(V,W)$:\, that the annihilator ideal $Ann ^W(t; V)$ is contained
in the ideal of the zero section of $Gr^W(Tot^F(V))$.

A {\em weight-filtered $n$-stack} is a pre-weight-filtered $n$-stack which
satisfies conditions ${\bf A1}(V,W)$ and ${\bf A2}(V,W)$.

\begin{center}
{\bf An example: deformation to the tangent space}
\end{center}

\mylabel{defts}

An important basic example of a pre-weight-filtered $n$-stack is the
following,
which puts a linearized filtration on any pointed geometric $n$-stack.
This discussion is closely related to Fulton's ``deformation to the normal cone''
see \cite{Fulton} \cite{Gerstenhaber} \cite{Rees}.

Let 
$$
p':Q' \rightarrow \Aone 
$$
be the closed subscheme of ${\bf A}^2$ defined by the equation
$(t-u)(t+u)=0$ with projection $p'$ given by the coordinate $t$.
Let $Q:= Q' /\Gm $ where the action of $\Gm$ on $Q'$ is by homotheties. 
This action covers the standard action on $\Aone$ so $p'$ induces a
projection
$$
p: Q\rightarrow \Aquot .
$$
Let $P' = \Aone \subset Q'$ be the subspace defined by $(t-u)=0$ (it
defines a section of $p'$) and let $P=P'/\Gm $. Thus $P$ is the image of
a section of $p$ (i.e. $P=\Aquot$). The projection $p: Q\rightarrow \Aquot$ 
is flat,
projective and the fibers are finite schemes. The fiber over $[1]$ is
a scheme with two points (and $P_{[1]}$ is one of the two points)
whereas the fiber over $[0]$ is the spectrum of the dual numbers $\cc [u
]/u ^2$ with $P_{[0]}$ being the basepoint and the $\Gm$-action being by
homotheties.

If $Y$ is a geometric $n$-stack then 
$$
\underline{Hom}\left(\frac{Q}{\Aquot},\frac{ Y\times \Aquot }{\Aquot}
\right)
$$
is a geometric $n$-stack mapping to $\Aquot$. 
To see this (which is an easy case of (Conjecture 1, \cite{geometricN}, p. 27)
or alternatively a relative version of (Theorem 8.1, {\em op cit})), 
note that the
construction of taking relative internal $\underline{Hom}$ from $Q$ to
an object, is compatible with homotopy fiber products and takes smooth
morphisms to smooth morphisms. In the case of $Y$ being a smooth scheme,
it is clear that the $\underline{Hom}$ stack defined above is again a
smooth
scheme over $\Aquot$. Now any scheme is obtained from smooth ones as a
homotopy fiber product, so we get the desired statement for $Y$ a
scheme;
finally, it follows for any geometric $n$-stack $Y$ by applying the
definition of geometric $n$-stack, which only involves notions of
smoothness, of being a scheme, and of homotopy fiber products. 

The operation of restricting a morphism to the section $P$ provides a
morphism
$$
\underline{Hom}(\frac{Q}{\Aquot},\frac{ Y\times \Aquot }{\Aquot} )
\rightarrow 
\underline{Hom}(\frac{P}{\Aquot},\frac{ Y\times \Aquot }{\Aquot} )
=Y\times \Aquot . 
$$
Suppose now that $y: \ast \rightarrow Y$ is a point. Put
$$
Tot ^{W(y)}(Y):= 
\underline{Hom}\left(\frac{Q}{\Aquot},\frac{ Y\times \Aquot }{\Aquot} \right)
\times _{Y\times \Aquot} \Aquot ,
$$
where the first morphism is the restriction morphism described above and
the second morphism is $(y,1): \Aquot \rightarrow Y\times \Aquot$. This
is a geometric $n$-stack mapping to $\Aquot$. From the previous descriptions of
the fibers of $Q$ over $[1]$ and $[0]$ respectively, we get immediately
that the fiber of $Tot ^{W(y)}(Y)$ over $[1]$ is $Y$ and that the fiber
over 
$[0]$ is the tangent stack $T_y(Y)$ (see \cite{geometricN} for the
definition of the tangent stack).  Thus, by definition $Tot^{W(y)}(Y)$
defines a filtration $W(y)$ of $Y$, whose associated-graded is the
tangent stack to $Y$ at $y$.

Furthermore, recall from \cite{geometricN} that the tangent stack
$T_y(Y)$ actually has a structure of ``spectrum'', i.e. it has a canonical
infinite delooping structure. Since we are in characteristic zero, a
spectrum is the same thing as a complex, thus $T_y(Y)$ has a structure
of perfect complex (one can check in the construction that this
structure is equivariant for the action of $\Gm$ by homotheties on 
$Spec (\cc [u]/u^2)$).  Letting $LGr^{W(y)}(Y)$ denote this perfect  complex
whose Dold-Puppe is $Gr^{W(y)}(Y)$, we obtain a linearization of the
filtration 
$W(y)$, which makes $(Y, W(y), LGr^{W(y)}(Y))$ into a pre-weight-filtered geometric
$n$-stack. 

We obtain an $n+1$-functor
$$
nGEOM ^{ptd} \rightarrow F^{DP}.nGEOM.
$$

{\bf Exercises}: 
\mylabel{defexos}
(1) Observe that condition ${\bf A2}$ is always satisfied
(indeed $Tot ^{W(y)}(Y)$ comes with a section.
\newline
(2) Show that if $Y$ is a scheme with quadratic singularities at $y$ (or if $Y$
is a geometric $n$-stack with a smooth surjection from such a scheme)
then $W(y)$ satisfies condition ${\bf A1}$. The subscheme defined by the
annihilator ideal is the normal cone.

While we don't imagine that all weight filtrations arise in this way, we
do think that this should provide an important construction of certain
types of weight filtrations.  The relationship between weight
filtrations
arising ``in nature'' and those arising from this construction, seems to
be an interesting question for further thought.

\begin{center}
{\bf Interaction with the Hodge filtration}
\end{center}

We apply the discussion from page \pageref{bifiltstruct} to the case of the
weight
filtration. A {\em weight-Hodge linearized bifiltration} of a geometric $n$-stack
is an object of $F.LF.nGEOM = F.F^{DP}.nGEOM$ where $DP : nHPERF
\rightarrow nGEOM$ is the Dold-Puppe linearization functor considered above. 
This object is denoted $(V,W,F)$ where $V$ is a geometric $n$-stack;
where $(V,W)$ is the ``weight-filtration'' on $V$ which is the
one which is linearized, i.e. it is a filtration plus
a ($\Gm$-equivariant) perfect complex $LGr^W(V)$ whose associated 
($\Gm$-equivariant) geometric $n$-stack is $Gr^W(V)$; and finally where
$F$ is a filtration of $(V,W)$ in the stack $F^{DP}.nGEOM$. As
described above, concretely this means that $Tot^F(V,W)$ is a
Dold-Puppe-linearized filtered geometric
$n$-stack over $\Aquot$. In turn this means that there is a
$\Gm$-equivariant filtered
perfect complex 
$$
(LGr ^W(V), F) : B\Gm \rightarrow F.nHPERF
$$
whose underlying filtered geometric $n$-stack (Dold-Puppe) is $(Gr^W(V), F)$. 

In order to clarify what is going on, when we look at both the weight filtration
and the Hodge filtration we will use two different
subscripts for the two copies of $\Aquot$ which come into the picture:
the copy which concerns the weight filtration will be denoted by $\Aquot
_{wt}$ and the copy concerning the Hodge filtration will be denoted by
$\Aquot _{hod}$. Thus a weight-Hodge linearized bifiltration is a morphism of geometric
$n$-stacks
$$
Tot ^{W,F}(V) \rightarrow \Aquot _{hod} \times \Aquot _{wt} ,
$$
together with a perfect complex $Tot^F(LGr^W(V))$ over $[0]_{wt} \times
\Aquot _{hod}$.

In the presence of both weight and Hodge filtrations, we would like to impose
our two conditions ${\bf A1}$ and ${\bf A2}$ on the weight filtration, plus some additional conditions
related to the Hodge filtration. These conditions are quite naturally imposed on
the biggest stack for which they make sense, in the following way.
\newline
---${\bf A1}(V,W,F)$:\, that the annihilator ideals $Ann^W(t^m; Tot^F(V))$
satisfy the condition
${\bf A1}(Tot^F(V), W)$
(they are all equal for $m\geq 1$); 
\newline
---${\bf A2}(V,W,F)$:\, that the annihilator ideal $Ann ^W(t; Tot^F(V))$
satisfies condition ${\bf A2}(Tot^F(V), W)$, namely that it is contained
in the ideal of the zero section of $Gr^W(Tot^F(V))$; and
\newline
---${\bf Fl}(V,W,F)$:\, that the Hodge filtration $F$ is flat, by which we mean that the
morphism corresponding to the Hodge filtration,
$$
Tot^{W,F}(V) \rightarrow \Aquot _{hod} 
$$
is flat.

These three conditions are  meant to insure that the biflitration
$(V,W,F)$ is reasonably well connected to $V$ itself. 
Condition ${\bf A2}$ prevents silly constructions such as in Counterexample
(1) of p. \pageref{counter1page}.
The flatness condition means that the only degeneracies
involved are the non-flatness of the weight filtration.
Note in particular that because of the flatness condition,
$Tot^FGr^W(V)\rightarrow \Aquot _{hod}$ is flat, which implies that
the hemiperfect linearization is actually a perfect one, i.e.
$Tot^FGr^W(V)$ is a perfect complex. 
Condition
${\bf A1}$ puts a limit on how complicated the degeneracy in the 
weight-filtration direction can be.

To bring things back down to earth a bit, we take note of the following
remarks which follow  immediately from the definitions. The homotopy
stacks
$\pi _i(Gr^W(V))$ are vector sheaves over $[0]_{wt} \cong B\Gm$, i.e. they are
finite dimensional complex vetor spaces with $\Gm$-action. This includes
the cases $i=0$ and $i=1$ because they are the cohomology sheaves of the
perfect complex $LGr^W(V)$:
$$
\pi _i(Gr^W(V)) = H^{-i}(LGr^W(V)).
$$
In the case $i=0$ this means that $\pi _0(Gr^W(V))$ is given the extra
structure
of a complex vector space with $\Gm$-action.  Furthermore, a perfect
complex over $B\Gm$ splits up (non-canonically)
as the product of its
cohomology sheaves, so $Gr^W(V)$ splits as a product of its
Postnikov-Eilenberg-MacLane
components (cf the exercise on p. \pageref{gmdecomppage}):
$$
Gr ^W(V) \cong \prod _{i=0}^n K(\pi _i(Gr^W(V))/B\Gm , i).
$$
Now the Hodge filtration induces a filtration on $LGr^W(V)$. 
This means that $(LGr^W(V),F)$ corresponds to a perfect complex
$$
Tot ^F(LGr^W(V)) \rightarrow \Aquot _{hod} .
$$
In particular the cohomology sheaves of this are vector sheaves on
$\Aquot _{hod}$, so they correspond to filtrations of the vector spaces 
$\pi _i(Gr^W(V))$ in the sense that these latter are objects of $F.VESH$ 
(cf p. \pageref{filtvectsheafpage}). This leads us to our next condition, stating that
these vector sheaves are actually vector bundles (and thus that the
filtration is a filtration in the classical sense i.e. an object of 
$F.VECT$). We call this
condition ``strictness'' (see Lemma \ref{strdef}): 
\newline
---${\bf Str}(V,W,F)$:\,  that the cohomology sheaves 
$$
H^{-i}(Tot ^F(LGr^W(V))) 
$$
over $\Aquot _{hod}$, are actually vector bundles.

We can view $(LGr^W(V),F)$ as a filtered complex, and 
by Lemma \ref{strdef} this condition
amounts to asking that it be a {\em strict} filtered complex in the
sense of \cite{hodge2} \cite{hodge3}. 
This condition means, in view of the discussion of $\Gm$-equivariant perfect
complexes, that $Tot^F(LGr^W(V))$ (which is a $\Gm$-equivariant perfect complex
over $\Aquot _{hod}$) decomposes as a direct sum of its cohomology
vector bundles. Recall from Lemma \ref{strdef} that 
it is equivalent to degeneration of the spectral sequence for
the Hodge filtration $F$ on $LGr^W(V)$, at the $E_1$-term. 

In view of the equivalence $Gr^W(V) = \Phi (LGr ^W(V))$ the above
strictness condition can be rewritten as saying that the 
$\pi _i (Tot^F(Gr ^W(V)))$ are vector bundles over $\Aquot$.

\begin{center}
{\large \bf  Analytic and real structures}
\end{center}
\mylabel{analyticrealpage}

In the main definitions to come in Part II, we shall 
need some basic definitions
about analytic objects and  real structures, particularly as they relate
to filtered objects.

\begin{center}
{\bf Algebraic to analytic objects}
\end{center}

Let $\Zz ^{\rm an}$ denote the site of paracompact complex analytic
spaces with the etale topology. 
We have a functor $\Phi : \Zz \rightarrow \Zz ^{\rm an}$
defined by $\Phi (X):= X^{\rm an}$. 
Define the operation $\Ff \mapsto \Ff ^{\rm an}$ from $n$-stacks on
$\Zz$ to $n$-stacks on $\Zz ^{\rm an}$, by the formula
$$
\Ff ^{\rm an}:= \Phi _! (\Ff ),
$$
where we put (from \cite{ahcs}) 
$$
\Phi _{!,pre}(\Ff )(Y):= {\rm colim} _{Y\rightarrow \Phi (X)} \Ff (X),
$$
and we define $\Phi _!(\Ff )$ to be the $n$-stack associated to the
$n$-prestack $\Phi _{!,pre}(\Ff )$.

The $colim$ in the definition or $\Phi _{!,pre}$ 
is of course the homotopy colimit of $n$-stacks
(and since we're talking about $n$-stacks of groupoids, this is the same
as the homotopy colimit of the diagram of spaces \cite{BousfieldKan}).

We have the standard notion of smooth morphism between complex analytic
spaces, and by the construction of \cite{geometricN} this gives rise to
the notion of {\em analytic geometric $n$-stack}. Note however that in the
analytic situation we don't make any analogue of the ``very presentable''
condition.  

The following lemma
gives the fundamental compatibility between the analytic and algebraic notions
of geometricity.

\begin{lemma}
The functor $\Ff \mapsto \Ff ^{\rm an}$ is compatible with homotopy
fiber products. Furthermore, it takes a scheme to its associated
analytic space. As a consequence, it takes geometric $n$-stacks to
analytic geometric $n$-stacks.
\end{lemma}
{\em Proof:}
The $hocolim$ involved in the definition of $\Phi _{!,pre}$
is filtering (because it is taken over a category which admits finite
products and fiber products), therefore it commutes with finite homotopy 
limits. The operation of taking the associated $n$-stack is also compatible
with finite homotopy limits, see \cite{ahcs} \cite{Rezk}, so the operation
$\Ff \mapsto \Ff ^{\rm an}$ is compatible with finite homotopy limits. The
second sentence is standard; and the third sentence is an immediate
corollary
of the fact the definition of ``geometric $n$-stack'' (resp.
``analytic geometric $n$-stack'') refers only to 
homotopy fiber products, schemes (resp. analytic spaces)
and smoothness of maps between schemes (resp. smoothness of maps between
analytic spaces).
\eop

To understand in a simple way what is going on in the above lemma,
suppose $\Ff$ is a geometric $n$-stack with a surjective morphism from a
scheme $Z\rightarrow \Ff$. Let $R= Z\times _{\Ff}Z$, which is a geometric
$(n-1)$-stack. Thus in a certain sense $\Ff$ is presented as the quotient
of $Z$ by the ``relation'' $R$. Now $Z^{\rm an}\rightarrow \Ff ^{\rm an}$
is a smooth surjection and 
$$
R^{\rm an}= Z^{\rm an}\times _{\Ff^{\rm an}}Z^{\rm an}.
$$
Thus $\Ff^{\rm an}$ is presented as the quotient of $Z^{\rm an}$ by
$R^{\rm an}$.

We define filtrations in the same way as in the algebraic case, using
the analytic stack 
$$
\Aquot ^{\rm an}:= (\Aone )^{\rm an} / (\Gm )^{\rm an}.
$$ 
(which is the analytic stack
associated to $\Aquot$). Thus for
example a filtered analytic geometric $n$-stack is a morphism of
analytic geometric $n$-stacks $T\rightarrow \Aquot ^{\rm an}$. The only
construction we need is that if $(V,W)$ is a filtered algebraic
geometric $n$-stack then by taking the analytic $n$-stack $T^{\rm an}$
associated to the total space $T=Tot^W(V)$ we obtain a filtered analytic
geometric $n$-stack $(V^{\rm an}, W)$.

An {\em analytic perfect complex} (resp. {\em analytic hemiperfect complex})
over a base $B$ is a morphism 
$B\rightarrow nPERF^{\rm an}$ (resp. $B\rightarrow nHPERF^{\rm an}$). 
If $B$ is a complex analytic space,
this is the same thing as a complex of
sheaves of $\Oo ^{\rm an}$-modules on $B$ which is locally
quasiisomorphic to a complex of vector bundles supported in $[-n,0]$
(resp. the same but with $C^0$ a vector scheme
\footnote{
We define an analytic vector scheme as being something which is locally
the kernel of a map
of finite rank vector bundles. Now the $1$-stack of analytic vector schemes
is the analytic stack associated to $VESCH$ because any analytic vector scheme
is pulled back from a sort of ``universal'' algebraic vector scheme over the
space of matrices, by an analytic map. One conjectures that this definition
of analytic vector scheme is equivalent to the notion of ``vector space
in the category of complex analytic spaces over the base'' but we don't have
a proof of that.
}). Again,
given an algebraic perfect complex over a base $X$ we obtain an
associated analytic perfect complex over $X^{\rm an}$. 

An important
observation 
\mylabel{observationpage}
in what follows is that an analytic perfect complex over
$\Aquot ^{\rm an}$ has a canonical algebraic structure (because of the
$\Gm$-action on $\Aone$ whose orbit goes out to the only point at
infinity). 
In particular, if $(U_1,W)$ and $(U_2,W)$ are algebraic filtered
perfect complexes, then any morphism between
the associated analytic  filtered complexes, is
actually algebraic. This remark will be  used in defining the notion of
linearization of a pre-namhs: whereas there is an analytic equivalence
involved in the notion of pre-namhs, for the ``linear'' object this
equivalence is taken to be algebraic because no greater generality would
be obtained by using an analytic equivalence.

\begin{center}
{\bf Real structures}
\end{center}

We treat here the extension of ground fields $[\cc : \rr ]$. Most of what we say
should be valid for an arbitrary extension of ground field (although in 
the case of infinite field extensions, some extra work is necessary because
we have to go out of the realm of schemes of finite type). In fact we state 
many things without proof, things which would be more appropriately dealt with in
the general context of an arbitrary field extension.

Let $\Zz _{\rr}$ denote the category of schemes of finite type over $\rr$
with the etale topology. We have a complexification functor
$\Psi : \Zz _{\rr} \rightarrow \Zz $ which sends a real scheme to its
associated complex scheme. On the level of schemes, this functor is defined
by
$$
\Psi (X ) := X\times _{Spec (\rr )} Spec (\cc ).
$$
As in the previous section this gives a functor
$\Psi _!$ from $n$-stacks on $\Zz _{\rr}$ to $n$-stacks on $\Zz$.
We denote this operation by $X\mapsto X_{\cc}$; often an object on
$\Zz _{\rr}$ will be denoted by $X_{\rr}$ in which case the complexified object
is denoted $X_{\cc}$. 

The formula for $\Psi _!$ is
$$
\Psi _! (X)(Y) := {\rm colim}_{Y\rightarrow \Psi (U)} X(U),
$$
where $X$ is a stack on $\Zz _{\rr}$ and $Y\in \Zz $ is a complex scheme.
The colimit is taken over real schemes $U$ with maps $Y\rightarrow \Psi (U)$.
Note in this case however (this is a difference with the analytification functor)
that the category over which the colimit is taken has a final object. Indeed
$\Psi $ acting on schemes has a left adjoint $\Psi ^!$ defined by putting
$\Psi ^!(Y)$ equal to the real scheme obtained by composing 
$Y\rightarrow Spec (\cc )\rightarrow Spec (\rr )$. A morphism of complex schemes
$$
Y\rightarrow \Psi (U) = U\times _{Spec (\rr )}Spec (\cc )
$$
is the same thing as
a morphism of real schemes $\Psi ^!(Y)\rightarrow U$ (take the first projection
of the previous).  Thus the colimit may be seen as a colimit over maps
$\Psi ^!(Y)\rightarrow U$ and we obtain the formula
$$
\Psi _!(X)(Y) := X(\Psi ^!(Y)).
$$

On stacks the functor $\Psi _!$ has a left adjoint again denoted 
$\Psi ^!$ which, heuristically,
takes a complex stack $X\rightarrow Spec (\cc )$ to the real stack obtained by
composing $X\rightarrow Spec (\cc )\rightarrow Spec (\rr )$. (Technically speaking
this definition doesn't make sense because $X$ is a stack over the site $Sch /
\cc $ and $Spec (\rr )$ doesn't exist in that context;  to be correct, 
$\Psi ^!$ must be defined by the adjoint
property.) There is a natural transformation $\Psi ^!(X) \rightarrow \Psi ^!Spec
(\cc )$ and this gives an equivalence of $n+1$-categories between the $n$-stacks
over $\Zz$ and the $n$-stacks over $\Zz _{\rr}$ provided with morphism to
$\Psi ^!Spec (\cc )$. Using this point of view (plus a slight abuse of notation)
it makes sense to write
$$
\Psi _! (X) = X \times _{Spec (\rr )} Spec (\cc ).
$$

If $X$ is an $n$-stack on $\Zz$ then we say that a {\em real structure} on $X$
is the specification of a real $n$-stack $Y_{\rr}$  plus an equivalence
$X\cong Y_{\cc}$.

There is of course a ``Galois'' version of this notion, namely if $X$ is a
complex stack then composing with the complex conjugation involution $\Zz
\rightarrow \Zz$ we obtain a new complex stack denoted $\overline{X}$; and
a real structure on $X$ is the same thing as a twisted $\zz / 2$-action
on $X$ consisting principally of an equivalence $X\cong \overline{X}$ plus some
higher coherence data. We don't wish to go into further detail on this at 
present. Note however that in terms of the ``Galois'' point of view,
$\Psi ^!(X)$ is the scheme $X\coprod \overline{X}$ with real structure
induced by the involution which exchanges the two pieces.

As usual there is a notion of smooth map between real schemes, and using this
notion we obtain following the procedure of \cite{geometricN}, the notion of
{\em real geometric $n$-stack} (which is an $n$-stack on $\Zz _{\rr}$). 
We impose our tacit hypothesis of the introduction, that included in the
condition of being a real geometric $n$-stack is the condition that
the associated complex geometric $n$-stack is also very presentable. 
As the condition is being made on the complexification, it concerns 
the homotopy group sheaves for all choices of
basepoint in the complexified stack.
Nonetheless, one can say that for real basepoints the homotopy group
sheaves of the complex stack are the complexifications of those of the real
stack. In the connected case when there is only one choice of basepoint
up to equivalence, this determines 
the homotopy group sheaves of the complexified stack.

\begin{lemma}
\mylabel{realgeometric}
The operation of complexification of a real $n$-stack is compatible with 
homotopy fiber products. Therefore it sends real geometric $n$-stacks to
geometric $n$-stacks. Furthermore, if $X_{\rr}$ is an $n$-stack on $\Zz _{\rr}$
then $X_{\rr}$ is a real geometric $n$-stack if and only if its complexification
$X_{\cc}$ is a geometric $n$-stack.
\end{lemma}
{\em Proof:} From the formula $\Psi _!(X)(Y)= X(\Psi ^!(Y))$ 
we obtain immediately the 
compatibility with homotopy fiber products. 
This yields compatibility with the notion of geometricity (compatibility with the
very presentable condition is tautological since we impose that condition on 
the complexified stack to begin with).
To prove
the last sentence, suppose $X_{\rr}$ is a real stack and suppose
$Y\rightarrow X_{\cc}$ is a smooth surjection. Then
$$
Y\times _{X_{\cc}}\overline{Y} \rightarrow X_{\cc}
$$
is a smooth surjection which admits a definition over $\rr$
(the involution exchanges the two factors  $Y$ and $\overline{Y}$). 
Using this
we can prove that if $X_{\cc}$ is geometric then so is $X_{\rr}$.
\eop

\begin{center}
{\bf Real structures on moduli stacks}
\end{center}

Let $nGEOM _{\rr}$ denote the stack of real geometric $n$-stacks over $\Zz
_{\rr}$, defined by setting $nGEOM _{\rr}(U)$ equal to the $n+1$-category
of real geometric $n$-stacks $X\rightarrow U$. 

There is also 
a notion of real perfect (resp. hemiperfect)
complex, which just means a perfect (resp. hemiperfect)
complex of real vector bundles on the real scheme. This gives an $n+1$-stack
$nPERF_{\rr}$ (resp. $nHPERF_{\rr}$) 
over $\Zz _{\rr}$. A perfect complex on a real $n$-stack
$X_{\rr}$ is by definition 
a morphism $X_{\rr}\rightarrow nPERF_{\rr}$
(resp. $X_{\rr}\rightarrow nHPERF_{\rr}$). 

\begin{theorem}
We have 
$$
(nGEOM _{\rr})_{\cc} = nGEOM \;\;\mbox{and}\;\;
(nHPERF _{\rr})_{\cc} = nHPERF . 
$$
In other words, the complex moduli stacks $nGEOM$ and $nHPERF$ have canonical
real structures.
\end{theorem}
{\em Proof:}
Suppose $Y$ is a complex scheme. Then $(nGEOM _{\rr})_{\cc}(Y)$ is by definition
the $n+1$-category of real geometric $n$-stacks over $\Psi ^!(Y)$. These are thus
stacks $X$ fitting into a diagram
$$
X\rightarrow \Psi ^!Y \rightarrow \Psi ^!Spec (\cc )\rightarrow Spec (\rr ).
$$
Looking only at the first three terms of this diagram and applying the discussion
above which gives an equivalence between $n$-stacks over $\Zz$ and $n$-stacks
over $\Zz _{\rr}$ provided with morphism to $\Psi ^!Spec (\cc )$, we obtain that
the data of $X$ is the same as the data of an object in $nGEOM (Y)$. This gives
the equivalence in question for $nGEOM$. For $nHPERF$ we can deduce it directly by
considering a hemiperfect complex as an $n+N$-stack which is $N$-connected,
for $N\geq n+1$, whose $\Omega ^N$ is geometric. 
\eop

\begin{center}
{\bf Real structures and filtrations}
\end{center}

Now we turn to filtrations. The stack $\Aquot := \Aone /\Gm$ is defined over $\rr$,
in other words we have a real $n$-stack 
$$
\Aquot _{\rr} := \Aone _{\rr} /{\bf G}_{m,\rr} . 
$$
A {\em real filtration of a real geometric $n$-stack} $V_{\rr}$, denoted
$(V_{\rr}, W)$, is a real geometric
$n$-stack $Tot^W(V_{\rr})$ over $\Aquot _{\rr}$.
The associated complex
$n$-stack is thus a filtration of $V_{\cc}$, we denote this associated
filtration by $(V_{\cc},W)$, with its total space
$$
Tot^W(V_{\cc}):= Tot^W(V_{\rr})_{\cc}.
$$
The associated-graded $Gr^W(V_{\rr})$ is a real geometric $n$-stack with
action of ${\bf G}_{m,\rr}$. Again, this complexifies to $Gr^W(V_{\cc})$.

We take note of the fact that the canonical decomposition  of
$\Gm$-equi\-va\-riant perfect complexes 
(see pages
\pageref{decomp}
or \pageref{decompfilt})
is compatible with the real structure.

A weight-filtered $n$-stack with real structure $(V,W)$ means a collection 
of real stacks just as in the definition of page \pageref{dplinpage}.
In this case, the $\pi _iGr^W(V)\cong H^{-i}(LGr^W(V))$ are real vector
spaces with action of ${\bf G}_{m,\rr}$. The components
$H^{-i}(LGr^W_k(V))$ are thus real subspaces. In the next part of the paper,
this will be combined with a ``Hodge filtration'' which induces a complex
filtration $F$ on the vector space $H^{-i}(LGr^W(V))_{\cc}$, i.e. an actual
filtration by subspaces $F^pH^{-i}(LGr^W(V))_{\cc}$. It thus makes sense to
take the complex conjugate of this filtration with respect to the real structure,
giving a filtration by subspaces $\overline{F}^qH^{-i}(LGr^W(V))$. This allows us
to impose the condition that the Hodge filtration and its complex-conjugate
be $k-i$-opposed. This is the only place where the real structures
which will be involved  in the definitions to follow,
are actually used.

\newpage

\begin{center}
{\Large \bf Part II: The main definitions and conjectures}
\end{center}

We now present the main definitions of pre-nonabelian mixed Hodge structure and
nonabelian mixed Hodge structure, as well as a number of the basic conjectures
which serve, we hope, to explain our motivation for making these definitions.

One small point to notice is that the notion of pre-namhs is not quite analogous
to the notion of pre-weight-filtered $n$-stack which we discussed above. For
pedagogical reasons, we included the notion of linearization in the notion of
pre-weight-filtered $n$-stack. However, for practical reasons the notion of
pre-namhs, without linearization, seems to be the most useful intermediate
version of this definition. 

\begin{center}
{\large \bf  Pre-nonabelian mixed Hodge structures}
\end{center}
\mylabel{pnamhspage}

We can now give our first main definition. A {\em pre-nonabelian mixed Hodge
structure} or {\em pre-namhs}
is a collection
$$
\Vv = \{ (V_{DR}, W,F); \,\, (V_{B,\rr }, W); \,\, \zeta _{\Vv }  \}
$$
where:
\newline
$V_{DR}$ is a geometric (hence by our convention, very presentable)
$n$-stack;
\newline
$(V_{DR}, W,F)$ is a geometric bifiltration of $V_{DR}$ with the filtrations
pa\-ra\-met\-ri\-zed by $\Aquot_{wt}\times \Aquot_{hod}$ (but note that the
weight filtration is
{\em not} linearized);
\newline
$V_{B,\rr}$ is a real geometric $n$-stack;
\newline
$(V_{B,\rr},W)$ is a geometric filtration on $V_{B,\rr}$; and
\newline
$\zeta _{\Vv } $ is an equivalence between the associated analytic
filtered objects
$$
\zeta _{\Vv}: (V_{DR},W)^{\rm an} \cong (V_{B,\cc},W)^{\rm an}.
$$

We obtain an {\em $n+1$-category of pre-namhs} denoted
$nPNAMHS$. Because of the presence of analytic morphisms in the
definition, it is not immediately useful to consider an $n+1$-stack of
these objects so we don't. Constructing $nPNAMHS$ is straightforward
following the definition of a pre-namhs. For example we have the
following diagram of $n+1$-categories
$$
\begin{array}{rcc}
F.F.nGEOM(\ast )  \rightarrow F.nGEOM(\ast )&& \\
& \searrow &\\
&& F.nGEOM^{[{\rm an}]}(\ast )\\
&\nearrow &\\
F.nGEOM_{\rr}(\ast )\rightarrow F.nGEOM (\ast ) &
\end{array}
$$
and $nPNAMHS$ is the limit or homotopy fiber product
$$
nPNAMHS := F.F.nGEOM(\ast )\times _{F.nGEOM^{\rm an}(\ast )}F.nGEOM_{\rr}(\ast
).
$$
Here $F.nGEOM^{[{\rm an}]}$ is the analyic 
moduli $n+1$-stack of filtered analytic
geometric $n$-stacks (which is different from the analytification
of $F.nGEOM$), and 
$F.nGEOM_{\rr}$ is the real moduli $n+1$-stack of filtered real geometric
$n$-stacks. 
The top horizontal morphism is ``forgetting the Hodge filtration'' and
the bottom horizontal morphism is complexification; the diagonal
morphisms are analytification.

Concretely a $1$-morphism 
$$
\varphi : \{ (V_{DR},W,F),\, (V_{B,\rr },W), \, \zeta _{\Vv}\} \rightarrow 
\{ (V'_{DR},W,F),\, (V'_{B,\rr },W), \, \zeta _{\Vv '} \}
$$
in $nPNAMHS$
consists of a morphism of bifiltered objects,
$$
\varphi _{Hod}:(V_{DR},W,F)\rightarrow (V'_{DR},W,F),
$$
a morphism of filtered real $n$-stacks
$$
\varphi _{B,\rr }:(V_{B,\rr },W)\rightarrow (V'_{B,\rr },W),
$$
and a homotopy or $2$-morphism (invertible up to equivalence) between
$$
(\varphi _{B,\cc })^{\rm an}\circ \zeta _{\Vv }  \;\; \mbox{and}\;\;
\zeta _{\Vv '}  \circ (\varphi _{DR})^{\rm an}.
$$

The $n+1$-categories entering into the fiber product defining
$nPNAMHS$ admit finite limits, and the morphisms in the fiber product
are compatible with these limits; therefore the fiber product $nPNAMHS$
also admits finite limits. In particular we can use homotopy fiber
products here. This will not be the case with the notion of ``namhs'',
which explains the utility of isolating out a notion of ``pre-namhs''.

\begin{center}
{\bf Pre-mixed Hodge objects}
\end{center}

Our definition of pre-namhs could be generalized in the following way.
\mylabel{premhobject}
Suppose 
$OBJ$ is an $n+1$-stack (not of groupoids though) with real form $OBJ
_{\rr}$
and with a morphism $OBJ ^{\rm an}\rightarrow OBJ ^{[{\rm an }]}$ of
analytic stacks.
Then a
{\em pre-mixed Hodge object in $OBJ$} is a collection
$$
\{ (V_{DR}, W, F), \;\; (V_{B, \rr}, W), \zeta \}
$$
where $V_{DR}: \ast \rightarrow OBJ$ is an object of type $OBJ$;
where $(W,F)$ is a bifiltration of $V_{DR}$ i.e. 
$$
Tot^{W,F}(V_{DR}) : \Aquot_{wt} \times \Aquot_{hod} \rightarrow OBJ,
$$
with fiber over $[1]\times [1]$ equal to $V_{DR}$; where $V_{B, \rr}:
Spec (\rr )\rightarrow OBJ _{\rr}$ is a real object, with filtration
$W$;
and where $\zeta$ is an equivalence in $OBJ^{[{\rm an}]}$ between
$Tot ^W(V_{DR})$ and the complexification of $Tot ^W(V_{B, \rr})$.
We obtain an $n+1$-category $PMH.OBJ$ of pre-mixed Hodge objects in
$OBJ$ (we don't consider an $n+1$-stack of these objects).

A pre-namhs is just a pre-mixed  Hodge object in $nGEOM$. Below we will
see that a pre-mixed Hodge structure is an object of 
$PMH.VECT$; a pre-mixed  Hodge complex is an object of $PMH.nPERF$; and
we shall have occasion to meet objects of $PMH.nHPERF$, $PMH.VESH$, etc.

\begin{center}
{\bf Split pre-namhs's}
\end{center}

Define a {\em split} pre-namhs
to be a $\Gm$-equivariant pre-namhs
such that the weight-filtrations arise from the construction of 
page \pageref{gradetofilt} which
to a $\Gm$-equivariant stack associates a $\Gm$-equivariant filtered
stack. Denote by $splPNAMHS$ the $n+1$-category of these objects. 
The operation of ``taking the associated-graded for the weight
filtration''
provides an $n+1$-functor
$$
Gr^W : nPNAMHS \rightarrow splPNAMHS.
$$

\begin{center}
{\large \bf  Mixed Hodge complexes and linearization}
\end{center}
\mylabel{mhcpage}

Our next step is the notion of ``linearization'' for the weight
filtrations.  For this, we need to go back to the abelian yet higher
cohomological case of ``mixed Hodge complexes''. In order to
establish notational conventions we first discuss the notion of mixed
Hodge structures. A {\em pre-mixed Hodge structure} or {\em
pre-mhs} (by which we always mean, in the current paper, pre-$\rr$-mhs) is an 
object of $PMH.VECT$, in other words it is an uple
$$
\Vv = \{ (V, W,F); \;\; (V_{\rr}, W), \zeta \} 
$$
where $V$ is a complex vector space with two filtrations $F$ and $W$, 
$V_{\rr}$ is a real vector space
with real filtration $W$, and $\zeta$ is an isomorphism of filtered
objects between $(V,W)$ and the complexification of $(V_{\rr}, W)$. In
particular, we can think of $(V,W,F)$ as being an object of $F.F.VECT$
and $(V_{\rr}, W)$ as an object of $F.VECT_{\rr}$. 

If $\Vv$ is a pre-mhs, we obtain the {\em associated-graded} $Gr^W(\Vv
)$ which is a {\em split} pre-mhs, i.e. it splits up as a direct sum 
$$
Gr^W(\Vv ) = \bigoplus Gr^W_k(\Vv )
$$
of pieces on which the $\Gm$-action is by the character $\lambda \mapsto
\lambda ^{-k}$. On these pieces, the weight filtrations are one-step
(concentrated in degree $k$) and we have Hodge filtrations and real
structures: we can denote these pieces by
$$
\{(Gr^W_k(V), F), \;\; Gr^W_k(V)_{\rr}\} .
$$
In particular, we can take the {\em complex conjugate} of the Hodge
filtration.
To do this, use the old-fashioned approach to a filtration and take the
complex conjugate of every subspace. According to our conventions 
(p. \pageref{vectorpage}), $F$
corresponds to an {\em increasing} filtration $F_p Gr^W_k(V)$, but we
can turn it into a {\em decreasing} filtration by the reindexation
$$
F^pGr^W_k(V):= F_{-p}Gr^W_k(V).
$$
The complex conjugate filtration $\overline{F}$ is defined by setting 
$\overline{F}^pGr^W_k(V)$ equal to the complex conjugate of $F^pGr^W_k(V)$ with
respect to the real structure $Gr^W_k(V)_{\rr}$.
Recall from \cite{hodge2} that one says that $F$ and $\overline{F}$ are
{\em $\mu$-opposed on $Gr^W_k(V)$} if 
$$
F^pGr^W_k(V) \cap \overline{F}^qGr^W_k(V) = 0, \;\;\; p+q >\mu ,
$$
and 
$$
Gr^W_k(V)=\bigoplus _{p+q=\mu}F^pGr^W_k(V) \cap \overline{F}^qGr^W_k(V).
$$
Note that in this condition, we are using the upper-reindexed versions
of $F$. 
In terms of the lower indexing, we get that $F$ and $\overline{F}$ are
$\mu$-opposed on $Gr^W_k(V)$ if 
$$
F_pGr^W_k(V) \cap \overline{F}_qGr^W_k(V) = 0, \;\;\; \mu + p+q <0 ,
$$
and 
$$
Gr^W_k(V)=\bigoplus _{p+q+\mu =0}F_pGr^W_k(V) \cap \overline{F}_qGr^W_k(V).
$$

Here is the place where the conventions are fixed! 
\mylabel{shiftconvpage}
If $s$ is
an integer, we say that a pre-mhs $\Vv$ is an {\em $s$-shifted mhs} if
for every $k$, the filtrations $F$ and $\overline{F}$ are
$k+s$-opposed  on $Gr^W_k(V)$. We say that $\Vv$ is a {\em mhs} if it
is a $0$-shifted mhs.

We now get to mixed Hodge complexes.
Fix $n\geq 0$ which will be the length of the interval of support for
the mixed Hodge complexes we are considering (we only consider ones
of finite length).
Define a {\em pre-mixed Hodge complex} or {\em
pre-mhc}
to be an uple
$$
\Cc = \{ (C_{DR}, W, F),\;\; (C_{B,\rr}, W), \;\; \zeta _{\Cc} \}
$$
where $C_{DR}$ is a perfect complex (i.e. a complex of finite
dimensional $\cc$-vector spaces); $W$ and $F$ form a bifiltration 
in the category of hemiperfect complexes (see page
\pageref{hemiperfectpage}),
so on the whole $(C_{DR}, W, F)$ is an object of
$F.F.nHPERF$; similarly $(C_{B,\rr},W)$ is a filtered real hemiperfect
complex,
and $\zeta _{\Cc}$ is an eqivalence between $(C_{B,\cc}, W)$ and
$(C_{DR}, W)$ as objects of $F.nHPERF$. 

The use of hemiperfect complexes rather than perfect ones is a technical
point which the first-time reader can safely ignore. In the definition
of ``mixed Hodge complex'' (i.e. when we take off the prefix ``pre-'') 
the hemiperfect complexes will be assumed to be perfect.

Notice here in particular that $\zeta _{\Cc}$ is an algebraic
equivalence---we don't speak of associated analytic objects, because in fact an
analytic equivalence would automatically be algebraic, see p. \pageref{observationpage}. 

Let $PMHC$ denote the $n+1$-category of pre-mixed Hodge complexes.
This definition is almost equivalent to the object studied by Deligne 
in \cite{hodge3}, the only difference is our use of hemiperfect rather
than perfect complexes, which comes from the fact that we may want to
consider
truncations of Deligne-type pre-mhc's.

We again have a notion of {\em split} pre-mhc, \mylabel{splmhc} which
means a $\Gm$-equivariant pre-mhc such that the $\Gm$-action splits the
weight filtration. Let $splPMHC$ denote the $n+1$-category of these
objects.
As before, there is a functor $\Cc \mapsto Gr^W(\Cc
)$,
$$
Gr^W: PMHC \rightarrow splPMHC.
$$

The canonical decomposition of $\Gm$-equivariant perfect complexes given
by Theorem \ref{decomp}
(compatible with filtrations and real structures) gives a decomposition
of split pre-mhc: if $\Cc$ is a split pre-mhc then it decomposes
canonically as a direct sum of pre-mhc with trivial (but shifted) weight
filtrations $\Cc _k$. If $\Cc$ is any pre-mhc then the associated-graded
$Gr^W(\Cc )$, as a split-pre-mhc, decomposes as a direct sum which we
write
$$
Gr^W(\Cc ) = \bigoplus _k Gr^W_k(\Cc ).
$$
Recall that the index $k$  means the piece on which $\Gm$ acts by the
character $\lambda \mapsto \lambda ^{-k}$.

The Dold-Puppe functor from $nHPERF$ to $nGEOM$
induces in an obvious way functors
$$
DP: PMHC \rightarrow PNAMHS
$$
and
$$
DP : splPMHC \rightarrow splPNAMHS.
$$

\begin{center}
{\bf Linearization of pre-namhs's}
\end{center}
\mylabel{lpnamhspage}

Using the previous notation we can explain what we mean by
``linearization'' of a pre-namhs. A {\em linearized pre-namhs}
is a triple $(\Vv , LGr^W(\Vv ),\varepsilon )$ where $\Vv$ is a pre-namhs and 
$LGr^W(\Vv )$ is a split pre-mhc, and 
$$
\varepsilon : DP(LGr^W(\Vv )) \stackrel{\cong}{\rightarrow} 
Gr^W(\Vv )
$$
is an equivalence in $splPNAMHS$ of split pre-namhs. Concretely this
means that for each weight-filtered object in $\Vv$, we have an linearized
associated-graded object. For example, $(LGr^W(V_{DR}), F)$ is a
$\Gm$-equivariant filtered hemiperfect complex whose Dold-Puppe is equivalent
to $(Gr^W(V_{DR}), F)$  (via
$\varepsilon$), and similarly for the real Betti
objects.

Define $nLPNAMHS$ to be the $n+1$-category of linearized pre-namhs's; 
\mylabel{nlpnamhspage} it
may be defined as a fiber product of $n+1$-categories in an obvious way,
with a diagram that we leave to the reader to draw.

\begin{center}
{\bf The mixed Hodge complex conditions}
\end{center}

We have a notion of {\em mixed Hodge
complex} which (following \cite{hodge3}) is a pre-mhc which satisfies
the following conditions:
\newline
---the hemiperfect complexes $Tot^{W,F}(\Cc _{DR})$ over $\Aquot
_{wt}\times
\Aquot_{hod}$ and $Tot^W(\Cc _{B,\rr })$ over $\Aquot _{wt, \rr}$ are
actually perfect complexes;
\newline
---strictness of the Hodge filtration on $Gr^W (\Cc )$, a condition
which
we denote by ${\bf Str}(Gr^W(\Cc ), F)$, which actually
means that $Tot^F(Gr^W(\Cc ))$ splits 
(non-canonically) as a product of
its
cohomology vector bundles over $\Aquot _{hod}$ (cf Lemma \ref{strdef}); and 
\newline
---the purity condition \mylabel{puritypage} ${\bf MHC} (Gr^W(\Cc ))$ which 
says that
on the $k$-th graded piece for the $\Gm$-action 
$$
H^{-i}(Gr^W_k(\Cc )) \subset H^{-i}(Gr^W(\Cc ))
$$
the Hodge filtration $F$ and its complex conjugate $\overline{F}$
 (with respect to the real structure on the vector space in question),
are
$k -i$-opposed. In other words, this condition says that the
split pre-mhs's $H^{-i}Gr^W(\Cc )$ are $-i$-shifted mhs's. 
This shifting depending on the cohomological degree, is the fundamental
point behind Deligne's definition \cite{hodge3}. See also a recent discussion in
Saito \cite{Saito}. We give a heuristic justification for the sign in
the next subsection below.

If one has the strictness condition, then the 
$H^{-i}Tot^F(Gr^W(\Cc ))$ are vector bundles over $\Aquot_{hod}$, thus
$$
(H^{-i} Gr^W_k(\Cc ),F) \in F.VECT
$$
are filtered vector spaces in the sense of the discussion of page 
\pageref{vectorpage}.

We obtain an $n+1$-category $MHC$ of mixed Hodge complexes.

More generally, we say that a pre-mhc $\Cc$ is an {\em $s$-shifted mhc}
if it satisfies the strictness condition and if 
the $H^{-i}Gr^W(\Cc )$ are $s-i$-shifted mhs's.

\begin{center}
{\bf The sign in Deligne's shift}
\end{center}

Getting the signs right in the condition ${\bf MHC}$ is a nontrivial
task.
\mylabel{footnotepage} 
Recall from the introduction that there are three numbers involved:
the {\em Hodge degree} which is the integer
$h$ such that the Hodge filtration and its complex conjugate are
$h$-opposed; the {\em weight} $w$ which is the degree in the weight
filtration, equal to $k$ on the graded piece $Gr^W_k(V)$; and 
the {\em homotopical degree} $i$ which indicates that we are
looking at $\pi _i$ or $H^{-i}$. Deligne's shift in the condition 
${\bf MHC}$ is the linear relationship
$$
h=w-i.
$$
We give a heuristic argument here which shows how to determine the sign
in this formula.

We know of course
that the sign of $h$ should be the same as that of $w$; the question is
whether to put $+$ or $-$ in front of $i$. In order to resolve this,
note that for a coefficient stack $V$ (such as an Eilenberg-MacLane stack
$V=K(U,m)$) we have that
$$
\pi _0 \underline{Hom}(X_M, V)
$$
contains (in a vague way) $H^j(X_M, \pi _j(V))$.  Thus we would like to
require that
$$
h(\pi _0 \underline{Hom}(X_M, V) ) =h(H^j(X_M, \pi _j(V)) )
$$
and 
$$
w(\pi _0 \underline{Hom}(X_M, V)) = w(H^j(X_M, \pi _j(V))).
$$
On the other hand, we will be proposing that $X_M$ be constant in the
weight-filtration direction, so 
$$
w(H^j(X_M, \pi _j(V))) = w(\pi _j(V)).
$$
Taking the $j$-th cohomology of $X_M$ has the effect
of adding $j$ to the Hodge degree:
$$
h(H^j(X_M, \pi _j(V))) = h(\pi _j(V)) +j.
$$
Finally, note that $i=0$ for
the $\pi _0$, so whatever the sign is, we will have
$$
h(\pi _0 \underline{Hom}(X_M, V) ) = w(\pi _0 \underline{Hom}(X_M, V)).
$$
Putting all of these together we obtain the formula:
$$
h(\pi _j(V)) +j=w(\pi _j(V))
$$
or
$$
h(\pi _j(V))=w(\pi _j(V))-j.
$$
As $j$ is the homotopical degree of the piece 
$\pi _j$, this justifies the sign in the formula for condition ${\bf MHC}$.

\begin{center}
{\bf Application to linearized pre-namhs's}
\end{center}

If $(\Vv , LGr^W(\Vv ),\varepsilon )$ is an linearized pre-namhs then we
say that it satisfies the {\em strictness condition} ${\bf
Str}(LGr^W(\Vv ), F)$ if the split pre-mhc $LGr^W(\Vv )$ satisfies the
strictness condition (and if the cohomology sheaves are vector bundles,
in particular this condition implies that the hemiperfect complexes
underlying
$LGr^W(\Vv )$ are actually perfect complexes);
and we say that it satisfies the
 {\em mixed Hodge complex condition} denoted 
${\bf MHC}(LGr^W(\Vv ))$ if the split pre-mhc $LGr^W(\Vv )$ satisfies the
``purity condition'' described above. 

If these two conditions are satisfied, then the $H^{-i}(LGr ^W(\Vv ))$
are
split mixed Hodge structures with weight shifted by $i$. In particular, 
for $i=0$ there is no shift and
$H^0(LGr ^W(\Vv ))$ is a split mixed Hodge structure.

\begin{center}
{\bf Truncated mixed Hodge complexes}
\end{center}
\mylabel{truncationmhcpage}

While the Hodge filtration is strict on a mixed Hodge complex (on the
$Gr^W$ this is an axiom, and on the whole mhc it is a result of
\cite{hodge3} Scholie 8.1.9), the weight filtration is not in general
strict.
This is seen by the fact that the spectral sequence for $W$ degenerates
at $E_2$ but not necessarily at $E_1$ (the objects occuring in Proposition
\ref{splitoff2} p. \pageref{splitoff2}) provide
an example of this behavior). A consequence of this is that a
shift and truncation of a mixed Hodge complex is no longer a mixed Hodge
complex. We will see what type of structure it does have, and this
provides motivation for the role of the annihilator ideals in the
definition of namhs. 

The definition of pre-mhc was made using the notion of hemiperfect
complex;
recall that this notion is stable under the truncated shift operation
$\Omega$; thus the notion of pre-mhc is stable under the truncated shift
operation.
The conditions ${\bf Str}$ and ${\bf MHC}$ are also stable under the
truncated shift operation. In particular, it follows that the
associated-graded
$(Gr^W(\Cc ), F)$ is a filtered object in $nPERF$. Thus the only
filtrations that occur in $nHPERF$ rather than in $nPERF$ are the weight
filtrations.

A {\em truncated mixed Hodge complex}  
is the same type of thing as a 
mixed Hodge complex but without the condition that the 
underlying complexes be perfect (thus they are only hemiperfect cf p. 
\pageref{hemiperfectpage}),
and with the following additional requirements. First, 
$Gr^W$ is required to be an actual perfect
complex and the Hodge filtration here is required to be strict. Then we
also require: ${\bf A1}$ that the annihilator ideals are all equal, 
$$
\forall m \geq 1,\;\;\; Ann ^W(t^m; \Cc ) = Ann ^W(t;\Cc );
$$
and ${\bf A3}$ that the linear subspace defined by this annihilator
ideal be a sub-mixed Hodge structure of $H^0(Gr^W(\Cc ))$.

\begin{corollary}
\mylabel{truncatingmhc}
Suppose $\Dd$ is a $k$-shifted mixed Hodge complex of length $n+k$ and let $\Cc
=\Omega ^{k}\Dd$. Then $\Cc$ is a truncated mixed Hodge complex.
\end{corollary}
{\em Proof:}
Recall that the terms 
$$
H^{-i}(Gr ^W_l(\Cc ))=
H^{-k-i}(Gr ^W_l(\Dd ))
$$ 
are pure Hodge structures
with weight $l-i-k+k=l-i$, the first $-i-k$ coming from the shift in the definition
of condition ${\bf MHC}$ and the second $k$ coming from the fact that
we supposed $\Dd$ was $k$-shifted. Thus we get that $\Cc$ satisfies
condition ${\bf MHC}$. The sheaves entering into condition ${\bf Str}$ are 
the same as those for $\Dd$ so $\Cc$ satisfies condition ${\bf Str}$.

Recall from \cite{hodge3} that the spectral sequence of the 
weight filtration on
a mixed Hodge complex such as $\Dd$, degenerates at $E_2$.
Therefore by Lemma \ref{e2} the truncated and shifted complex $\Cc$
satisfies condition ${\bf A1}$. 
For condition ${\bf A3}$ the subspace in question is the 
direct sum over weights
$l$ of the kernel of the
differential
$$
d: H^{-k}(Gr ^W_l(\Dd ))\rightarrow H^{1-k}(Gr ^W_{l-1}(\Dd )).
$$
This differential is a morphism of $l$-shifted pure Hodge structures,
so its kernel is again an $l$-shifted pure Hodge structure. This provides
condition ${\bf A3}$ and completes the proof.
\eop

Given a truncated mixed Hodge complex as defined above, we leave as an
exercise to show that the operation of taking the flat subobject of the
total spaces of the weight filtrations, gives canonically and
functorially a mixed Hodge complex.  Combining this exercise
with the previous corollary, one obtains the result of 
Hain and Zucker \cite{HainZucker} showing how to truncate a mixed Hodge
complex into another mixed Hodge complex.
Certainly the key conditions ${\bf A1}$ and ${\bf A3}$ must appear
in disguised form in the argument of \cite{HainZucker}.

\begin{center}
{\large \bf  Nonabelian mixed Hodge structures}
\end{center}
\mylabel{namhspage}

A {\em nonabelian mixed Hodge structure} or {\em namhs} is a linearized
pre-namhs
(see page \pageref{nlpnamhspage})
which satisfies a number of conditions which we shall now explain.
Before getting to the conditions, we say that the {\em $n+1$-category of
namhs} denoted $nNAMHS$, is defined as the full sub-$n+1$-category
of $nLPNAMHS$ (see page \pageref{nlpnamhspage})
consisting of those objects which satisfy the following
conditions. In particular, a morphism of namhs is defined as a morphism
of linearized pre-namhs. 
\medskip
\newline
{\bf Conditions:} These conditions are divided into parts (I), (II)  and (III).

\begin{center}
(I)
\end{center}

We impose the conditions that we discussed above, namely:
\newline
${\bf A1}(V_{DR},W,F)$ (resp. ${\bf A1}(V_{B,\rr },W)$), 
that the annihilator ideals \newline
$Ann^W(t^m; Tot^F(V_{DR}))$ (resp. $Ann^W(t^m; V_{B,\rr})$)
are all equal for $m\geq 1$;
\newline
${\bf A2}(V_{DR},W,F)$ (resp. ${\bf A2}(V_{B,\rr },W,F)$), 
that the annihilator ideal \newline
$Ann ^W(t; Tot^F(V_{DR}))$ (resp. $Ann ^W(t; V_{B,\rr})$)
is contained
in the ideal of the zero section of
$Gr^W(Tot^F(V_{DR}))$ (resp. $Gr^W(V_{B,\rr})$);
\newline
${\bf Fl}(V_{DR},W,F)$\, that the Hodge filtration $F$ is flat.

\begin{center}
(II)
\end{center}

These are the mixed Hodge complex conditions. First we ask that: 
\newline
${\bf Str}(LGr^W(V_{DR}),F)$ the sheaves 
$$
H^{-i}(Tot ^F(LGr^W(V_{DR}))/\Aquot _{hod})\cong \pi
_i(Tot^FGr^W(V_{DR})/\Aquot ) 
$$
be vector bundles over $\Aquot _{hod}$. In particular, by Lemma \ref{strdef}, the
linearization of $(Tot^F(V_{DR}), W)$ is a perfect one cf p.
\pageref{perfectlinpage}.
\newline
Now the object 
$$
LGr^W(\Vv )=\{ (LGr^W(V_{DR}), F), LGr^W(V_{B,\rr }), LGr^W(\zeta _{\Vv})\}
$$
is basically a graded complex with filtration $F$ and real structure. We
can turn this into a complex with a real structure, a real filtration $W$ and
a complex filtration $F$, by setting $W$ equal to the filtration
associated to the grading, see p. \pageref{gradetofilt}. 
\newline
${\bf MHC}(LGr ^W(\Vv  ))$\, 
We ask that the resulting object be a mixed Hodge complex \cite{hodge3},
see also \cite{Saito}. Note that
this includes the strictness condition ${\bf Str}$ above, plus a
``purity'' condition. By the strictness condition, the cohomology
objects $H^{-i}(LGr^W(\Vv ))$ are vector spaces provided with a Hodge
filtration in the stack $VECT$, which just  means filtered vector spaces
in the usual sense. These vector spaces also have a grading by the
$\Gm$-action and a real structure, and
the purity condition says that on the $k$-th graded piece 
$H^{-i}(LGr^W_k(\Vv ))$ of $H^{-i}(LGr^W(\Vv ))$, the Hodge filtration
$F$ and its complex-conjugate $\overline{F}$ (taken with respect to the
real structure), are $m$-opposed with
$$
m = k-i.
$$
This shift by the homotopical degree $i$ 
originates with Deligne \cite{hodge3}. In our notation, see
the subsection on p. \pageref{footnotepage}
for an explanation of the sign. 

In particular the above condition means that  the vector space 
$$
\Pp := H^0(LGr^W(\Vv ))
$$
is a split mixed Hodge structure. To fix notations for below,
we have
$$
\Pp = \{ (P_{DR}:= H^0LGr^W(V_{DR}) , F), P_{B,\rr }:= 
H^0LGr^W(V_{B,\rr}), \zeta _P\}
$$
together with an action of $\Gm$ on this collection of stuff. 
Here note that $\zeta _P$ is actually an algebraic isomorphism 
of vector spaces between $P_{DR}$ and $P_{B,\cc}=P_{B,\rr}\otimes 
_{\rr}\cc $. Thus $\Pp$ consists of a graded real vector space with a
complex filtration $F$. Turning the grading back into a weight filtration $W$ we
obtain a pre-mixed Hodge structure, and one of the consequences of
conditions ${\bf Str}$ and ${\bf MHC}$ is that this object $\Pp$ is a mixed Hodge
structure.

\begin{center}
(II)
\end{center}

We finally get to our third condition. Recall that $Ann^W(t, \Vv
)$ is a sheaf of ideals on $\pi _0(Gr^W(\Vv ))= \Pp$. The isomorphism
$\zeta _P$ transforms the de Rham version  $Ann^W(t,V_{DR})$ to  
the complexified Betti version $Ann^W(t,V_{B,\cc })$, so we don't need
to distinguish between these. Our sheaf of ideals has a Hodge
filtration,
a real structure, and it is compatible with the $\Gm$-action on $\Pp$.
The $\Gm$-action may be transformed into a weight filtration on $\Pp$.
Thus 
$$
Ann ^W(t; \Vv ) \subset Sym ^{\cdot}(\Pp ^{\ast})
$$
is a sub-pre-mixed Hodge structure (i.e. subobject with two filtrations
and a real structure preserving the weight filtration). The condition of
part (III) is stated as our {\em third condition on the annihilator
ideals:} \mylabel{a3}
\newline
${\bf A3}(\Vv )$, that $Ann ^W(t; \Vv ) \subset Sym ^{\cdot}(\Pp
^{\ast})$
be a sub-mixed Hodge structure. 

This completes our definition of ``nonabelian mixed Hodge structure''.

{\em Remark:} Conditions ${\bf Str}$ and ${\bf MHC}$ are conditions on
$LGr^W(V)$, in other words they are conditions only on the ``abelian'' part;
conditions ${\bf A1}$, ${\bf A2}$, ${\bf A3}$ and ${\bf Fl}$ concern the
nonabelian part.

\begin{center}
{\large \bf Homotopy group sheaves}
\end{center}

\mylabel{htygrpshpage}

A property which we would clearly like to have is that the homotopy
group
sheaves of a nonabelian mixed Hodge structure, should carry mixed Hodge
structures in a natural way. This turns out to be less obvious than one
might
think, and requires a somewhat technical treatment. We will give a
sketch here. The more technical parts of this 
section are optional and are not used in the
computations in Part III. On the other hand, the first set of remarks---including
Proposition \ref{firstversion} and up to the 
section on Whitehead products p. \pageref{whitehead}---should be
illuminating for understanding what is going on in Part III, so we suggest to
read this section lightly the first time.

There is probably a close correspondence between $1$-connected
nonabelian mixed Hodge structures and the type of ``mixed Hodge dga's''
considered by Morgan \cite{Morgan} and Hain \cite{Hain}. In this sense,
everything that we say in the present section should probably be
considered as already known by  \cite{Morgan} and \cite{Hain}; 
and one should consider that our work in this section consists of
transcribing those constructions into the context of our present
definitions of nonabelian mixed Hodge structure using $n$-stacks. 

Throughout this section, $\Vv$ will denote a nonabelian mixed  Hodge
structure which is simply connected, i.e. the $n$-stacks involved are
$1$-connected relative to their base stacks $\Aquot$ or $\Aquot \times
\Aquot$.  In particular, this implies that $\Vv$ is smooth so  
the conditions on the annihilator ideals and the flatness assumption
on the Hodge filtration are automatic and we don't use them.
It is clear that eventually one would like to treat the case of a
non-simply
connected $\Vv$ but this is left for a future paper.

Locally over their bases, all of the $n$-stacks involved in the definition
of $\Vv$ have basepoints, and the basepoints are unique up to 
locally-defined
homotopy which is locally unique  up to a $2$-homotopy which itself is
not unique; all of this is because of our hypothesis that $\Vv$ is
simply connected. In particular there is no need to specify a basepoint
section when speaking of the homotopy group sheaves $\pi _i(\Vv )$.

\begin{center}
{\bf The vector bundle case}
\end{center}

The first situation we consider (and this is the one which is
encountered in our construction of an example in Part III) is the case
when the homotopy group sheaves $\pi _i(\Vv )$ are vector bundles. To
explain this condition recall that $\Vv$ consists essentially of the
data of geometric $n$-stacks 
$$
Tot ^{W,F}(V_{DR})\rightarrow \Aquot
_{wt}\times \Aquot _{hod}, \;\;\;\; Tot ^W(V_{B,\rr })\rightarrow \Aquot
_{wt, \rr}
$$
together with an analytic equivalence and a linearization. We define
$\pi _i(\Vv )$ as being the object consisting of sheaves in the big site
$$
Tot ^{W,F}(\pi _i(\Vv )_{DR}):= \pi _i( Tot ^{W,F}(V_{DR})/  \Aquot
_{wt}\times \Aquot _{hod}) \;\;\; \mbox{over} \;\;\; \Aquot
_{wt}\times \Aquot _{hod}
$$
and
$$
Tot ^{W}(\pi _i(\Vv )_{B,\rr}):=
\pi _i(Tot ^W(V_{B,\rr })/ \Aquot
_{wt, \rr})  \;\;\; \mbox{over} \;\;\; \Aquot _{wt, \rr}
$$
together with an (analytic but which is automatically algebraic)
equivalence between
$Tot ^{W}(\pi _i(\Vv )_{DR})$ and $Tot ^{W}(\pi _i(\Vv )_{B,\cc})$.
In general $\pi _i(\Vv )$ is a pre-mixed Hodge object (p. \pageref{premhobject})
in the stack $VESH$ of vector
sheaves, see below. However,
our
hypothesis in this subsection is that the objects 
$Tot ^{W,F}(\pi _i(\Vv )_{DR})$ and $Tot ^{W}(\pi _i(\Vv )_{B,\rr})$ are vector
bundles over their respective base objects. 

With this hypothesis, $(\pi _i(\Vv )_{DR}, W,F)$ becomes an object of
$F.F.VECT$ (i.e. just a bifiltered vector space) and 
$(\pi _i(\Vv )_{B, \rr}, W)$ becomes an object of $F.VECT_{\rr}$ (i.e. a
filtered real vector space); and the equivalence is just 
$$
(\pi _i(\Vv )_{DR}, W)\cong 
(\pi _i(\Vv )_{B, \cc}, W).
$$
All in all this means that $\pi _i(\Vv )$ is a pre-mixed Hodge structure
in the classical sense, i.e. an object of $PMH.VECT$. 

A first version of our main result of this section is the following proposition.

\begin{proposition}
\mylabel{firstversion}
If $\Vv$ is a simply connected namhs such that the sheaves underlying
the object $\pi _i(\Vv )$ as defined above are vector bundles, then this
pre-mixed Hodge structure is a $-i$-shifted mixed Hodge structure. 
\end{proposition}
{\em Proof:}
Note that
$$
\pi _iGr^W(\Vv ) = Gr ^W (\pi _i(\Vv )),
$$
indeed these two objects are obtained by restricting all of the sheaves
involved, to $[0]\subset \Aquot _{wt}$. Note also that 
$Gr ^W (\pi _i(\Vv ))$ is the actual associated-graded of the pre-mixed
Hodge structure $\pi _i(\Vv )$, in other words it is a vector space with
real structure and Hodge filtration obtained by taking the
associated-graded
of the vector space $\pi _i(\Vv )$ with respect to its weight-fil\-tra\-tion
by subspaces. The indexation of the components determined by the
$\Gm$-action
is compatible with the usual one, i.e.
$$
\pi _iGr^W_k(\Vv ) = Gr ^W_k (\pi _i(\Vv )) = W_k\pi _i(\Vv )/ W_{k-1}\pi _i(\Vv )
$$
(see p. \pageref{signmotpage}).

The hypothesis that $\Vv$ is a namhs means that it is linearized and that
the pre-mhc
$LGr^W(\Vv )$ a mixed Hodge complex. In particular, this means that the
Hodge filtration and its complex conjugate are $k-i$-opposed on
$$
H^{-i}(LGr^W_k(\Vv )) = \pi _i(Gr^W_k(\Vv ))= Gr ^W_k(\pi _i(\Vv )).
$$
But this is exactly the definition of $\pi _i(\Vv )$ being a
$-i$-shifted mixed Hodge structure.
\eop

In the context of these arguments we can briefly say what complicates
the general case:  in general the $\pi _i(Tot ^{W,F}(V_{DR})/\Aquot
_{wt}\times \Aquot _{hod})$ are just vector sheaves rather than vector
bundles over the base $\Aquot
_{wt}\times \Aquot _{hod}$. Thus while they yield actual filtrations of
the underlying vector space $\pi _i(V_{DR})$ by subspaces (cf p.
\pageref{filtdefpage}),
one doesn't recover the total vector sheaves from these filtrations
and in particular the fiber of say $\pi _i(Tot ^W(V_{DR})/\Aquot _{wt})$
over $[0]$ does not in general correspond to the associated-graded for
the filtration by subspaces. The simple argument given for
Proposition \ref{firstversion} no longer works.

\begin{center}
{\bf Normalized homotopy objects}
\end{center}

It is convenient to {\em normalize} the pre-mhs's $\pi _i(\Vv )$ which are
$-i$-shifted mhs's, to turn them into actual mixed Hodge structures. To do this,
re-index the weight filtration shifting the indexation by $i$. Call the resulting
mixed Hodge structures the {\em normalized homotopy groups} of $\Vv$,
denoted $\pi _i ^{\nu} (\Vv )$.

\begin{center}
{\bf Whitehead products}
\end{center}

\mylabel{whitehead}

Before getting on to the treatment of the general problem raised in the
preceding paragraph, we discuss the ``Whitehead products'' in the easy
case represented by Proposition \ref{firstversion}.

If $(X,x)$ is a pointed space then we have the functorial {\em Whitehead product}
which is a bilinear map 
$$
Wh_{i,j}:\pi _i(X,x)\times \pi _j(X,x) \rightarrow \pi _{i+j-1}(X,x).
$$
The same  holds if $X$ is an $n$-groupoid by the equivalence of
\cite{Tamsamani}.
Suppose now that $X$ is an $n$-stack over a base $B$ (which could be a
scheme, a stack or an $n$-stack itself). Suppose $x: B\rightarrow X$
is a basepoint section. Then we obtain a bilinear map
of sheaves over $B$
$$
Wh_{i,j}:\pi _i(X/B,x)\times \pi _j(X/B,x) \rightarrow \pi _{i+j-1}(X/B,x).
$$
If $X$ is relatively $1$-connected over $B$ then locally over $B$ there
is a basepoint section and the homotopy groups don't depend on choice of
basepoint section up to unique isomorphism, so we obtain relative $\pi
_i(X/B)$ which are sheaves  over $B$ (which if necessary means morphisms
from $B$ to the $1$-stack $SH$ of sheaves). 
These exist even though there need not exist a global
section $x: B\rightarrow X$. In
this case the Whitehead product maps again give bilinear maps of sheaves
over $B$,
$$
Wh_{i,j}:\pi _i(X/B)\times \pi _j(X/B) \rightarrow \pi _{i+j-1}(X/B).
$$

Now we return to the case of a $1$-connected namhs $\Vv$. We obtain
bilinear maps of the sheaves occuring in the definitions of the $\pi
_i(\Vv )$, for example
$$
Wh_{i,j}:\pi _i(Tot ^{W, F}(V_{DR})/\Aquot _{wt}\times \Aquot _{hod})
\times 
\pi _j(Tot ^{W, F}(V_{DR})/\Aquot _{wt}\times \Aquot _{hod})
$$
$$
\rightarrow 
\pi _{i+j-1}(Tot ^{W, F}(V_{DR})/\Aquot _{wt}\times \Aquot _{hod}).
$$
Idem for $V_{B,\rr}$, and these maps are compatible with the
equivalences
$\zeta$. 

In the case we are currently considering, the terms in the above bilinear map are
assumed to be vector bundles over $\Aquot _{wt}\times \Aquot _{hod}$.
Thus the bilinear map becomes a morphism from the tensor product of the
two bundles, to the third one. 
The correspondence between filtered vector spaces in the usual sense and 
bundles over $\Aquot$ is compatible with tensor product. Taking the
tensor product of all of the filtrations in question, we obtain
a notion of tensor product of pre-mhs. The Whitehead product becomes
a morphism of
pre-mhs
$$
Wh ^{\rm pre}_{i,j}(\Vv ):
\pi _i(\Vv ) \otimes \pi _j (\Vv ) \rightarrow \pi _{i+j-1}(\Vv ).
$$
The alert reader will have noticed that there is a problem with
the shifts here: on the left we have tensored a $-i$-shifted mhs
with a $-j$-shifted mhs, which gives a $-i-j$-shifted mhs. On the
right is a $1-i-j$-shifted mhs. 

Here is where our main idea, of requiring that the associated-graded be a
linear object, comes in. The object $Gr^W(\Vv )$ is linear, so its Whitehead
products vanish. In other words, the associated-graded 
of the above Whitehead product for $\Vv$ vanishes:
$$
Gr ^W(Wh ^{\rm pre}_{i,j}(\Vv ))=0.
$$
In other words, $Wh ^{\rm pre}_{i,j}(\Vv )$ lowers the degree in the weight
filtration by one. Denote by $\pi _{i+j-1}(\Vv )(-D_{wt})$ the pre-mhs
obtained by shifting the weight filtration by one, in the direction such
that there is a morphism of pre-mhs inducing zero on the associated-graded,
$$
\pi _{i+j-1}(\Vv )(-D_{wt})\rightarrow \pi _{i+j-1}(\Vv ).
$$
Then $\pi _{i+j-1}(\Vv )(-D_{wt})$ is a $-i-j$-shifted mhs, and
the Whitehead product factors through a morphism
$$
Wh_{i,j}(\Vv ):
\pi _i(\Vv ) \otimes \pi _j (\Vv ) \rightarrow \pi _{i+j-1}(\Vv )(-D_{wt}).
$$
This is now a morphism of $-i-j$-shifted mixed Hodge structures.

If we normalize the homotopy group sheaves, and include an extra shift by 
one to normalize $\pi _{i+j-1}(\Vv )(-D_{wt})$ (which gives the same answer as if
we normalize $\pi _{i+j-1}(\Vv )$), then we have ended up shifting by $i+j$ on
both sides of the map $Wh_{i,j}(\Vv )$, so it normalizes to
a morphism of mixed Hodge structures 
(the {\em Whitehead product}):
$$
Wh_{i,j}^{\nu}(\Vv ):
\pi _i^{\nu}(\Vv ) \otimes \pi _j^{\nu} (\Vv ) \rightarrow 
\pi _{i+j-1}^{\nu}(\Vv ).
$$

The fact that the Whitehead product goes from degree $i$ times degree $j$
to degree $i+j-1$ means that in order to get a correctly shifted morphism
of mixed Hodge structures, we need to pick up a shift by $1$ somewhere. This
shift by $1$ comes from the fact that by assumption the Whitehead products
vanish on $Gr^W(\Vv )=DP (LGr ^W(\Vv ))$. This situation may inversely be taken
as one of the motivations \mylabel{whitemotivatepage}
for the idea of requiring $Gr^W(\Vv )$ to be a linear
object.

\begin{center}
{\bf The general case}
\end{center}

We now turn to the general case where we don't make the hypothesis that
the sheaves underlying $\pi _i(\Vv )$ are vector bundles. Recall from
\cite{geometricN} that if $X\rightarrow B$ is a relatively $1$-connected
geometric $n$-stack over a base $B$ then the $\pi _i (X/B)$ are vector
sheaves over $B$ (which if necessary means morphisms $B\rightarrow VESH$
to the $1$-stack of vector sheaves). Recall also from pages
\pageref{filtdefpage} and \pageref{filtvectsheafpage} that if $U$ is a vector
sheaf over $\Aquot$ then we obtain a filtration $W_{\cdot}$ of the vector space
$U_{[1]}$
by subspaces defined by the condition that $u\in W_k$ if and only if
$t^{-k}\tilde{u}$
extends to a section of $U$ over $\Aone$. This functor from vector
sheaves over $\Aquot$ to filtered vector spaces is unfortunately not an
equivalence of categories, and the fiber $U_{[0]}$ is no longer
identifiable with the usual associated-graded of the filtration (i.e.
the direct sum of the subspaces $W_k/W_{k-1}$).  This is what makes the
present situation more complicated.

We continue to assume that $\Vv$ is
a nonabelian mixed Hodge structure which is $1$-connected. Thus the
sheaves underlying $\pi _i(\Vv )$ are vector sheaves over the
appropriate base stacks $\Aquot_{wt}\times \Aquot _{hod}$ or 
$\Aquot _{wt, \rr }$. We obtain in particular, vector sheaves
$$
Tot ^F(\pi _i(V_{DR})) \;\;\; \mbox{over} \;\; \Aquot _{hod},
$$
and
$$
Tot ^W(\pi _i(V_{B, \rr })) \;\;\; \mbox{over} \;\; \Aquot _{hwt, \rr}.
$$
The fibers of these over $[1]$ are respectively the finite dimensional 
complex vector space $\pi _i(V_{DR})$ and
the finite dimensional 
real vector space
$\pi _i(V_{B,\rr })$; and the former is isomorphic via $\zeta$ to the 
complexification of the latter.  Let $\pi _i(V)$ denote these isomorphic
spaces, together with its real structure. From the above vector sheaves we obtain
filtrations $F^{\cdot}$ (complex) and $W_{\cdot}$ (real) of $\pi _i(V)$
by subspaces (again see pages \pageref{filtdefpage} and \pageref{filtvectsheafpage}).

This yields a pre-mixed Hodge structure 
$$
\pi _i^{\mu }(\Vv ) := \{ \pi _i(V), \pi _i(V)_{\rr}, F^{\cdot},
W_{\cdot} \} .
$$
Our goal is to obtain the following result.

\begin{theorem}
\mylabel{gencase}
The pre-mixed Hodge structure $\pi _i^{\mu }(\Vv )$ which we have defined above 
is a $-i$-shifted
mhs. Denote by $\pi _i^{\nu}(\Vv )$ 
its version where the weight filtration is shifted by $i$
places so that it becomes a mixed Hodge structure. 
The Whitehead products on the homotopy
groups of $V_{DR}$ give morphisms
of mixed Hodge structures 
$$
Wh^{\nu}_{i,j} : \pi _i^{\nu}(\Vv ) \otimes 
\pi _j^{\nu}(\Vv )  \rightarrow \pi _{i+j-1}^{\nu}(\Vv ) .
$$
\end{theorem}

Along the way we will also obtain a somewhat more precise structure
theorem about the vector sheaves underlying the original object $\pi
_i(\Vv )$.

\begin{center}
{\bf The case of mixed Hodge complexes}
\end{center}

Suppose $\Cc$ is a mixed Hodge complex. Its Dold-Puppe $\Vv = DP(\Cc )$ is a
namhs, and the homotopy groups of the latter are the same as the
cohomology groups of $\Cc$. Taking the cohomology sheaves of all of the
perfect complexes underlying $\Cc$ we obtain a ``pre-mixed Hodge vector
sheaf''
$H^{-i}(\Cc )$ and the Dold-Puppe compatibility says that this object is
equivalent
to the object $\pi _i(\Vv )$ considered above. In particular, we
encounter the same problem as that which was outlined above, for the
$H^{-i}(\Cc )$. On the other hand, Deligne shows in \cite{hodge3} that
the
cohomology groups of a mixed Hodge complex are (shifted) mixed Hodge
structures.
Our strategy for proving Theorem \ref{gencase} will be to reduce to the
case of a mixed Hodge complex. For this, we need a structure theorem
which says that in a certain sense one can ``split off'' the behavior at
the rightmost part of the mixed Hodge complex (which corresponds to the
lowest homotopy degrees). From there, we will be able to use the
Hurewicz fact that homotopy and homology coincide in low degrees to
apply the case of mixed Hodge complexes to the general case. 

The following two propositions and the resulting theorem
represent our ``splitting-off'' results. 
As they are results about mixed Hodge complexes, we leave their proofs
as exercises for the reader, to be done using what is known e.g. in
Deligne \cite{hodge3}, Hain and Zucker \cite{HainZucker},
Saito \cite{Saito} \cite{Saito2} \cite{Saito3} etc.

\begin{proposition}
\mylabel{splitoff1}
Suppose $\Cc$ is a mixed Hodge complex which is supported in $(-\infty , -k]$.
Let $C$ denote the underlying complex (i.e. one of the two isomorphic
complexes $C_{DR}$ or $C_{B,\cc}$). If $H^{-k}(C)\neq 0$ then $H^{-k}(C)$
provided with its Hodge and weight filtrations and real structure, is a
$-k$-shifted mixed Hodge structure. Furthermore, if we let $H^{-k}_{\mu}(\Cc )$ denote
this shifted mixed Hodge structure considered as a mixed Hodge complex sitting
in degree $-k$, then there is a morphism of mixed Hodge complexes
inducing the identity on $H^{-k}(C)$:
$$
\Cc \rightarrow H^0_{\mu}(\Cc ).
$$
\end{proposition}
\eop

\begin{proposition}
\mylabel{splitoff2}
Suppose $\Cc$ is a mixed Hodge complex  which is supported in $(-\infty , -k]$.
Let $C$ denote the underlying complex (i.e. one of the two isomorphic
complexes $C_{DR}$ or $C_{B,\cc}$). Suppose that  $H^{-k}(C)=0$.
Then there is a split $-k$-shifted
mixed Hodge structure $K= H^{-k}(Gr^W(\Cc ))$,
and a mixed Hodge complex $\Kk$ supported in $[-k-1,-k]$ such that 
$H^{-k}(Gr^W(\Kk ))= K$, such that $H^{-1-k}(Gr^W(\Kk ))=K'$ is the split
$-1-k$-shifted mixed Hodge structure obtained by shifting the weight
filtration of $K$ by one, and such that the complex underlying $\Kk$ is
exact. The differential of $\Kk$ is an isomorphism of pre-mixed Hodge
complexes but which shifts the weight filtration by one.
There is a morphism of mixed Hodge complexes 
$$
f:\Cc \rightarrow \Kk 
$$
which induces an isomorphism $H^{-k}(Gr^W(\Cc ))\cong H^{-k}(Gr^W(\Kk ))$.
Furthermore, $Cone (f)$ is quasi-isomorphic to a mixed Hodge complex
supported in $[-\infty , -1-k]$.
\end{proposition}
\eop

We show how to use these two propositions to understand the structure of
$H^i(\Cc )$. Suppose $\Cc$ is supported in $(-\infty , -k]$.
If $H^{-k}(C)\neq 0$ then we apply Proposition \ref{splitoff1} to obtain
a mixed Hodge complex $\Ll$ supported in degree $-k$ and a morphism 
$$
g:\Cc \rightarrow \Ll .
$$
Let $\Cc _1:= Cone(g)$, and let $C_1$ denote its underlying complex.
(If $H^{-k}(C)=0$ then set $\Cc _1=\Cc$ and pass directly to the
next step.)
Think of the cone as a kernel of a map, so it is supported in $(-\infty
, -k]$.
By the property obtained in Proposition \ref{splitoff1},
$H^{-k}(C_1)=0$. Therefore we can apply Proposition \ref{splitoff2} to
obtain a mixed Hodge complex $\Kk$ supported in $[-1-k, -k]$ and a
morphism
$$
f:\Cc _1 \rightarrow \Kk ,
$$
such that $\Cc _2 := Cone (f)$ is quasiisomorphic to a mixed Hodge
complex
supported in $(-\infty , -1-k]$. Thus in two steps we have split off two
types of mixed Hodge complexes: the first one equal to a mixed Hodge structure
supported in a single degree; and the second one supported in two adjacent degrees
and whose underlying complex is exact. After doing this, we move the
rightmost bound for the interval of support, back by one. 

We obtain the
following
structure theorem for the objects $H^{-i}(\Cc )$.

\begin{theorem}
\mylabel{mhcstructure}
Suppose $\Cc$ is a mixed Hodge complex supported in $(-\infty , 0]$.
Let $H^{-i}(\Cc )$ denote cohomology of $\Cc$ considered as a 
pre-mixed Hodge object in the category of
vector sheaves. Let $H^{-i}_{\mu }(\Cc )$ denote the bifiltered vector
space with real structure (i.e. pre-mhs) obtained by applying the
construction of filtrations of p. \pageref{filtdefpage} as above. Then there are 
objects $\Kk$ and $\Kk '$ (pre-mixed Hodge objects in the category of
vector sheaves) and morphisms
$$
a: H^{-i}(\Cc ) \rightarrow \Kk \rightarrow 0,
$$
$$
0 \rightarrow \Kk ' \stackrel{b}{\rightarrow} H^{-i}(\Cc )
$$
such that $\ker (a) / {\rm im}(b) \cong H^{-i}_{\mu}(\Cc )$.
Furthermore, the object $Tot ^{W,F}(\Kk )$ is a vector scheme over
$\Aquot _{wt}\times \Aquot _{hod}$, satisfying property ${\bf A1}$, 
whose fiber over $[1]\times \Aquot
_{hod}$
is $0$ and whose fiber over $[0]\times \Aquot _{hod}$ corresponds to a
split $-i$-shifted mixed  Hodge structure. Dually, the object 
$Tot ^{W,F}(\Kk ')$ is a coherent sheaf over $\Aquot _{wt}\times \Aquot
_{hod}$
whose dual vector scheme satisfies property ${\bf A1}$ (this means that
the coherent sheaf is annihilated by the coordinate $t$ on $\Aquot
_{wt}$).
Again its fiber over  $[1]\times \Aquot
_{hod}$
is $0$ and its fiber over $[0]\times \Aquot _{hod}$ corresponds to a
split $-i$-shifted mixed  Hodge structure.  Finally, the differential of
the mixed Hodge complex is an isomorphism between the split shifted mhs
corresponding to $\Kk$ for degree $-1-i$, and the split shifted mhs
corresponding to $\Kk '$ for degree $i$, except that this isomorphism
shifts the weight filtration by one. 
\end{theorem}
\eop

\begin{center}
{\bf Homology and cohomology of a geometric $n$-stack}
\end{center}

Suppose $X$ is a $1$-connected geometric $n$-stack relative to a base
$B$. Then we obtain its {\em cohomology object} ${\bf H}^{\ast}(X/B, \Oo )$
which is a perfect complex over $B$, supported in $[0, \infty )$. This 
may be seen as the higher direct image complex ${\bf R}p_{\ast}(\Oo )$
where $p: X\rightarrow B$ is the projection. 

Denote the dual perfect complex by 
$$
{\bf H}_{\ast}(X/B) := \underline{Hom}({\bf H}^{\ast}(X/B, \Oo ), \Oo ).
$$
It is a perfect complex over $B$, supported in $(-\infty , 0]$.
Its Dold-Puppe is an $\infty$-stack over $B$ such that any truncation is
``almost geometric''
in the sense of \cite{aaspects} p. 118. We have a morphism
$$
s:X \rightarrow DP({\bf H}_{\ast}(X/B )).
$$
This morphism is the ``stabilization map'', and one could write
$$
DP({\bf H}_{\ast}(X/B )) = \Omega ^{\infty} S^{\infty}(X/B).
$$
(While we don't go into that here, one could make a precise statment of
a lemma represented by the previous equation.)

The version of Hurewicz which applies to the present situation is
the following lemma:

\begin{lemma}
\mylabel{hurewicz}
Suppose $X$ is $k-1$-connected relative to $B$, for $k\geq 2$.
Then the above stabilization map $s$ induces an isomorphism between the
homotopy vector sheaf $\pi _k(X/B)$ and the cohomology vector sheaf 
$H^{-k}({\bf H}_{\ast}(X/B))$. In particular, the homotopy fiber of $s$
is $k$-connected.
\end{lemma}
\eop

\begin{center}
{\bf Homology of a namhs}
\end{center}

We now turn back to consideration of a $1$-connected nonabelian mixed
Hodge structure $\Vv$. We apply the discussion of the previous
subsection to the geometric $n$-stacks underlying $\Vv$, for example
$Tot^{W,F}(V_{DR})\rightarrow \Aquot_{wt}\times \Aquot_{hod}$.

We obtain homology objects such as 
$$
{\bf H}_{\ast}(Tot^{W,F}(V_{DR})/ \Aquot_{wt}\times \Aquot_{hod})
$$
which is a perfect complex over $\Aquot_{wt}\times \Aquot_{hod}$
supported in $(-\infty , -2]$. Doing this for $Tot^W(V_{B,\rr })$ and
the equivalence $\zeta$ too, we obtain all in all a pre-mixed Hodge
object in the $\infty$-stack of perfect complexes supported in 
$(-\infty , -2]$. Denote this object by ${\bf H}_{\ast}(\Vv )$. 
It is a pre-mixed Hodge complex supported in $(-\infty , -2]$.

\begin{lemma}
\mylabel{homology}
Suppose $\Vv$ is a $1$-connected namhs. 
Then the pre-mixed Hodge complex ${\bf H}_{\ast}(\Vv )$
constructed above is a mixed Hodge complex. 
\end{lemma}
{\em Proof:}
The associated-graded of ${\bf H}_{\ast}(\Vv )$ is the homology object
of $Gr^W(\Vv )$. But since this latter is by
hypothesis linearized, and since the Hodge filtration satisfies the
strictness hypothesis ${\bf Str}$, we obtain a noncanonical
decomposition of $Gr^W(\Vv )$ into a product of Eilenberg-MacLane
objects
for split shifted mixed Hodge structures. The homology object (which is
a coalgebra) is the tensor product of the homology objects for the
Eilenberg-MacLane factors. Each factor looks like $K(\Uu , i)$ where
$\Uu$ is a  split $-i$-shifted mixed Hodge structure. We have the formulae
$$
{\bf H}_{\ast}(K(\Uu , i)) = Sym^{\ast} (\Uu [i]) \;\;\;
\mbox{for} \; i \; \mbox{even},
$$
$$
{\bf H}_{\ast}(K(\Uu , i)) = \bigwedge^{\ast} (\Uu [i]) \;\;\;
\mbox{for} \; i \; \mbox{odd},
$$
where $\Uu [i]$ is the object 
$\Uu$ placed in degree $-i$. For these formulae, see the ``Breen
calculations'' in Part III below, page \pageref{breen1}. From these, 
it is immediate that 
${\bf H}_{\ast}(Gr ^W(\Vv ))$ is a split mixed Hodge complex. The mixed
Hodge complex conditions concern only $Gr^W$
so this implies that ${\bf H}_{\ast}(\Vv )$ is a mixed  Hodge complex.
\eop

The $i$-th cohomology of the mixed Hodge complex of cohomology
${\bf H}^{\ast}(\Vv )$ is a
pre-mixed Hodge object in $VESH$, and as usual we can extract from it 
a pre-mhs by taking the induced filtrations on the underlying vector
spaces. This yields in fact an $i$-shifted mixed Hodge structure which
we can then normalize by shifting the weight filtration to give a mixed
Hodge structure denoted $H^i_{\nu}(\Vv )$. Cup product is a morphism of
mixed Hodge structures, and we expect a statement about compatibility
between the higher-order Massey products and this mixed Hodge structure,
but we don't have a precise statement at the present.

\begin{center}
{\bf Homotopy fibers of certain morphisms of namhs}
\end{center}

In general, the condition of being a namhs is not stable under homotopy
fiber products. To see this one just has to look at the case of a mixed
Hodge complex equal to a single $-i$-shifted mhs placed in degree $-i$. Taking the
loop-space
of the associated namhs (which is a homotopy fiber product of two times
the punctual namhs, over the namhs in question) yields the same
$-i$-shifted
mhs but now placed in degree $1-i$. It is no longer a mixed Hodge
complex. 

Nonetheless, we need to have a result about preservation of the  namhs
condition when we take fibers of certain morphisms. A quick look at the
long exact sequence for homotopy groups (on the level of $Gr^W$, say) 
leads to the appropriate
condition:
the connecting morphisms in the long exact sequence are morphisms
between shifted mhs's with different shifts---there will be a problem if
these connecting morphisms are nonzero, so we make a hypothesis that
guarantees that the connecting  morphisms are zero. We obtain the
following statement.

\begin{lemma}
\mylabel{fiprod}
Suppose $f: \Vv \rightarrow \Uu $ is a morphism of $1$-connected
nonabelian mixed
Hodge structures. Suppose that $\Uu$ is pointed by a point $u$. Suppose
that for all $i$
the morphism induced by $Gr^W(f)$ is a surjection of homotopy group objects:
$$
Gr^W(f) : \pi _i Gr^W(\Vv )\rightarrow \pi _iGr^W(\Uu ) \rightarrow 0.
$$
Then the homotopy fiber $\Ff := \Vv \times _{\Uu}u$ of $f$ over $u$ is a
nonabelian mixed Hodge structure.
\end{lemma}
{\em Proof:}
The $n$-stacks underlying $\Ff$ are geometric and relatively $1$-connected;
therefore they are smooth over the respective base schemes so the annihilator
ideals vanish and we have automatically conditions ${\bf A1}$, ${\bf A2}$, 
${\bf A3}$. Also locally any basepoint section provides a smooth surjection
from a scheme so the flatness condition ${\bf Fl}$ holds too. Thus we just
have to verify that $LGr^W(\Ff )$ is a mixed Hodge complex, which contains
conditions ${\bf Str}$ and ${\bf MHC}$. The long exact sequence for the fiber
of a morphism yields a long exact sequence
(of split pre-mixed Hodge objects in the stack of vector sheaves)
$$
\ldots \rightarrow H^{-i}(LGr ^W(\Ff ))
\rightarrow H^{-i}(LGr ^W(\Vv ))\rightarrow
H^{-i}(LGr ^W(\Uu ))\rightarrow \ldots .
$$
The morphisms going from the places occupied by $\Vv$ to 
those occupied by $\Uu$, are surjections by hypothesis. Thus
the connecting morphisms (the leftmost and rightmost arrows in the above
diagram) are zero and 
$$
H^{-i}(LGr ^W(\Ff ))
=\ker \left( H^{-i}(LGr ^W(\Vv ))\rightarrow
H^{-i}(LGr ^W(\Uu )) \right) . 
$$
The total objects of the Hodge filtrations 
$$
Tot ^F(H^{-i}(LGr ^W(\Vv )) \;\;\; \mbox{and} \;\;\; 
Tot ^F(H^{-i}(LGr ^W(\Uu )))
$$ 
are vector bundles over $\Aquot _{hod}$, therefore
the kernel is a vector bundle too. This gives condition ${\bf Str}$ for $\Ff$.
Furthermore the split pre-mhs corresponding to $H^{-i}(LGr ^W(\Ff ))$
is the kernel of a morphism of $-i$-shifted mhs's, so it is again a $-i$-shifted
mhs. This gives condition ${\bf MHC}$. We obtain that $\Ff$ is a namhs.
\eop

\begin{center}
{\bf Conclusion: the proof of Theorem \ref{gencase}}
\end{center}

We now  put together all of the above statements to obtain a
proof of Theorem \ref{gencase}.  Along the way we will see that the structure
result of Theorem \ref{mhcstructure} also holds for the homotopy group objects
$\pi _i (\Vv )$. 

Suppose that $\Vv$ is $k-1$-connected for $k\geq 2$. The proof will be by
descending induction on $k$, so we assume that Theorem \ref{gencase} and indeed
the structure result analogous to Theorem \ref{mhcstructure} are known for
any $k$-connected namhs $\Ff$. (To be technically correct in the induction, 
we add to the inductive hypothesis that 
at the bottom degree $k+1$ there is no quotient object in the structure result.)

Let $\Cc := {\bf H} _{\ast} (\Vv )$. It is a mixed Hodge complex supported in 
$(-\infty , -k)$. Combining the constructions of Propositions \ref{splitoff1} and
\ref{splitoff2}, we obtain a mixed Hodge complex $\Nn$ supported in 
$[-1-k, -k]$ and a morphism $\Cc \rightarrow \Nn$ inducing surjections
on the cohomology objects of $Gr ^W(\Cc )$ and $Gr ^W(\Nn )$, and indeed inducing
an isomorphism of cohomology objects in degree $-k$. 

The complex
underlying $\Nn$ is exact in degree $-1-k$ and its cohomology in degree 
$-k$ carries a mixed Hodge structure $H^{-k}_{\mu} (\Cc )$.
The cohomology objects $H^{-i}(\Nn )$ are as follows: in degree
$-1-k$, an object such as was denoted $\Kk$ in Theorem \ref{mhcstructure};
and in degree $-k$ an extension with subobject as was denoted $\Kk '$ in 
Theorem \ref{mhcstructure} and quotient object $H^{-k}_{\mu} (\Cc )$.

The Dold-Puppe $DP (\Nn )$ has a structure of namhs (it is linearized by setting
$LGr ^W(DP (\Nn )):= Gr ^W(\Nn )$). Note that, whereas $DP (\Cc )$ is an
$\infty$-stack, $DP (\Nn )$ is an $n$-stack (or at most an $n+1$-stack if 
we are dealing with the limit case $k=n$). 

Now the fact that $k\geq 2$ means that the nonlinear terms in the Breen
calculations for the homology of $Gr ^W(\Vv )$ start intervening only in degree
$k+2$ or more. In particular the morphism
$$
Gr ^W(\Vv )\rightarrow DP (Gr ^W(\Cc ))
$$
induces a surjection on $\pi _k$ and $\pi _{k+1}$. Therefore the
morphism  $f:\Vv \rightarrow DP (\Nn )$ induces surjections on the homotopy group
objects of the associated-graded, in particular 
$f$ satisfies the hypotheses of Lemma \ref{fiprod}. That lemma says that 
the fiber $\Ff$ of 
$f$ over the zero-section of $DP(\Nn )$ is a nonabelian mixed Hodge structure.
The long exact sequence of homotopy objects gives a long exact sequence of
pre-mixed Hodge objects in the category of vector sheaves,
$$
\ldots
\rightarrow \pi _i(\Ff )\rightarrow \pi _i (\Vv )
\rightarrow H^{-i}(\Nn )\rightarrow \pi _{i-1}(\Ff )\rightarrow \ldots .
$$
Our inductive hypothesis says that for $i\geq k+2$ the objects 
$\pi _i (\Ff )$ satisfy the structure result of Theorem \ref{mhcstructure}
(and in particular that the pieces corresponding to bifiltered vector spaces with
real structure, are $i$-shifted mhs's); and it says that for $i=k+1$ the same
holds but with the quotient being zero in the structure result. 

Note that for $i\geq k+2$, $H^{-i}(\Nn )=0$ so $\pi _i(\Vv )=\pi _i (\Ff )$
and we obtain the desired structure result for $\pi _i (\Vv )$. 
For $i=k+1$ we have an exact sequence
$$
0\rightarrow \pi _{k+1}(\Ff ) \rightarrow \pi _{k+1}(\Vv )
\rightarrow H ^{-1-k}(\Nn ) \rightarrow 0.
$$ 
The fact that the connecting map is zero comes from the fact that
we constructed $\Nn$ so that $\pi _{k+1}(\Vv )
\rightarrow H ^{-1-k}(\Nn )$ is surjective. From 
this exact sequence we obtain the desired structure result for $\pi
_{k+1}(\Vv )$: the inductive hypothesis gives the structure result but without
quotient object, for $\pi _{k+1}(\Ff )$; and the piece $H^{-1-k}(\Nn )$
provides exactly a quotient object as forseen in Theorem \ref{mhcstructure}.
Finally we have by construction an isomorphism 
$$
\pi _{k}(\Vv )
\cong H ^{-k}(\Nn )
$$
(which justifies why $\Ff$ was $k$-connected). This gives the structure
result for $\pi _{k}(\Vv )$ with no quotient object, and in particular 
with $\pi _k^{\mu} (\Vv )$ being a $-k$-shifted mixed Hodge structure.

By induction we obtain the structure result of Theorem \ref{mhcstructure}
for the homotopy group objects $\pi _i(\Vv )$. This gives the first  part
of the statement of Theorem \ref{gencase}, that the $\pi _i ^{\mu }(\Vv )$
are $i$-shifted mixed Hodge structures.
It also gives the statement that the
associated-graded objects $Gr ^W\pi _i ^{\mu }(\Vv )$ are subquotients 
of the $\pi _i(Gr ^W(\Vv ))$. Now by the existence of a linearization,
the Whitehead products between the $\pi _i(Gr ^W(\Vv ))$ vanish. This implies
that the Whitehead products between the $Gr ^W\pi _i ^{\mu }(\Vv )$ vanish.
The same argument as above (starting page \pageref{whitehead}) 
gives that on the normalized mixed Hodge
structures $\pi _i ^{\nu}(\Vv )$, the Whitehead products are morphisms of mixed
Hodge structures. This completes the proof of Theorem \ref{gencase}.
\eop

We conclude this section by remarking that the above Lemma
\ref{fiprod} and its counterexample show that the the $n+1$-category $nNAMHS$ is
not closed under  homotopy limits, but only under very special kinds of them.
We wonder about the following:

\noindent
{\bf Question:} Is the $n+1$-category $nNAMHS$ closed under a reasonable class of
homotopy {\em colimits}?

\begin{center}
{\large \bf  The basic construction}
\end{center}
\mylabel{basicconstructionpage}

The basic construction which we would like to describe is the
``nonabelian cohomology'' of a smooth projective variety $X$ with
coefficients in a nonabelian mixed Hodge structure. In this section we
define the cohomology with coefficients in a pre-namhs $\Vv$ (eventually
linearized), denoted
$$
\Hh = \underline{Hom}(X_M, \Vv ).
$$
This will be an object of the same type as a pre-namhs, i.e. an object
in $PMH.nSTACK$.  We show that
it is in fact a pre-namhs (the meaning of this condition being that
the $n$-stacks involved are geometric or equivalently that it is an
object of $PMH.nGEOM$) under certian reasonable
hypotheses on $\Vv$. If $\Vv$ was linearized then $\Hh$ will be
linearized.

In what follows $X$ will be a smooth projective variety over $\cc$ and
$\Vv$ will be a pre-namhs (with the same notation as in the definition
page \pageref{pnamhspage}). As a matter
of notation, we will denote by $\Hh$ the object $\underline{Hom}(X_M,
\Vv )$ which we will construct, and as above we will have
$$
\Hh = \{ (H_{DR}, W,F),\;\; (H_{B,\rr } , W ), \;\; \zeta _{\Hh} \} .
$$
We first describe what is meant by $X_M$. This is an object which
is somewhat like a pre-namhs, however the $n$-stacks involved are not
geometric, and the morphism $\zeta$ is not an equivalence.

Recall from \cite{santacruz} \cite{aaspects} that we have defined a
formal category $X_{Hod}$ over $\Aquot _{hod} $ which gives a
$1$-stack
$$
X_{Hod}\rightarrow \Aquot _{hod}.
$$
Technically speaking in \cite{santacruz}
we denoted by $X_{Hod}$ the $\Gm$-equivariant
formal category over
$\Aone _{hod}$; for our purposes here it is more convenient 
\mylabel{convenientpage} to take
the quotient and denote by $X_{Hod}$ the object over $\Aquot _{hod}$.

The fiber of $X_{Hod}$ over $[1]$ is $X_{DR}$ and the fiber over $[0]=B\Gm$ is
$X_{Dol}$ with its action of $\Gm$. On the other hand, we define
$X_{B,\rr}$ to be the constant $n$-stack over $\Zz _{\rr}$ whose value
is the Poincar\'e $n$-groupoid corresponding to the $n$-truncation of the homotopy
type of $X^{\rm top}$. Note that the complexified stack associated to a
constant real stack is again a constant stack on $\Zz$ with the same
value, in particular $X_{B,\cc}$ is the constant stack on $\Zz$ whose
value is the Poincar\'e $n$-groupoid of $X^{\rm top}$. 

In order to define the
morphism
of analytic $n$-stacks
$$
\zeta _X: X_{DR}^{\rm an} \rightarrow X_{B,\cc}^{\rm an},
$$
we need the following lemma.

\begin{lemma}
Define the $n$-prestack $C^{\rm pre}(-, X^{\rm top})$ on $\Zz ^{\rm an}$ 
to be the prestack which to each $Y$ associates the
Poincar\'e $n$-groupoid of the space of continuous maps from $Y^{\rm
top}$ to $X^{\rm top}$. Then the associated $n$-stack
which we denote by $C(-, X^{\rm top})$, is naturally equivalent to 
$ X_{B,\cc}^{\rm an}$. 
\end{lemma}
{\em Proof:}
Let $X_{B, pre, \cc}$ be the constant $n$-prestack with values 
equal to the Poincar\'e $n$-groupoid of $X^{\rm top}$, and let
$X^{\rm an}_{B, pre, \cc}$ be the associated analytic prestack which is again a
constant prestack. Then $X_{B,\cc }$ (resp. $X_{B,\cc}^{\rm an}$) is the 
$n$-stack associated to $X_{B, pre, \cc}$ (resp. $X^{\rm an}_{B, pre, \cc}$).
We have a morphism of prestacks obtained by considering a point of $X^{\rm top}$
as a constant map from any $Y^{\rm top}$ to $X^{\rm top}$,
$$
X^{\rm an}_{B, pre, \cc} \rightarrow C^{\rm pre}(-, X^{\rm top}).
$$ 
This morphism is an equivalence, over any $Y$ such that $Y^{\rm top}$ is
contractible. Since the contractible open sets form a base for the topology
of any $Y\in \Zz ^{\rm an}$, this property implies that the above 
morphism of prestacks
induces an equivalence on associated $n$-stacks
$$
X^{\rm an}_{B,  \cc} \stackrel{\cong}{\rightarrow} C(-, X^{\rm top}).
$$ 
\eop

Recall that 
$$
X_{DR}^{\rm an}(Y) = X(Y^{\rm red})
$$
where $Y^{\rm red}$ is the reduced complex analytic space. A morphism
from $Y^{\rm red}$ to $X$ gives a continuous map from 
$Y^{\rm top}$ to $X^{\rm top}$, so we get a morphism 
$$
X_{DR}^{\rm an}(Y) \rightarrow 
C(-, X^{\rm top})(Y) \cong X_{B,\cc}^{\rm an} (Y).
$$
Our $\zeta _X$ is defined to be the composed morphism
$X_{DR}^{\rm an}\rightarrow X_{B,\cc}^{\rm an}$.

The object which we denote informally by $X_M$ is the collection
$$
X_M = \{ X_{Hod} \rightarrow \Aquot _{hod} , X_{B,\rr }, \zeta _X \} .
$$
Note that $X_{Hod}$ is an object of $F.nSTACK$. We can consider it as an
object of $F.F.nSTACK$ by making it constant in the other ``weight''
direction;
and similarly the $n$-stack $X_{B,\rr}$ may be thought of as an object
of $F.nSTACK$ constant over $\Aquot _{wt}$. In this sense the object
$X_M$ is somewhat similar to a pre-namhs. The important differences are:
\newline
--the $n$-stacks involved are not geometric; and
\newline
--the morphism $\zeta _X$ is a morphism but not an equivalence.

In spite of these differences, we still define an object
$\Hh = \underline{Hom}(X_M, \Vv )$ as we shall now start to do.

\begin{center}
{\bf The nonlinearized case}
\end{center}
\mylabel{nonlincase}

We will first construct $\Hh$ in the case where the pre-namhs $\Vv$ is
not linearized.
Put:
$$
H_{DR} := \underline{Hom}(X_{DR}, V_{DR}).
$$
It is a geometric $n$-stack by \cite{kobe}, \cite{geometricN}.

As indicated above we obtain objects which are constant in the weight
direction
$$
X_{Hod}\times \Aquot _{wt} \rightarrow \Aquot _{hod} \times \Aquot
_{wt}  ,
$$
and
$$
X_{B,\rr } \times \Aquot _{wt, \rr} \rightarrow 
\Aquot _{wt,\rr} .
$$
Restricting the first to $[1]\in \Aquot _{hod}$ we obtain
$$
X_{DR}\times \Aquot _{wt} \rightarrow  \Aquot _{wt}  .
$$

Now we define the weight and Hodge
filtrations on this $n$-stack by defining the total space as
$$
Tot ^{W,F}(H_{DR}) := \underline{Hom}^{{\rm se}, 0}\left( 
\frac{X_{Hod}\times \Aquot _{wt} }{ \Aquot _{hod} \times \Aquot _{wt}} ,
 \frac{Tot^{W,F}(V_{DR})}{\Aquot _{hod} \times \Aquot
_{wt}}\right) .
$$
Here $\underline{Hom}(\frac{A}{B} , \frac{C}{B})$ is the relative
internal $Hom$ of $n$-stacks over $B$. It can also be viewed as the
internal $Hom$ in the stack $(n+1)STACK _{1/}(B, nSTACK )$ (there is a
compatibility between these two things, which we don't get into). 

The member on the left of
the internal $\underline{Hom}$ is a formal category of smooth type over
the base $B= \Aquot _{hod} \times \Aquot
_{wt}$. 

The superscript means that we take the 
substack of {\em semistable morphisms with vanishing Chern classes}, 
see the discussion in \cite{aaspects} \S 10.7 
and on page \pageref{semistablepage} below.
This semistability condition is automatic for Betti or de Rham
cohomology, so it is only necessary for the Dolbeault
cohomology,
i.e. at the associated-graded for the Hodge filtration, over $[0]_{hod}
\subset \Aquot_{hod}$. A morphism $X_{Dol}\rightarrow V$ is said to be {\em
semistable with vanishing Chern classes} 
if the morphism $X_{Dol}\rightarrow \pi _0(V)$ is constant, say with
values at a point $t\in \pi _0(V)$, and if the morphism 
$X_{Dol}\rightarrow \tau _{\leq 1}(V)$ which becomes essentially
$$
X_{Dol}\rightarrow K(\pi _1(V,t), 1),
$$
is a semistable principal Higgs $\pi _1(V,t)$-torsor on $X$ with vanishing
rational Chern classes. 
One thing that we need to know here is that the condition of
semistability with vanishing
Chern classes is an open condition in the $1$-stack of morphisms toward $\tau _{\leq 1}V$.
This is Proposition \ref{openness} below. 

If $\Vv$ is simply connected, i.e. the fundamental group object is trivial
$$
\pi _1(Tot ^{W,F}(V_{DR})/\Aquot _{wt} \times \Aquot _{hod}) = \{ 1\},
$$ 
then
the semistability condition is automatic. In this case we can apply
Theorem 10.3.3 of \cite{aaspects}
to get that $Tot ^{W,F}(H_{DR})$ is a geometric $n$-stack
over $\Aquot _{wt}\times \Aquot _{hod}$.  For now we will restrict to the
case where $\Vv$ is simply connected in our statement (Theorem
\ref{basicconstruction}) 
about $\Hh$. 

Later on
in Proposition \ref{boundedness} we show that if the fundamental group object
of $\Vv$ is a flat linear group scheme then $Tot ^{W,F}(H_{DR})$ is
geometric; this leads to the statement of Theorem \ref{flcase} which extends Theorem 
\ref{basicconstruction}
to cover this case.

The total space of the weight filtration on
$H_{DR}$
has a formula similar to the previous one, namely:
$$
Tot ^W (H_{DR}) := \underline{Hom} \left(
\frac{X_{DR} \times \Aquot _{wt}  }{\Aquot _{wt} } , 
\frac{Tot^W(V_{DR})}{\Aquot _{wt} }
\right) .
$$

Next, we define similarly 
$$
H_{B, \rr} := \underline{Hom}(X_{B,\rr}, V_{B,\rr}).
$$
It is a real geometric $n$-stack cf \cite{geometricN} Proposition 2.1 and the
comment on page 12, which refers to \cite{kobe} Corollary 5.6;  also in the
present
real case we use Lemma
\ref{realgeometric} above saying it suffices to check that  $H_{B,\cc}$
is geometric.

As before the weight filtration on $H_{B,\rr}$ is defined by the formula
$$
Tot ^W (H_{B,\rr}) := \underline{Hom} \left(
\frac{X_{B,\rr} \times \Aquot _{wt}  }{\Aquot _{wt} } , 
\frac{Tot^W(V_{B,\rr})}{\Aquot _{wt}}
\right) .
$$
This is geometric by \cite{geometricN} Proposition 2.1, \cite{kobe}
Corollary 5.6,  and Lemma \ref{realgeometric} above.

The analytic stacks associated to the algebraic $n$-stacks
$Tot^W(H_{DR})$ and $Tot ^W(H_{B,\cc})$ are given again by the same
formulae, namely:
$$
Tot ^W (H_{DR})^{\rm an} := \underline{Hom} \left(
\frac{X_{DR}^{\rm an} \times \Aquot _{wt}^{\rm an}  }{\Aquot _{wt}^{\rm an}  } , 
\frac{Tot^W(V_{DR})^{\rm an}}{\Aquot _{wt}^{\rm an} }
\right) 
$$
and
$$
Tot ^W (H_{B,\cc})^{\rm an} := \underline{Hom} \left(
\frac{X_{B,\cc}^{\rm an} \times \Aquot _{wt}^{\rm an} }{\Aquot _{wt}^{\rm an} } , 
\frac{Tot^W(V_{B,\cc})^{\rm an}}{\Aquot _{wt}^{\rm an} }
\right) .
$$

Using these last, our map $\zeta _X$ gives rise to a map
$$
\zeta _{\Hh} : Tot ^W (H_{B,\cc})^{\rm an}  \rightarrow 
Tot ^W (H_{DR})^{\rm an},
$$
and the ``GAGA theorem'' of \cite{kobe} says that this map is an
equivalence. (Note that here we come across, for the first time, a case
where the techniques of \cite{kobe} are needed for nonabelian cohomology
with coefficients in stacks where $\pi _0$ is nontrivial, in this case
the $\pi _0$ is the base of the weight filtration, $\Aquot _{wt}$.)

We have now constructed an object $\Hh$ in $PMH.nSTACK$, and if $\Vv$ is 
simply connected then $\Hh$ is in $PMH.nGEOM$. In Theorem \ref{flcase}
we will get $\Hh$ in $PMH.nGEOM$ if 
$$
\pi _i(Tot ^{W,F}(V_{DR})/\Aquot _{wt} \times \Aquot _{hod})
$$
is a flat linear group scheme over $\Aquot _{wt} \times \Aquot _{hod}$.

\begin{center}
{\bf The case of coefficients in a pre-mhc}
\end{center}
\mylabel{mhccase}

In order to treat the case of coefficients in a Dold-Puppe linearized 
object $\Vv$, we have to first treat the construction of an object
$\underline{Hom}(X_M, \Cc )$ when $\Cc$ is a split pre-mhc (this will then be
applied to $\Cc = LGr^W(\Vv )$). In fact it is
obviously interesting to do this even when $\Cc$ is not split. Thus in this
section we will define a pre-mhc $\underline{Hom}(X_M, \Cc )$ whenever
$\Cc$ is a pre-mhc.

We start with the following general construction. Suppose $B$ is a base $n$-stack
(of $n$-groupoids).
Suppose $\Ff \rightarrow B$ is an $n$-stack over $B$ and suppose that $\Cc$ is 
a complex of sheaves of abelian groups (which are assumed to be
rational vector spaces) on $B$ supported in $[-n, 0]$ 
(i.e. $\Cc$ is a morphism from $B$ to the $n$-stack $nCPX$ of complexes 
of sheaves of abelian groups). Then we can define the relative morphism complex
$$
\Mm = \underline{Hom} \left(\frac{\Ff}{B}, \frac{\Cc}{B}\right)
$$
which is a complex of sheaves of abelian groups (again, rational vector spaces)
over $B$. 
One way of doing this is to set $\Mm$ equal to the higher direct image complex of
the restriction $\Cc |_{\Ff}$ via the morphism $p:\Ff \rightarrow B$, which is then
truncated in positive degrees to give a complex supported in degrees $[-n,0]$:
$$
\Mm := \tau _{\leq 0} {\bf R} p_{\ast} (p^{\ast} \Cc ).
$$
Another way of doing it is to let $\Cc$ correspond to a relatively $N$-truncated 
$n+N$-stack $\tilde{\Cc}$ over $B$, to put 
$$
\tilde{\Mm}' := \underline{Hom}\left(\frac{\Ff}{B},
\frac{\tilde{\Cc}}{B}\right) ,
$$
and to let 
$$
\tilde{\Mm} := {\rm fib}(\tilde{\Mm}' \rightarrow \tau _{\leq n}(\tilde{\Mm}')).
$$
The fiber is taken over the zero-morphism as basepoint section, where the 
zero-morphism is the constant morphism with values the basepoint section of 
$\tilde{\Cc}$.  Now $\tilde{\Mm}$ is a relatively $N$-connected $n+N$-stack
over $B$ so it corresponds to a complex of sheaves of abelian groups.

The above construction is compatible with
the usual $\underline{Hom}$ of $n$-stacks via Dold-Puppe:
$$
\underline{Hom}\left(\frac{\Ff}{B}, \frac{DP(\Cc )}{B}\right)=
DP\underline{Hom}\left(\frac{\Ff}{B}, \frac{\Cc}{B}\right).
$$

The case which will interest us is when $\Ff$ is a formal category over $B$.
See for example \cite{aaspects} for the definitions concerning this
notion.
Note that the definitions of formal category and the various properties
of such, were defined in \cite{aaspects} for the case where the base $B$
was a scheme. These are immediately extended to the case where $B$ is
any $n$-stack by the usual trick of looking at a formal category over
$B$ (resp. one having certain properties)
as a morphism from $B$ into the $2$-stack of formal categories
(resp. the substack of those having the properties in question). 

The following lemma is basically the same as Theorem 10.2.5 of
\cite{aaspects}, except that we have added the extension to the case of
hemiperfect complexes. 

\begin{lemma}
\mylabel{cohfclem}
Suppose $\Ff$ is a formal category over $B$ which is projective and of smooth
type relative to $B$. Suppose $\Cc$ is a hemiperfect complex over $B$. Then the morphism complex
$\Mm = \underline{Hom} (\frac{\Ff}{B}, \frac{\Cc}{B})$ is a hemiperfect complex over
$B$.
\end{lemma}
{\em Proof:} To check that $\Mm$ is hemiperfect it suffices (by definition)
to check that this is the case over any scheme $Z\rightarrow B$. Thus we may 
assume that $B$ is a scheme, and in fact we may when necessary replace $B$ by an
open subset (part of an open cover). With this reduction we may assume that there
is a perfect complex $\Kk$ supported in $[-n,1]$ such that 
$\Cc = \tau _{\leq 0}\Kk$ (put $\Kk ^i=\Cc ^i$ for $i<0$ and let $d:\Kk
^0\rightarrow \Kk ^1$ be a morphism which has the vector scheme $\Cc ^0$ 
as its kernel).  Recall that $p$ denotes the morphism $\Ff \rightarrow B$.
We  have
$$
\Mm := \tau _{\leq 0} {\bf R}p_{\ast} (p^{\ast} \Cc ) = 
\tau _{\leq 0} {\bf R}p_{\ast} (p^{\ast} \Kk ).
$$
Thus it suffices to prove that 
$\tau _{\leq 0} {\bf R}p_{\ast} (p^{\ast} \Kk) $ is a perfect complex 
supported in $[-n,\infty )$. On the other hand, $\Kk$ is, up to quasiisomorphism,
made up by taking successive cones of shifts of vector bundles. Perfection is
preserved by cones, and higher direct image preserves triangles so up to
quasiisomorphism it takes cones to cones. Shifts also commute with higher direct
image and preserve perfection. Finally, by replacing $B$ by an open cover we may
assume that the vector bundles are trivial bundles, so we are reduced to claiming
that 
$$
{\bf R}p_{\ast} (\Oo )
$$
is a perfect complex on $B$. Since $\Ff$ is of smooth type relative to
$B$ (i.e. the morphism $\Ff \rightarrow B$ is of smooth type in the
terminology of \cite{aaspects} 8.6),
its underlying scheme
$X$ is flat over $B$, and the 
higher direct image is calculated as the higher direct image
of the ``complex of differentials'' for $\Ff$ (cf \cite{aaspects} 10.2.5): 
$$
{\bf R}p_{\ast} (\Oo )= {\bf R} p_{X/B,\ast }( \Omega ^{\cdot} _{X/\Ff} ).
$$
The condition that $\Ff$ be of smooth type (cf \cite{aaspects})
means that 
$\Omega ^{\cdot} _{X /\Ff }$ is a differential complex of locally free sheaves on
$X$. Since $X$ is by hypothesis projective over $B$, it is well-known that the
higher direct image is a perfect complex. This completes the proof.
\eop

We have an analogue of the previous lemma designed to be applied to the
``Betti''
case.

\begin{lemma}
\mylabel{cohbetlem}
Suppose $\Ff$ is an $n$-stack over $B$ which is locally constant with
fiber the homotopy type of a finite CW complex.
Suppose $\Cc$ is a hemiperfect complex over $B$. Then we can define as
previously
a morphism complex
$\Mm = \underline{Hom} (\frac{\Ff}{B}, \frac{\Cc}{B})$ which is  hemiperfect over
$B$.
\end{lemma}
{\em Proof:}
Left as an exercise.
\eop

Suppose $\Cc$ is a pre-mhc, and suppose $X$ is a smooth projective
variety. By using exactly the same formulae as in the previous
construction starting on page \pageref{nonlincase} we can define a
pre-mhc $\Hh := \underline{Hom}(X_M, \Cc )$. In doing this translation, replace
$\Vv$ by $\Cc$ (thus, $V_{DR}$ by $C_{DR}$ and so forth), and 
use the morphism complexes of Lemmas \ref{cohfclem} and
\ref{cohbetlem} in place of the usual internal $\underline{Hom}$ stacks
used in the construction starting on page \pageref{nonlincase}. We don't
write out the formulae here because they are exactly the same as in the
previous construction.  Again the map $\zeta _X$ gives rise to the
required
equivalence $\zeta _{\Hh}$ between $H_{DR}$ and $H_{B,\cc }$. Note that in this case,
analytic equivalences between hemiperfect complexes are algebraic
so $\zeta _{\Hh}$ is actually an algebraic equivalence. 

If $\Cc$ is a split pre-mhc then $\underline{Hom}(X_M, \Cc )$ is also a
split
pre-mhc. 
Indeed, the action of $\Gm$ splitting the weight filtration on $\Cc$ induces an
action on $\underline{Hom}(X_M, \Cc )$ which again splits the weight
filtration; this is because the object $X_M$ is constant in the
weight-filtration direction (i.e. in the direction of $\Aquot _{wt}$)
so it doesn't contribute to changing the degree in the weight
filtration. Thus, 
$\underline{Hom}(X_M, \Cc )$ is split as a pre-mhc.

\begin{center}
{\bf The case of coefficients in a linearized pre-namhs}
\end{center}

We now turn to the case which combines the two previous ones. 
Suppose $\Vv$ is a linearized pre-namhs. Recall that this means that
$\Vv$ is a pre-namhs and that we have a split pre-mhc $LGr^W(\Vv )$
whose Dold-Puppe is equivalent to $Gr^W(\Vv )$ as split pre-namhs's. 
We now define a pre-namhs 
$$
\Hh = \underline{Hom}(X_M, \Vv )
$$
as follows. The underlying pre-namhs is the one which was constructed
starting on page \pageref{nonlincase}. We set
$$
LGr^W(\Hh ):=  \underline{Hom}(X_M, LGr^W(\Vv )),
$$
which is a split pre-mhc. The above constructions are compatible with
Dold-Puppe
(i.e. the morphism complex defined in Lemmas \ref{cohfclem} and
\ref{cohbetlem}
are compatible with the
usual relative internal $\underline{Hom}$ stacks via Dold-Puppe), so 
we have a canonical equivalence 
$$
DP (LGr^W(\Hh )) \cong Gr^W(\Hh ).
$$
This completes our construction of a linearized pre-namhs $\Hh$.

We sum up the results of this and the preceding paragraphs:

\begin{theorem}
\mylabel{basicconstruction}
Suppose $X$ is a smooth projective variety and $\Vv$ is a pre-namhs
(resp. pre-mhc, resp. linearized pre-namhs).
Suppose that $\Vv$ is simply connected (resp. no additional hypothesis
in the pre-mhc case, resp. the same hypothesis in the linearized
pre-namhs case).
Then the object 
$$
\underline{Hom}(X_M, \Vv ) = \Hh = \{ (H_{DR}, W,F),\;\; (H_{B,\rr}, W), \;\; \zeta _{\Hh} \}
$$
defined above is a pre-namhs (resp. pre-mhc, resp. linearized pre-namhs).
\end{theorem}
\eop

\begin{center}
{\bf The mixed Hodge complex conditions}
\end{center}

The following proposition comes from Deligne's
construction in``Hodge III'' \cite{hodge3}.

\begin{theorem}
\mylabel{hodge3plus}
Suppose $\Cc$ is a truncated mixed Hodge complex (see page
\pageref{truncationmhcpage}) and
suppose that $X$ is a smooth projective variety. Then the pre-mixed
Hodge complex
$\underline{Hom}(X_M, \Cc )$ that we have constructed above 
is a truncated mixed Hodge complex, i.e. it satisfies conditions
${\bf Str}$ (p. \pageref{strdef}), ${\bf MHC}$ (p. \pageref{puritypage}), 
${\bf A1}$ (p. \pageref{a1a2}) and ${\bf A3}$ (p. \pageref{a3}). 
\end{theorem} 
{\em Proof:}
There is an actual mixed Hodge complex $\Dd$ (for example supported in 
$[-n, 1]$) such that $\Cc$ is the truncation of $\Dd$ into degrees $\leq 0$.
Then $\underline{Hom}(X_M, \Cc )$ is the truncation of $\underline{Hom}(X_M, \Dd
)$. If we show that the latter is a truncated mixed Hodge complex, then by 
Corollary \ref{truncatingmhc}
its truncation will be a truncated mixed Hodge complex. Thus we are reduced to
the case where $\Cc$ is a mixed Hodge complex.  Furthermore one shows $(\ast )$,
in fact, that the full cohomology complex of $X_M$ with coefficients in $\Cc$ is 
a mixed Hodge complex supported in $[-n, \infty )$; the complex we denote
$\underline{Hom}(X_M, \Cc )$ is the truncation of this, so it is
a truncated mixed  Hodge complex. 

To show $(\ast )$ that the cohomology is a mixed Hodge complex, it suffices to treat
the associated-graded; but $Gr^W(\Cc )$ decomposes as a direct sum of
shifted split mixed  Hodge structures. Thus it suffices to treat the 
case of cohomology with coefficients in a pure Hodge structure. It is
well-known that it gives back a pure Hodge structure; the only problem is to
verify that the shifts of weights work out appropriately. 
For this, we refer to the
discussion on page \pageref{footnotepage}. 
\eop

\begin{corollary}
If $\Cc$ is a split mixed Hodge complex (cf p. \pageref{splmhc}), and
$X$ is a smooth projective variety, then the pre-mhc
$\underline{Hom}(X_M, \Cc )$ is a split mixed Hodge complex.
\end{corollary}
{\em Proof:}
As remarked above, the splitting carries over to 
$\underline{Hom}(X_M, \Cc )$.
We just have to show that the truncated mixed Hodge complex is actually
a mixed Hodge complex. 
The only direction in which the hemiperfect complexes involved in 
$\underline{Hom}(X_M, \Cc )$
might not be perfect, is in the $\Aquot _{wt}$ direction (in the $\Aquot
_{hod}$ direction they are perfect complexes because of the strictness
conditon ${\bf Str}$). However, the splitting of the weight filtration
means that in the  $\Aquot _{wt}$ direction the complexes are constant
so they are perfect complexes. Thus $\underline{Hom}(X_M, \Cc )$ is a
split mixed  Hodge complex.
\eop

\begin{corollary}
\mylabel{mhcatleast}
Suppose in the situation of Theorem \ref{basicconstruction} that $\Vv$
is a linearized pre-namhs and satisfies the strictness condition ${\bf
Str}$
(Lemma \ref{strdef} p. \pageref{strdef}) and
the mixed Hodge complex condition ${\bf MHC}$ (p. \pageref{puritypage}).
Then the linearized pre-namhs $\Hh = \underline{Hom}(X_M, \Vv )$ also satisfies 
the strictness condition ${\bf Str}$ and the condition
${\bf MHC}$.
\end{corollary}
{\em Proof:}
These conditions are measured on the $LGr^W(\Vv )$ and
$LGr^W(\Hh )$, so Theorem \ref{hodge3plus} and the previous corollary apply.
\eop

\begin{center}
{\bf A couple of further remarks}
\end{center}

The construction which to $X$ and
$\Vv$ associates $\Hh$ is functorial in $\Vv$ and contravariantly
functorial in $X$ (i.e. it gives an
$n+1$-functor from the product of the category of smooth projective
varieties with the  $n+1$-category of connected pre-namhs's (resp.
connected linearized pre-namhs's), to the
$n+1$-category of pre-namhs's (resp. linearized pre-namhs's)).

It is relatively clear from the above construction 
(namely, from the fact that $X_M$ was taken to be constant in the
$\Aquot_{wt}$-direction) that $\Vv
\mapsto \Hh$
depends only on the data of the {\em very presentable shape} of
$X_{Hod}/\Aquot _{hod}$ together with the very presentable shape of
$X_B$ and the equivalence between the Betti shape and the de Rham shape
after composing with the analytification functor. In the future one
hopes for a construction dealing with singular and open varieties, in
which case there will surely be weight data contained in the appropriate
generalization of $X_M$.

\begin{center}
{\bf Semistable morphisms}
\end{center}

\mylabel{semistablepage}

In this section, we look more closely at $\underline{Hom}(X_M, \Vv )$
when $\Vv$ is a pre-namhs which is not necessarily simply connected.
This allows us to give a statement like Theorem \ref{basicconstruction} for the case of
$\Vv$ not simply connected. The results of this section are not used in
Part III.

For the purposes of this discussion and to reduce the volume of
notation,
we denote $B:= \Aquot _{wt}\times \Aquot _{hod}$ and
$$
\Ff := Tot ^{W,F}(X_{DR})=X_{Hod}\times \Aquot _{wt} \rightarrow B
$$
which is a formal category, projective and of smooth type over $B$.
The important properties of this situation are:
\newline
(1) \, There is a closed substack $B_1\subset B$ such that 
$$
\Ff |_{B_1} = X_{Dol}\times B_1,\;\;\; 
\Ff |_{B-B_1} = X_{DR}\times (B-B_1);
$$
(2) For cohomology with vector sheaf coefficients, $\Ff /Z$ is of finite
cohomological dimension; and
\newline
(3) There are only finitely many equivalence classes of $Spec (\cc
)$-valued
points in $B$.

We shall make use of these properties to get around what seem otherwise
to be some delicate questions. 

To continue with our notation, we fix a pre-namhs $\Vv$ and we let $V:=
Tot^{W,F}(V_{DR})$
denote the component which is a geometric $n$-stack over $B$. We assume
that this is connected relative to $B$. Suppose $Z$ is a scheme mapping
to $B$.
Let $\Ff _Z:= \Ff \times _BZ$ and $V_Z:= V\times _BZ$; by going to an
etale  neighborhood in $Z$ if necessary, we may assume that there is a
section
$v: Z\rightarrow V_Z$. Now put $G_Z:= \pi _1(V_Z/Z, v)$. This is a very
presentable group sheaf over $Z$. 

\noindent
{\bf Hypothesis:} That the fundamental group object of $\Vv$ is
a flat linear group scheme; this means that for every
scheme $Z$ mapping to $B$ and basepoint section $v: Z\rightarrow V_Z$
as above, the group sheaf $G_Z= \pi _1(V_Z/Z, v)$ is represented
by a flat linear group scheme over $Z$.

(Note that there might not exist a basepoint of $\Vv$ over $B$ so the
fundamental group object doesn't exist {\em per se}; the condition on
$G_Z$ for
every scheme $Z$ is thus the definition of the above hypothesis.)

We recall the definition of the sub-$n$-stack
``of semistable morphisms with vanishing Chern classes'',
$$
\underline{Hom}^{{\rm se}, 0} (\Ff / B, V/B )\subset 
\underline{Hom} (\Ff / B, V/B ).
$$
If $Z\stackrel{f}{\rightarrow} B$ is a scheme mapping to $B$, 
we shall say when a point 
$\varphi : \Ff _Z \rightarrow V_Z$ of the right hand $n$-stack, is a
member of the substack. For the purposes of stating this condition we
may replace $Z$ by an etale open neigborhood and assume that $V_Z$
admits
a basepoint section $v$ (the condition will be independent of the choice
of
basepoint section so the condition glues to give a condition over any
$Z$).

Given a map $\varphi $ as above, for any closed point $z\in Z$ we obtain
a morphism $\varphi _z: \Ff _z \rightarrow V_z$. The first truncation 
in the Postnikov tower for $V_z$ is 
$$
V_z \rightarrow \tau _{\leq 1}(V_z)=K(G_z, 1).
$$
Thus $\varphi _z$ projects to give a morphism $\rho _z: \Ff
_z\rightarrow K(G_z,1)$, i.e. a principal $G_z$-bundle over $\Ff _z$.

Recall from (1) above that there are two cases. If $f(z)\in (B-B_1)$
then
$\Ff _z = X_{DR}$ and $\rho _z$ is a principal $G_z$-bundle with flat
connection over $X$. In this case the condition of being semistable with
vanishing Chern classes is automatic, and we admit $\varphi _z$. If
$f(z)\in B_1$ then $\Ff _z=X_{Dol}$ and $\rho _z$ is a principal
$G_z$-Higgs bundle over $X$. In this case we say that $\varphi _z$ is
``semistable with vanishing Chern classes'' if the rational Chern
classes of the principal bundle $\rho _z$ vanish for all invariant
polynomials on $G_z$, and if for every linear representation $E$ of
$G_z$,
the associated Higgs bundle $E(\rho _z)$ is semistable. 

Now we say that a point $\varphi : \Ff _Z \rightarrow V_Z$ is
``semistable with vanishing Chern classes'' if for every point $z\in Z$
the restriction $\varphi _z$ is semistable with vanishing Chern classes
in the above sense. We define 
$$
\underline{Hom}^{{\rm se}, 0} (\Ff / B, V/B )(Z)\subset 
\underline{Hom} (\Ff / B, V/B )(Z)
$$
to be the full sub-$n$-groupoid of such points. This defines the
substack.

\begin{proposition}
\mylabel{openness}
Under the above notations for $\Ff /B$ and $V/B$ and with the hypothesis
that the $G_Z$ are flat linear group schemes over $Z$,
the sub-$n$-stack $\underline{Hom}^{{\rm se}, 0} (\Ff / B, V/B )$
is an open substack of $\underline{Hom} (\Ff / B, V/B )$.
\end{proposition}
{\em Proof:}
Suppose $Z\rightarrow B$ is a morphism from a scheme, and suppose that
$\varphi : \Ff _Z \rightarrow V_Z$ is a point. We have to show that the
subset of points $z\in Z$ such that $\varphi _z$ is semistable with
vanishing Chern classes, is an open subset of $Z$. For this we can
localize on $Z$ so we assume that there is a basepoint section $v$ of $V_Z$ and
let $G_Z$ be the fundamental group as above. Let $\rho $ be the
principal $G_Z$-bundle over $\Ff _Z$ obtained by the first Postnikov
projection of $\varphi$. We have to prove that the set of points $z\in
Z$ such that $\rho _z$ is semistable with vanishing Chern classes, is
open.  Let $Z_1$ be the inverse image of $B_1$; it is a closed subscheme
of $Z$ and outside of $Z_1$ all of the points are semistable with
vanishing Chern classes. The complement of the set in question is thus
contained in $Z_1$; since $Z_1$ is closed in $Z$, it suffices to prove
that this complement is closed in $Z_1$, in other words it suffices to
prove that the subset of points $z\in Z_1$ where $\rho _z$ is semistable
with vanishing Chern classes, is open in $Z_1$. We are reduced to this
question on $Z_1$ so we may assume $Z=Z_1$. With this assumption $\Ff
_Z=X_{Dol}\times Z$ and $\rho $ is a principal Higgs bundle with
structure group $G_Z$.

When $G_Z$ is the constant sheaf of groups $GL(n)$ over $Z$, the
openness of the set of points which are semistable with vanishing Chern
classes, is well-known. The difficulty here is to extend this to the
case of structure group sheaf $G_Z$ which is flat and linear but not
necessarily constant. 

Let $G_Z\rightarrow GL(E/Z)$ be the faithful representation given by the
condition that $G_Z$ is linear. Let $E(\rho )$ be the associated Higgs
bundle on $X_{Dol}\times Z /Z$. If $z\in Z$ is a point then we obtain
the Higgs bundle $E(\rho )_z$ which is associated to the principal Higgs
$G_z$-bundle $\rho _z$ by the faithful representation $G_z\rightarrow
GL(E_z)$.  Recall that $\rho _z$ is said to be semistable if the Higgs
bundles associated to {\em every} representation of $G_z$ are
semistable;
however, recall from \cite{hbls} p. 86 that if $\rho _z$ has vanishing
Chern classes, then it suffices to check this for one faithful
representation of $G_z$. Therefore the set of points $z$ such that $\rho
_z$ is semistable with vanishing Chern classes, is the same as the set
of points where $\rho _z$ has vanishing Chern classes and where $E(\rho
)_z$ is semistable. The set of points where $E(\rho )_z$ is semistable,
is open cf \cite{moduli}. Thus we may work inside this set, so we may
assume from now on that $E(\rho )_z$ is semistable for all $z\in Z$. 

Let $P\rightarrow X\times Z$ be the geometric principal $G_Z$-bundle
underlying the principal Higgs bundle $\rho$. We have to prove that the
set of points $z\in Z$ where $P_z=P|_{X\times \{ z\} }$ has vanishing
Chern classes, is open.

It is possible to define the Chern classes in de Rham cohomology in an
algebraic way, so the set of points where $P _z$ has vanishing Chern
classes is a constructible set. To prove openness it suffices therefore
to prove that it is open in the usual topology on $Z$ (recall that we
work over $Spec (\cc )$). So chose a point $z$ such that $P _z$ has
vanishing Chern classes. In view of the semistability hypothesis, we
obtain that the principal Higgs bundle $\rho _z$ corresponds to a flat
principal $G_z$-bundle via the correspondence between Higgs bundles and
local systems \cite{hbls}. This correspondence preserves the
$C^{\infty}$ type of the bundle, so we get that the $C^{\infty}$
principal $G_z$-bundle on $X$ underlying $\rho _z$, i.e. $P_z$, is flat. 
By Deligne and Sullivan
\cite{DeligneSullivan},
there is a finite etale covering $X'\rightarrow X$ such that $P
_z|_{X'}$ is a trivial principal bundle, i.e. has a $C^{\infty}$
section $s_z$.

Here is where we  use the flatness hypothesis on $G_Z$. In
characteristic zero flatness implies smoothness for group schemes, so
$G_Z$ is smooth over $Z$. Since $P$ is modelled locally on $G_Z$, we
obtain that the morphism $P\rightarrow X\times Z$ is smooth.
Let $Z'$ be a neighborhood of $z$ in the usual topology in $Z$, such
that there exists a retraction of $Z'$ to $\{ z\}$ where the
trajectories are smooth curves in $Z'$. Assume that the retraction is
obtained by following a vector field along these curves; and assume that
over $Z'$
we have lifted the vector field to a vector field on $P$ (projecting to
zero in the direction $X$ of the base $X\times Z$). Such a neighborhood
exists. 

Let $P '$ be the pullback of $P $ to $X'\times Z'$. Following  our lifts
of vector fields on $P'$ gives a trivialization of the underlying
differential manifold, $P'\cong (P_z\times _XX') \times Z'$. Via this
trivialization, we can extend our section $s_z$ to a section $s$ of the
bundle $P'$. In particular, for every point $y\in Z'$ the $C^{\infty}$
principal $G_y$-bundle $P'_y$ on $X'$ is trivial. In particular it has
vanishing Chern classes in rational cohomology. The morphism $H^i(X, \qq
)\rightarrow H^i(X', \qq )$ is injective, so we obtain that for every
point $y\in Z'$ the principal $G_y$-bundle $P_y$ has vanishing Chern
classes on $X$. Thus $Z'$ provides the neighborhood of $z$ necessary to
show that the set of points where $P_y$ has vanishing Chern classes, is open.
\eop

{\em Counterexample:} 
\mylabel{counterflatness}
The  hypothesis of flatness (i.e. smoothness) of
$G_Z$ is essential in the above argument; here is a counterexample
showing that we cannot remove it. Fix a stable vector bundle $E_0$ with
vanishing Chern classes for structure group $GL(n)$ over a curve $X$.
The bundle $E$ has a Borel reduction, i.e. a flag of sub-line-bundles,
such that the degrees of the line bundles are nonzero. Let $B$ denote
the Borel subgroup which is the structure group for this reduction. Let
$Z=\Aone$
and let $G_Z$ be a group scheme over $Z$, contained in $GL(n)\times Z$,
with generic fiber $B$ (over every point except the origin) 
and special fiber $GL(n)$ over the origin. Note that $G_Z$ isn't flat.
Let $E$ be the bundle $E_0 \times Z$ over $X\times Z$, but considered as
a principal $G_Z$-bundle using the $B$-reduction for $E_0$ outside the
origin. The fiber over the origin is the principal $GL(n)$-bundle $E_0$
which we started with; it is semistable with vanishing Chern classes.
However for any other point $z\neq 0$ the principal $G_z=B$-bundle $E_z$
doesn't have vanishing Chern classes (it isn't semistable either).
Thus the set of points of $Z$ where the bundle is semistable with
vanishing Chern classes, is reduced to $\{ 0\}$ which isn't open.

Next we turn to the problem of local geometricity. This part of the discussion
doesn't use semistability or our special situation for $\Ff /B$; however
we do use the hypothesis that the fundamental group sheaf is a flat linear
group scheme.

\begin{proposition}
\mylabel{localgeometricity}
Suppose $\Ff /B$ is a formal category which is projective and of smooth type
over a base $n$-stack $B$, with finite cohomological dimension for
vector sheaf coefficients. Suppose that $V/B$ is a relatively connected
geometric $n$-stack such that for any scheme $Z\rightarrow B$ with
basepoint section $v: Z\rightarrow V_Z$, the fundamental group sheaf
$G_Z=\pi _1(V_Z/Z,v)$ is represented by a flat linear group scheme.
Then $\underline{Hom}(\Ff /B, V/B)$ is a locally geometric $n$-stack.
\end{proposition}
{\em Proof:}
Suppose $Z\rightarrow B$ is a morphism from a scheme. By going to an
etale neighborhood we may assume that $V_Z$ has a basepoint section $v$,
so we obtain the group sheaf $G_Z$ (which by hypothesis is a
flat linear group scheme over $Z$) and the
first stage in the Postnikov tower for $V_Z$ is $K(G_Z/Z,1)$. Let 
$$
M_Z := \underline{Hom}(\Ff _Z/Z, K(G_Z/Z,1)/Z),
$$
which is the moduli $1$-stack for principal $G_Z$-bundles over $\Ff _Z$.
By \cite{aaspects}, Proposition 10.4.4, $M_Z$ is locally geometric (i.e.
it is an Artin algebraic $1$-stack which is locally of finite type). We
will show that the morphism 
$$
p:\underline{Hom}(\Ff /B, V/B)\times _BZ \rightarrow M_Z
$$
is geometric, yielding the desired conclusion. For this, we may suppose
that
$\rho _Z: Z\rightarrow M_Z$ is a point, and it suffices to show that the
fiber of $p$ over $\rho _Z$ is a geometric $n$-stack over $Z$. The point
$\rho _Z$ may also be considered as a morphism $\rho _Z: \Ff _Z
\rightarrow K(G_Z/Z,1)$ or equivalently a principal $G_Z$-bundle over
$\Ff _Z$. Let $U_Z$ be the fiber in the fibration sequence
$$
U_Z \rightarrow V_Z \rightarrow K(G_Z/Z,1).
$$
It is a simply connected geometric $n$-stack with action of the group
$G_Z$.
Twisting this $n$-stack by the principal $G_Z$-bundle $\rho _Z$ we
obtain a relatively $1$-connected geometric $n$-stack
$$
U_Z(\rho _Z) \rightarrow \Ff _Z.
$$
The fiber of $p$ over $\rho$ is the relative section stack
$$
\Gamma (\Ff _Z/Z, U_Z(\rho _Z)).
$$
We have to prove that this is a geometric $n$-stack over $Z$. This is a
twisted version of the statement 10.3.3 of \cite{aaspects}, and we
follow the idea of the proof that was sketched there. 

The considerations we discuss here are fundamentally well-known
to algebraic topologists.

In order to write the proof in an efficient way, we work with
$\infty$-stacks (which we assume without further stating this to be
stacks of $\infty$-groupoids or equivalently simplicial
presheaves). Suppose $f:R\rightarrow Z$ is a simply connected $\infty$-stack
over a base scheme $Z$. We say that $R$ is {\em cohomologically 
geometric} if ${\bf R}f_{\ast}(\Oo )$ is a perfect complex (supported in
$[0, \infty )$) over $Z$. We say that $R$ is {\em residually geometric}
if for any $N$ the relative truncation $\tau _{\leq N}(R/Z)$ is an
$n$-almost geometric $n$-stack in the terminology of \cite{aaspects} \S 7.3.
Recall the result 7.3.5 of  \cite{aaspects} which, in terms of our
present
definitions, implies that if
$R$ is cohomologically geometric then it is residually geometric. (The
converse is probably also true but we don't have a proof.) Recall also
\cite{aaspects} 7.3.1 saying that if $R$ is $n$-truncated relative to $Z$
and if $R$ is residually geometric then $R$ is in fact geometric.

We note the following canonical construction
(this construction is the ``Goodwillie tower of the identity functor''
see \cite{Goodwillie}). Define 
$$
\Sigma _0:= \Omega ^{\infty} S^{\infty}(R/Z)
$$
to be the Dold-Puppe of the dual  
${\bf R}f_{\ast}(\Oo )^{\ast}$. There is a natural map 
$$
R \rightarrow \Omega ^{\infty}S^{\infty}(R/Z)=\Sigma _0.
$$
If $R$ is cohomologically geometric then $\Sigma _0$
is the Dold-Puppe of a perfect complex so it again is cohomologically
geometric. Let $R_1$ be the fiber of the above map. It fits into a
fibration sequence
$$
\Omega  \Sigma _0
\rightarrow T \rightarrow R,
$$
so an easy spectral sequence argument shows that $T$ is also
cohomologically geometric. Iterating this construction we obtain a
sequence of cohomologically geometric $\infty$-stacks $R_i$ with
fibration sequences
$$
R_{i+1}\rightarrow R_i \rightarrow \Sigma _i
$$
where the $\Sigma _i$ are the Dold-Puppe of perfect complexes.
Furthermore the $R_i$ eventually become more and more connected,
i.e. there is a sequence $k_i\rightarrow \infty$ such that 
$R_i$ are $k_i$-connected. 

The above constructions relative to a base scheme are canonical (with
canonical homotopy coherence at all orders) so they relativize: if $B$
is any $\infty$-stack and if $R/B$ is an $\infty$-stack which is
cohomologically geometric (i.e. for any scheme $Z\rightarrow B$ the
pullback $R_Z=R\times _BZ$ is cohomologically geometric over $Z$) then
we obtain a sequence of $\infty$-stacks $R_i/B$ fitting into fibration
sequences as above, with $\Sigma _i$ again being the relative Dold-Puppe
of perfect complexes relative to $B$, the $R_i$ are
cohomologically geometric, and they are more and more connected relative
to $B$. 

We now apply this construction to the study of relative section stacks.
Suppose $\Ff _Z \rightarrow Z$ is a morphism and suppose $R$ is a
cohomologically geometric $\infty$-stack over $\Ff _Z$. 
We obtain a morphism 
$$
\Gamma (\Ff _Z /Z, R)\rightarrow \Gamma (\Ff _Z /Z, \Sigma _0).
$$
We will treat below the case of sections of the Dold-Puppe of a perfect
complex;
assume for now that these are known to be residually geometric. Thus 
$\Gamma (\Ff _Z /Z, \Sigma _0)$ is residually geometric. The fiber of the
above morphism over any section of $\Gamma (\Ff _Z /Z, \Sigma _0)$ 
is of the form $\Gamma (\Ff _Z /Z, R'_1)$ for an $\infty$-stack $R'_1$ over
$\Ff _Z$ which, over any scheme mapping to $\Ff _Z$, locally looks like the
fiber $R_1$ defined above. Thus the fiber $R'_1$ is cohomologically
geometric relative to $\Ff _Z$ and we can again iterate the construction.
We assumed that $\Ff _Z$  has finite cohomological dimension relative to $Z$
for vector sheaf coefficients. From this assumption, for any fixed $N$ 
the later terms in
the iteration eventually become irrelevant to $\tau _{\leq N}\Gamma (\Ff _Z /Z,
R)$,
and up to these later terms we have expressed $\Gamma (\Ff _Z /Z,
R)$ as the result of a process of fibrations where the base is always a 
residually geometric $\infty$-stack (this remains to be treated below),
and the fiber is again decomposed further in the same way. If the base
and fibers are residually geometric then the total space is residually
geometric
(this is stated with only an indication of the proof in 
\cite{aaspects} Lemma 7.3.3 but we 
leave it as an exercise here too). Finally we obtain that 
$\tau _{\leq N}\Gamma (\Ff _Z /Z,
R)$ is a truncation of a residually geometric $\infty$-stack, thus it is
an $N$-almost geometric $N$-stack. 

If the original $R$ was $n$-truncated relative to $\Ff _Z$ i.e. if it was
an $n$-stack relative to $\Ff _Z$ then choosing $N\geq n+1$ and applying
\cite{aaspects} Lemma 7.3.1, we obtain finally that 
$$
\Gamma (\Ff _ Z /Z,
R)=\tau _{\leq N}\Gamma (\Ff _Z /Z,
R)
$$
is a geometric $n$-stack relative to $Z$.
This applies to our twisted $U_Z(\rho _Z)$ to yield that
$$
\Gamma (\Ff _Z/Z, U_Z(\rho _Z))
$$
is a geometric $n$-stack over $Z$.

To complete the proof we have to treat the case of coefficients in a
perfect complex. If $C$ is a perfect complex over $\Ff _Z$ supported in $(-\infty ,
0]$ and if $\Sigma = DP(C/\Ff _Z )$ then 
$$
\Gamma (\Ff _Z /Z, \Sigma ) = DP \tau ^{\leq 0} \Gamma (\Ff _Z /Z, C).
$$
Thus, if we can show that $\Gamma (\Ff _Z /Z, C)$ is a perfect complex
(supported in $(-\infty , \infty )$) then we obtain that 
$\Gamma (\Ff _Z /Z, \Sigma )$ is residually geometric as desired. 
In fact it suffices to show that $\Gamma (\Ff _Z  /Z, C)$ is
{\em residually perfect} cf \cite{aaspects} Definition 7.3.4.

Now $\Gamma (\Ff _Z /Z, C)$ is calculated by the relative de Rham complex
for $\Ff _Z /Z$, twisted by $C$. There is a spectral sequence whose
starting terms are the cohomology of $X/Z$ with coefficients in $C$
tensored with the vector bundles appearing in the de Rham complex for
$\Ff _Z $ (here $X$ denotes the scheme underlying the formal category
$\Ff _Z$). Thus, in order to obtain the desired result it suffices to
consider a flat morphism of schemes $g:X\rightarrow Z$ and to prove that
if $D$ is a perfect complex on $X$ supported in $(-\infty , m]$
then the direct image $\Gamma (X/Z,
D)= {\bf R}g_{\ast}(D)$ is a residually perfect complex on $Z$. The
residually perfect condition means that for any $N$ the truncation 
$\tau ^{[ -N,N]}{\bf R}g_{\ast}(D)$ into the interval $[-N,N]$
should be quasiisomorphic to the
truncation of a perfect complex. To prove this suppose $k$ is the
cohomological dimension of $X/Z$ and choose an explicit resolution $L$ of $D$ by
vector bundles in the interval $[-N-k, m]$. In other words choose a
complex of vector bundles $L$ (supported in $[-N-k-1, m]$ say)
with a morphism $L\rightarrow D$ which is
exact in the interval $[-N-k, m]$. This is possible since $X/Z$ is flat
and projective. Now 
$$
\tau ^{[ -N,N]}{\bf R}g_{\ast}(L)\rightarrow 
\tau ^{[ -N,N]}{\bf R}g_{\ast}(D)
$$
is a quasiisomorphism. Since $L$ is a  complex of vector bundles over
$X$, the well-known fact that the higher direct image of a vector bundle
under a flat morphism is a perfect complex, implies that  
$\tau ^{[ -N,N]}{\bf R}g_{\ast}(D)$ is the truncation of a perfect
complex.
This being true for any $N$ we obtain that ${\bf R}g_{\ast}(D)$ is
a residually perfect complex. This completes the last part of the proof.
\eop

Finally we go back to our original situation and use the openness of
semistability as well as condition (3) to get a geometric $n$-stack of
semistable morphisms.

\begin{proposition}
\mylabel{boundedness}
Suppose $\Ff /B$ is a formal category which is projective and of smooth type
over a base $1$-stack $B$. Suppose that we are in the situation of (1)-(3)
above (with $B$ geometric). 
Suppose $V/B$ is a geometric $n$-stack which is relatively
connected, and whose relative fundamental group object is a flat linear
group scheme (as above this condition is made relative to every scheme
$Z\rightarrow B$ such that a basepoint section $v$ of $V_Z$ exists).
Then 
$\underline{Hom}^{{\rm se}, 0}(\Ff /B, V/B)$ is a geometric $n$-stack
over $B$.
\end{proposition}
{\em Proof:}
By Proposition \ref{localgeometricity}, the $n$-stack 
$\underline{Hom}(\Ff /B, V/B)$ is locally geometric.
Furthermore, the proof of that proposition actually gives the stronger statement
that
if $\tau _{\leq 1}(V/B)$ denotes the first stage in the Postnikov tower,
then
$$
\underline{Hom}(\Ff /B, \tau _{\leq 1}(V/B)/B)
$$
is locally geometric, and the morphism
$$
\underline{Hom}(\Ff /B, V/B)
\rightarrow 
\underline{Hom}(\Ff /B, \tau _{\leq 1}(V/B)/B)
$$
is geometric. Proposition \ref{openness} says that 
$$
\underline{Hom}^{{\rm se}, 0}(\Ff /B, V/B)\subset 
\underline{Hom}(\Ff /B, V/B)
$$
is an open substack; thus it too is locally geometric. Furthermore 
it is pulled back from the open substack
$$
M:= \underline{Hom}^{{\rm se}, 0}(\Ff /B, \tau _{\leq 1}(V/B)/B)\subset 
\underline{Hom}(\Ff /B,\tau _{\leq 1}(V/B)/B).
$$
Thus we have a geometric morphism
$$
\underline{Hom}^{{\rm se}, 0}(\Ff /B, V/B)
\rightarrow 
\underline{Hom}(\Ff /B,\tau _{\leq 1}(V/B)/B) =M
$$
between locally geometric $n$-stacks,
and it suffices to prove that the target $M$ of this morphism  is geometric,
i.e. of finite type. The reader may be happy to note that $M$ 
is only a $1$-stack---thus it is an Artin algebraic stack
locally of finite type over $B$ and we would like to say that it is of finite
type over $B$. 

We know that over every point $b\in B$, the fiber $M_b$ (which is just a
moduli stack of semistable principal bundles with vanishing Chern classes
on $X_{Dol}$ or $X_{DR}$) is of finite
type. Normally, this type of pointwise information for each fiber isn't sufficient to
conclude that the total stack is of finite type over the base.
However,  we have the  hypothesis (3) that there are only finitely many
equivalence classes of points in the base $B$, and this allows us to
obtain the desired conclusion.

 Let $Z=\coprod _{i\in I} Z_i$ be a disjoint union of schemes $Z_i$
of finite type, with a smooth surjection
$$
Z\rightarrow M.
$$
Let $B(Spec (\cc ))$ denote the finite set of points of $B$, and for
each
$b\in B(Spec (\cc ))$ let $M_b$ denote the preimage in $M$. 
We know (see \cite{moduli}) that each $M_b$ is of finite type, since it
is a moduli stack of semistable principal Higgs $G_b$-bundles of vanishing
Chern classes on $X$ (if $b\in B_1$) or of principal $G_b$-bundles with flat
connection (if $b\in B-B_1$). In
particular, there is a finite subset of indices $I'\subset I$ such that
for any point $b$, $M_b$ is covered by the images of the $X_i$ for $i\in
I'$. Let $X' := \coprod _{i\in I'}X_i$ (it is a scheme of finite type). 
The map $X'\rightarrow M$ is
smooth, and for any $Spec (\cc )$-valued point $m\in M(Spec (\cc ))$ we
have that the preimage of $m$ in $X'$ is nonempty. Thus, for any scheme
$Y$ mapping to $M$, the map $X'\times _MY \rightarrow Y$ has is smooth
with the same pointwise surjectivity property; but $X'\times _MY$ is
a geometric $n$-stack so this means that the map $X' \times
_MY\rightarrow Y$ is a smooth surjection. Therefore $X'\rightarrow M$ is
a smooth surjection, and $M$ is of finite type.
This completes the proof of the proposition.
\eop

As a corollary we obtain the following statement.

\begin{theorem}
\mylabel{flcase}
Suppose $\Vv$ is a pre-namhs (resp. linearized pre-namhs) 
which is connected, and such that the
fundamental group objects such as 
$$
\pi _1(Tot ^{W,F}(V_{DR})/\Aquot
_{wt}\times \Aquot _{hod})
$$
are flat linear group schemes. Suppose $X$
is a smooth projective variety. Then the object $\Hh =
\underline{Hom}(X_M, \Vv )$ is a pre-namhs (resp. linearized pre-namhs).
\end{theorem}
{\em Proof:}
In the construction of $\Hh$ given at the start of this chapter, the
only open question was whether $Tot ^{W,F}(H_{DR})$ was a geometric
$n$-stack over $\Aquot
_{wt}\times \Aquot _{hod}$. Recall that
$$
Tot ^{W,F}(H_{DR}):= 
\underline{Hom}^{{\rm se},0}\left(\frac{X_{Hod}\times \Aquot _{wt}}{
\Aquot_{wt}\times \Aquot_{hod}},
\frac{Tot ^{W,F}(V_{DR})}{\Aquot_{wt}\times \Aquot_{hod}}\right) .
$$
The hypothesis of the theorem says exactly that the hypothesis of 
Proposition \ref{boundedness} applies, so that proposition gives
geometricity of $Tot ^{W,F}(H_{DR})$. The rest of the construction of
$\Hh$, including the treatment of the
linearization, is as indicated at the start of the chapter.
\eop

\begin{center}
{\large \bf  The basic conjectures}
\end{center}
\mylabel{basicconjecturespage}

We can now formulate our basic conjectures concerning nonabelian mixed
Hodge structures.  These three conjectures \ref{nonabcoh}, \ref{internalhom}
and \ref{universal} are all statements we expect that one could prove
within a medium-range amount of time. We decided to post the current version of
the paper as is, without proofs of these conjectures, because it seems clear
(e.g. from the difficulty which we already encountered just treating the
mixed Hodge structures on the homotopy groups $\pi _i(\Vv )$) that their proofs
will require a significant amount of technical material and it didn't seem
necessary to burden the reader---or ourselves---with that at present.

\begin{conjecture}
\mylabel{nonabcoh}
Suppose $X$ is a connected smooth projective variety. Suppose $\Vv$ is a
nonabelian mixed Hodge structure. Then the linearized pre-namhs 
$\Hh = \underline{Hom}(X_M, \Vv )$ is a nonabelian mixed Hodge structure.
\end{conjecture}

Some of the main properties (${\bf Str}$ and ${\bf MHC}$ above) 
are given by Proposition
\ref{mhcatleast}, at least under a reasonable hypothesis on $\Vv$. 
The main things which remain to be shown are: (I) flatness of
the Hodge filtration (condition ${\bf Fl}$) and the conditions
${\bf A1}$ and  ${\bf A2}$ on the annihilator ideals;
and (III) the third Hodge-theoretic condition ${\bf A3}$ on the annihilator
ideals.

A subsidiary conjecture is the following one.

\begin{conjecture}
\mylabel{internalhom}
Suppose $\Vv$ and $\Yy$ are nonabelian mixed Hodge structures, and
suppose that $\Yy$ is $1$-connected (i.e. its $\pi _0$ and $\pi _1$ are
trivial). Then there is a nonabelian mixed Hodge structure 
$\underline{Hom}(\Yy , \Vv )$ (with natural choices for the Hodge and
weight filtrations), again functorial in $\Vv$ and (contravariantly) in
$\Yy$.
\end{conjecture}

Denote by $\ast$ the trivial nonabelian mixed  Hodge structure where all
objects are a single point (objects relative to a base such as 
$\Aquot_{wt} \times \Aquot _{hod}$ are equal to that base). If $\Vv$ is
a namhs then a morphism $\ast \rightarrow \Vv$ is the nonabelian analogue of a
Hodge class of type $(0,0)$ in $\Vv$. Thus we shall call such a point
a {\em Hodge class in $\Vv$}.

The above internal $\underline{Hom}$
should be compatible with the notion of morphism, in that a morphism
from
$\Yy$ to $\Vv$ should be the same thing as a Hodge class in
$\underline{Hom}(\Yy , \Vv )$. Getting back to the first conjecture,
in cases where it holds we can define a {\em  morphism from $X_M$ to
$\Vv$} as being a morphism $\ast \rightarrow \underline{Hom}(X_M, \Vv )$.

We now state another conjecture about representability in the simply
connected case. This conjecture (if it turns out to be true) 
states how to define the ``mixed Hodge
structure on the homotopy type of $X$''.

\begin{conjecture}
\mylabel{universal}
Suppose $X$ is a simply connected smooth projective variety. Then there
is a {\em universal morphism} to a simply connected nonabelian mixed
Hodge structure 
$$
X_M \rightarrow \Yy = {\cal MHS}(X)
$$
with the property that for any nonabelian mixed Hodge structure 
$\Vv$ the resulting morphism
$$
\underline{Hom}(\Yy , \Vv ) \rightarrow \underline{Hom}(X_M, \Vv )
$$
is an equivalence. Furthermore this representing object specializes to
the  $n$-stack representing the de Rham (resp. Betti) shape of $X$ (cf
\cite{aaspects}),
i.e.
\newline
$Y_{DR}$ is the $n$-stack representing the very presentable shape of $X_{DR}$;
\newline
$Tot ^F(Y_{DR})$ is the Hodge filtration on the de Rham representing
object, as constructed in \cite{aaspects}; and 
\newline
$Y_{B,\cc}$ is the $n$-stack representing the very presentable shape of
$X_B$.
\newline
Finally, the mixed Hodge structures on the homotopy vector spaces $\pi
_i^{\rm nu}(\Yy )= \pi _i(X_B)\otimes \cc$ 
coincide with those defined by Morgan and Hain. 
\end{conjecture}

It is a consequence of the functoriality in the first two conjectures
and the universality in this last conjecture that ${\cal MHS}(X)$ would
be functorial in $X$. In other words, if $X\rightarrow Z$ is a morphism
of smooth projective varieties then we would get a morphism of
nonabelian  mixed Hodge structures 
$$
{\cal MHS}(X)\rightarrow {\cal MHS}(Z)
$$
or equivalently a Hodge class in $\underline{Hom}(X_M, {\cal MHS}(Z))$.
It is natural to ask which morphisms of nonabelian mixed Hodge
structures come from morphisms of varieties. This is a subtle
question which could be viewed as a {\em nonabelian analogue of the  Hodge
conjecture}: 
\newline
---do all Hodge classes in 
$\underline{Hom}(X_M, {\cal MHS}(Z))$ come from morphisms of varieties
$X\rightarrow Z$?

In the non-simply connected case, there cannot exist a universal map such as
given by Conjecture \ref{universal}. Nonetheless, one might try
to get back from the nonabelian cohomology to a mixed Hodge structure on the 
``higher Malcev completion'' such as was proposed at the end of \cite{limits}.

Finally we state a conjecture which is intended to show how the condition ${\bf
A1}$ fits in with the result of Deligne-Griffiths-Morgan-Sullivan \cite{dgms}
that smooth projective varieties are formal. 

\begin{conjecture}
\mylabel{formal}
Suppose $X$ is a simply connected finite CW complex, and let $X_B$ denote the constant
$\infty$-stack with values the Poincar\'e $\infty$-groupoid of $X$. 
Then the following two properties are
equivalent: 
\newline
(1) for all $n$ and for any $1$-connected namhs $\Vv$, the pre-weight filtered
$n$-stack
$$
\underline{Hom}(X_B, (V_{B,\cc }, W))
$$
satisfies condition ${\bf A1}$ (cf page \pageref{a1a2});
\newline
(2) the dga associated to $X$ by the theory of Sullivan and Morgan, is formal.
\end{conjecture}

Some vague calculations with dga's seem to support this conjecture. 

Recall in a similar vein the result of the exercise on p. \pageref{defexos}, 
saying that if $Y$
is a quadratic cone at a point $y$ then $(Y, W(y))$ satisfies ${\bf A1}$. One
could also conjecture that that result is ``if and only if'', and extend the
above conjecture to the non-simply connected case. We hope to treat this type of
question in a future paper.

If these statements turn out to be substantially true, then one could say that
the properties of formality; of degeneration of spectral sequences at $E_2$;
and the property ${\bf A1}$ are all sort of the same thing.

\begin{center}
{\large \bf Variations of nonabelian mixed Hodge structure}
\end{center}

\mylabel{fvnamhspage}

Suppose $S$ is a smooth base scheme. We will sketch a definition of 
{\em variation of nonabelian mixed Hodge structure (vnamhs) over $S$}.
The basic idea is to combine the notion of Griffiths transversality for
the Hodge filtration that was defined using the object $S_{Hod}$ in
\cite{santacruz},
with the weight filtration as we are defining it in the present paper.

Most aspects of what we say in this section will be highly conjectural.
Nonetheless we feel that the reader might legitimately be
interested to know about a likely direction of
future development. 

Recall that $S_{Hod}\rightarrow \Aquot _{hod}$ is a formal category  
whose fiber over $[1]$ is $S_{DR}$ and whose fiber
over $[0]$ is $S_{Dol}$ (recall that starting from the $\Gm$-equivariant
formal category over $\Aone _{hod}$
defined in \cite{santacruz} \cite{aaspects}, we took the
quotient
object by $\Gm$ to get a formal 
category over $\Aquot _{hod}$ cf p. \pageref{convenientpage}).
Let $OBJ$ be one of the $n+1$-stacks we consider in the present paper,
and look at objects in $OBJ$.
An object $\underline{\Uu} \rightarrow S_{Hod}$, which may also be seen as a
functor 
$$
\underline{\Uu}:S_{Hod} \rightarrow OBJ,
$$
restricts over $S_{DR}$ to an object $\underline{U}_{DR}$ (or 
$\underline{U}_{DR}:S_{DR}\rightarrow OBJ$),
which is the same thing as an object over $S$ together with a
flat connection.
On the other hand, for every point $s\in S$ we obtain a morphism
$$
\Aquot _{hod} = \{ s\} _{Hod} \rightarrow S_{Hod},
$$
so the restriction of $\Uu$ to $\{ s\} _{Hod}$ gives a filtration $F(s)$ of the fiber
$U_{DR}(s)$. The fact that $\Uu$ is a family over $S_{Hod}$ encapsulates
{\em Griffiths transversality} for this family of filtrations together
with the connection. This is explained in somewhat more detail in
\cite{santacruz}.

Apply the above discussion to the case where $OBJ = F^{DP}.nGEOM$ is the
$n+1$-stack of pre-weight-filtered geometric $n$-stacks. An object 
$$
S_{Hod}\rightarrow F^{DP}.nGEOM
$$
gives a family of linearized weight filtrations with a flat connection
over $S$, and with a varying family of Hodge filtrations (including
Hodge filtrations on the perfect complexes of the linearizations) such
that the Hodge filtrations satisfy Griffiths transversality for the
weight filtration.  Denote the above morphism by
$$
(\underline{U}_{Hod}, W,LGr^W) :S_{Hod}\rightarrow F^{DP}.nGEOM.
$$
The restriction over $S_{DR}$ is denoted $(\underline{U}_{DR}, W, LGr^W)$. 

Note that the ``Hodge filtration'' which we would have called $F$, is
integrated into the object $U_{Hod}$. Thus it is reasonable to think of
the object $\underline{U}_{Hod}$ as being $(\underline{U}_{DR}, F)$. 
On a technical level,
however, $F$ is not a filtration of the object $\underline{U}_{DR}$ 
as this would
imply that $F$ were flat with respect to the connection. The object
$\underline{U}_{Hod}$ represents a family of filtrations $F(s)$ on the fibers 
$\underline{U}_{DR}(s)$, satisfying
Griffiths transversality with respect to the connection. We will however
sometimes use the abuse of notation and call this $(\underline{U}_{DR}, F)$.
This takes care of the part 
$(\underline{U}_{DR}, W,F)$ together
with its linearization $L$, in the
notion of vnamhs.

For the notion of {\em pre-vnamhs} we will want to drop the existence of a
linearization, in this case we just look at a morphism 
$$
(\underline{U}_{Hod}, W) :S_{Hod}\rightarrow F.nGEOM.
$$
If $s\in S$ is a point, then the restriction of $(\underline{U}_{Hod}, W)$ to
$\{ s\} _{Hod} \cong \Aquot _{hod}$ becomes an object 
$$
\Aquot _{hod}\rightarrow F.nGEOM,
$$
i.e. an object of $F.F.nGEOM$. This is just the part which was denoted by
$(\underline{U}_{DR}(s), F(s), W(s))$ in the notion of pre-namhs. 
In the linearized
case we get the same plus a linearization  $(LGr^{W(s)}
(\underline{U}_{DR}(s)), F(s))$.

We next look at the ``Betti'' part of the definition. Recall again that
$S_{B, \rr }$ is the constant stack with values the Poincar\'e groupoid
of $S^{\rm top}$. Since this will be used as the base for a family of
$n$-stacks, we need to go to the $n+1$-truncation of $S^{\rm top}$, in
other words we let $S_{B, \rr }$ denote the constant $n+1$-stack
on the real site, whose values are $\Pi _{n+1}(S^{\rm top})$. Now the
piece 
$(\underline{U}_{B, \rr}, W, L)$ 
in the  notion of vnamhs, corresponds to a functor
of $n+1$-stacks 
$$
(\underline{U}_{B, \rr}, W, LGr ^W): S_{B, \rr} \rightarrow F^{DP}.nGEOM.
$$
This should be seen as a local system over $S^{\rm top}$ of linearized
filtered objects. Again, for the notion of pre-vnamhs we will drop the
existence of the linearization and just look at 
$$
(\underline{U}_{B, \rr}, W): S_{B, \rr} \rightarrow F.nGEOM.
$$
If $s\in S$ is a point then $\{ s\} _{B, \rr }$ is just $\ast$ (i.e.
$Spec (\rr )$) and the restriction of our functor to $\{ s\} _{B, \rr }$
is just a real filtered geometric $n$-stack $(\underline{U}_{B, \rr}(s), W(s))$. 

Finally, recall that we have a morphism 
$$
\zeta _S: S_{DR}^{\rm an} \rightarrow S_{B, \cc }^{\rm an}.
$$
Thus we can ask to have an equivalence (between pre-weight filtered
$n$-stacks over $S_{DR}^{\rm an}$) of the form
$$
\zeta _{\underline{\Uu}}: 
(\underline{U}_{DR}, W, LGr^W)^{\rm an} \cong \zeta _S^{\ast}
(\underline{U}_{B, \rr}, W, LGr ^W)^{\rm an}.
$$
As above, for the notion of pre-vnamhs we will drop the linearization.

We have now indicated all of the elements going into the notion of a {\em
pre-vnamhs} $\underline{\Uu}$ over $S$: this consists of morphisms
$$
(\underline{U}_{Hod}, W) :S_{Hod}\rightarrow F.nGEOM
$$
and 
$$
(\underline{U}_{B, \rr}, W): S_{B, \rr} \rightarrow F.nGEOM,
$$
together with an analytic equivalence of filtered $n$-stacks over $S_{DR}^{\rm an}$,
$$
\zeta _{\underline{\Uu}}: 
(\underline{U}_{DR}, W)^{\rm an} \cong \zeta _S^{\ast}
(\underline{U}_{B, \rr}, W)^{\rm an}.
$$
A {\em linearized pre-vnamhs} $\underline{\Uu}$ is a 
pre-vnamhs with a linearization
in the sense that the object $Gr^W(\underline{\Uu} )$ is provided with an
equivalence to the Dold-Puppe of a ``split pre-vmhc'' (i.e. an object
similar to the above but with $nGEOM$ replaced by $nHPERF$), or
equivalently we include the linearizations in the above data and ask for
$$
(\underline{U}_{Hod}, W,LGr^W) :S_{Hod}\rightarrow F^{DP}.nGEOM
$$
and 
$$
(\underline{U}_{B, \rr}, W, LGr^W): S_{B, \rr} \rightarrow F^{DP}.nGEOM,
$$
again with an analytic equivalence 
$$
\zeta _{\underline{\Uu}}: 
(\underline{U}_{DR}, W, LGr^W)^{\rm an} \cong \zeta _S^{\ast}
(\underline{U}_{B, \rr}, W, LGr ^W)^{\rm an}.
$$

Suppose 
$\underline{\Uu}$ is a pre-vnamhs (resp. linearized pre-vnamhs) over $S$.
As indicated above, if $s\in S$ is a point then the restriction of $\Uu$
to the various objects $\{s\}_{Hod}$ or $\{s\}_{B, \rr }$ yields a
pre-namhs (resp. linearized pre-namhs) $\underline{\Uu} (s)$. 

We think that the following definition is adequate. We say that a {\em
variation of nonabelian mixed Hodge structure $\underline{\Uu}$ over $S$} is a
linearized pre-vnamhs $\underline{\Uu}$ such that for every point $s\in S$ the
restriction $\underline{\Uu} (s)$ is a nonabelian mixed Hodge structure (i.e.
satisfies the axioms of p. \pageref{namhspage}), and such that the
family
of geometric $n$-stacks $Tot^{W}(\underline{U}_{Hod})\rightarrow S_{Hod}$,
pulls back via $S\rightarrow S_{Hod}$ to a flat family over $S$. 
It is possible that it would also be necessary to apply condition ${\bf
A1}$ ``in a family'' over $S$, but we conjecture that that is a
consequence of the simple flatness hypothesis on the whole family.

\begin{center}
{\bf Relative section stacks}
\end{center}

Suppose $Y\rightarrow X \rightarrow B$ is a pair of morphisms of
$n$-stacks (with at least $X$ and $B$ being stacks of groupoids). 
We can define the ``sections of $Y/X$ relative to $B$''
which will be an $n$-stack over $B$, by the formula
$$
\Gamma (X/B, Y) := 
$$
$$
\underline{Hom}\left(\frac{X}{B}, \frac{Y}{B}\right)
\times 
_{\underline{Hom}(\frac{X}{B}, \frac{X}{B})}
B
$$
where the first morphism in the fiber product is the composition with
$Y\rightarrow X$ and the second morphism is the section 
$$
1_{X/B}: B\rightarrow \underline{Hom}(\frac{X}{B}, \frac{X}{B})
$$
corresponding to the identity in each fiber of $X/B$. 

This has the following universal property. If $B'\rightarrow
B$
is a morphism of $n$-stacks, let $X':= X\times _BB'$ and $
Y':= Y\times _BB'$. Then a lifting
$$
B'\rightarrow \Gamma (X/B, Y) \rightarrow B
$$
is the same thing as a section $X'\rightarrow Y'$ i.e. a morphism
composing to the identity of $X'$. 

We can do the same thing from a slightly different viewpoint. Suppose
$X$ and $B$ are $n+1$-stacks (of groupoids) and suppose that
$Y$ is a cartesian family of $n$-stacks over $X$, i.e. a functor
$$
Y: X\rightarrow nSTACK.
$$
Then we again  would like to obtain an object $\Gamma (X/B, Y)$ which, this
time,
should be a
cartesian family of $n$-stacks over $B$. We sketch the construction
here; it runs up against a limitation in \cite{aaspects}, so a rigorous
treatment would require further work. The reader may skip this
discussion;
it is only intended to record for future reference what types of things
we would like to know in the direction of future developpment of the
category-theory side of things. 

Recall  from \cite{aaspects} 
that for any $n+1$-stack $A$ we have an {\em arrow family} which is
a morphism of $n+1$-stacks 
$$
Arr _A : A^o \times A \rightarrow nSTACK.
$$
The further information we need is that this object is {\em functorial
in $A$}. To phrase this precisely, let 
$$
P, S: (n+1)STACK \rightarrow (n+1)STACK
$$
be the functors which respectively associate to an $n+1$-stack $A$, 
$$
P(A):= A^o \times A
$$
and 
$$
S(A):= nSTACK
$$
(thus $S$ is the constant $n+1$-stack with values $nSTACK$). The
statement we need is that the arrow family construction provides a
natural transformation $Arr : P\rightarrow S$. 

Now suppose that $f:A\rightarrow B$ is a morphism of $n+1$-stacks and
$\alpha , \beta : B\rightarrow A$ are sections of $f$. 
Then we obtain the ``morphisms in $A$ from $\alpha$ to $\beta$, relative
to $B$'' which is a cartesian family
$$
(A/B) _{1/}(\alpha , \beta ) : B\rightarrow nSTACK,
$$
constructed as follows using the previous functoriality of $Arr$. 
Think of $A/B$ as a cartesian family which we denote 
$$
\underline{A}: B\rightarrow (n+1)STACK,
$$
and which could be heuristically denoted $b\mapsto A(b)$.
Composing with the functors $P$ (resp. $S$) we obtain cartesian families
$$
P\underline{A} : B\rightarrow (n+1)STACK,\;\;\; b\mapsto A(b)^o\times
A(b),
$$
$$
S\underline{A} : B\rightarrow (n+1)STACK,\;\;\; b\mapsto nSTACK,
$$
and $Arr \circ \underline{A}$ is a natural transformation from 
$P\underline{A}$ to $S\underline{A}$. On the other hand 
$\ast _B$ is the cartesian family which is constant with values the
one-point stack, and $\alpha$ and $\beta$ correspond to natural
transformations
$$
\alpha , \beta : \ast _B \rightarrow \underline{A}.
$$
We obtain a natural transformation $(\alpha , \beta ): \ast _B
\rightarrow P \underline{A}$ and composing with the arrow family gives a
natural transformation
$$
Arr \cdot (\alpha , \beta ) : \ast _B \rightarrow S\underline{A}.
$$
Now go back from the point of view of cartesian families of $n+1$-stacks
over $B$, to the point of view of $n+1$-stacks mapping to $B$: the
cartesian family $\ast _B$ corresponds to the identity $B\rightarrow B$
and the cartesian family $S\underline{A}$ (which is constant with values
$nSTACK$)
corresponds to $nSTACK \times B \rightarrow B$. Our natural
transformation $Arr \cdot (\alpha , \beta )$ 
corresponds to a morphism of $n+1$-stacks
from $B$ to $nSTACK \times B$, and composing with the first projection we obtain
the desired
$$
(A/B) _{1/}(\alpha , \beta ) : B\rightarrow nSTACK.
$$
This construction seems somewhat complicated, for an object which is
``supposed to exist'' by general principles; and our proposed treatment
rests upon the functoriality of $Arr$ which remains to be proven. 
On the other hand, it does seem to be useful for what we are planning to
do below. Thus it would seem worthwhile to take a closer look at this
situation in the future. 

Now we get back to our proposed second construction of the relative
sections stack. If $X\rightarrow B$ is a morphism of $n+1$-stacks (of
groupoids)
then a cartesian family $Y: X\rightarrow nSTACK$ corresponds to a
section 
$$
B\stackrel{Y}{\rightarrow} A:= \underline{Hom}\left(\frac{X}{B}, 
\frac{nSTACK\times
B}{B}\right)\rightarrow B.
$$
The one-point cartesian family $\ast _X: X\rightarrow nSTACK$ similarly
corresponds to a section of $A/B$, and we can put
$$
\Gamma (X/B, Y):= (A/B)_{1/}(\ast _X, Y): B\rightarrow nSTACK
$$
using the above construction. 

There should, in fact, be a third construction of $\Gamma (X/B, Y)$
in the case where $B$ is an $n+2$-stack, $X$ is a cartesian family over
$B$, and $Y/X$ is a family over $B$ of cartesian families over $X$. This
situation isn't formalized yet.

As usual, one would like to have a compatibility between these 
constructions of $\Gamma (X/B, Y)$, but we don't do that here. For the
present, use the first construction which is mathematically sound. The
only problem with that is that it takes place when all the stacks are
$n$-stacks;
but in the Betti case we had to use an $n+1$-stack $X_{B, \rr}$ for the
base of a vnamhs, and technically speaking if we want to take the direct
image along $X\rightarrow S$ then $S_{B,\rr}$ should be an $n+2$-stack. It
suffices to  think of everybody as $n+2$-stacks and apply the first
construction.

\begin{center}
{\bf Extra structure}
\end{center}

Suppose now that $Y$ has extra structure relative to $X$, for example a
filtration or bifiltration, eventually with linearization of one of the
filtrations. Then $\Gamma (X/B, Y)$ will have the same extra structure.
If $(Y,W,F)$ is a bifiltered $n$-stack over $X$ for example, 
then we obtain a bifiltered section stack denoted either of two ways,
$$
(\Gamma (X/B, Y), W,F) = \Gamma (X/B, (Y,W,F)): B\rightarrow F.F.nSTACK.
$$ 
If $Y: X\rightarrow nGEOM$ is a geometric $n$-stack over $X$ then the
relative section stack  might not necessarily be geometric. Under some
circumstances it will be geometric. For example if $X/B$ comes from a formal
category which is projective and of smooth type over $B$, then 
\cite{aaspects} sections 8-10 give some results saying when
$\underline{Hom}(X/B,Y/B)$ is geometric. The treatment of
\cite{aaspects} concerns the case where $Y$ comes from a family over
$B$;
one would like to improve this to treat similar cases where $Y$ is a
family over $X$. This corresponds to the case of local systems along the
fibers of $X/B$. 

The result of this lack of knowledge is that we can't say that the
cartesian families which will appear in our proposed construction of the
relative section object of a vnamhs, will be geometric or not.

\begin{center}
{\bf Relative sections of a vnamhs}
\end{center}

We now come to the situation which interests us. Suppose $f:X\rightarrow
S$ is a projective morphism of smooth varieties. We obtain a morphism of
objects $X_M\rightarrow S_M$, in other words morphisms between all of
the various objects that enter into $X_M$ and $S_M$, with the
appropriate homotopies of compatibility between them. 

We can now state the variational version of our basic Conjecture
\ref{nonabcoh}. 

\begin{conjecture}
\mylabel{vnamhs}
Suppose $\underline{\Uu}$ is a (linearized) pre-vnamhs over $X$.
Then we obtain a (linearized) pre-namhs 
$$
\underline{\Gg} = \Gamma (X_M/S_M, \underline{\Uu} )
$$
over $S$. The objects going into $\underline{\Gg}$ 
are supposed to be the relative
section stacks of the corresponding objects in $\underline{\Uu}$, for example
$$
(\underline{G}_{Hod}, W) := \Gamma ^{\rm se}( X_{Hod}/S_{Hod}, 
(\underline{U}_{Hod}, W)),
$$
and
$$
(\underline{G}_{B, \rr}, W) := 
\Gamma ( X_{B, \rr}/S_{B, \rr}, (\underline{U}_{B, \rr}, W)).
$$
In the first case the superscript `se' indicates that 
we need to restrict to semistable sections over the
Dolbeault points $S_{Dol}\subset S_{Hod}$. 

If $\underline{\Uu}$ is a vnamhs over $X$ then 
$\underline{\Gg}$ should be a vnamhs over $S$.
\end{conjecture}

\newpage

\begin{center}
{\Large \bf Part III: Computations}
\end{center}

Now that we have made our main definitions and described the basic
conjectures, we get to stating what we will actually {\em do} in the
last part of the paper. Our work breaks into two parts:
\newline
(A)\, First we will construct a nonabelian mixed Hodge structure $\Vv$
whose underlying homotopy type is the complexified $2$-sphere
$S^2_{\cc}$
(cf \cite{secondaryKS}). This construction is inspired by the idea
that we would like to have
$$
\Vv = {\cal MHS}(\pp ^1),
$$
but we don't actually prove this latter fact. A little bit more
precisely, we can obtain a morphism $\pp ^1_M \rightarrow \Vv$ and a
number of the properties required in Conjecture \ref{universal} will hold, 
however
we don't treat the main universal property. 
\newline
(B)\, We show the first basic Conjecture \ref{nonabcoh} for the case of
$\underline{Hom}(X_M, \Vv )$ for any smooth projective variety $X$ and
the namhs $\Vv$ constructed in (I) above. We hope that this proof will
provide some ideas for proving Conjecture \ref{nonabcoh} in the general case.

\begin{center}
{\large \bf  Some morphisms between Eilenberg-MacLane pre-namhs}
\end{center}
\mylabel{basicpage}

In the next section we will construct a namhs $\Vv$. It will be the fiber of a
morphism between Eilenberg-MacLane pre-namhs's. Thus, in this section 
we will do a general study of morphisms between Eilenberg-MacLane pre-namhs's.
This is basically the ``Breen calculations'' (see \cite{Breen}, 
\cite{secondaryKS},
\cite{aaspects}) for pre-namhs's.

In this section, which prepares for the construction of $\Vv$ in the
next section, we study some basic morphisms between ``abelian''
pre-namhs's. Recall that a {\em pre-mhs} is a real vector space
$U_{\rr}$
with real weight filtration $W$ and with a ``Hodge filtration'' on the
complexification $(\Vv _{\cc}, F)$. By applying the construction $\xi$ of
\cite{Rees} \cite{naht} \cite{santacruz} we obtain a pre-namhs (for $n\geq 0$)
$$
\Uu = \{ (U_{\cc}, W, F),\;\; (U_{\rr}, W), \;\; 1_{U_{\cc}}\} .
$$
Furthermore, this object is an abelian group object in $0PNAMHS$, which 
allows us to ``deloop'' it and define for any $m\leq n$ the
pre-namhs denoted $K(\Uu , m)$. One way of constructing
this is to let $\Uu [m]$ be the complex with $\Uu$ in degree $-m$ and
zeros elsewhere; this complex is a pre-mhc, and 
$$
K(\Uu , m) := DP (\Uu [m]).
$$
Recall (p. \pageref{shiftconvpage}) that $\Uu$ is a {\em shifted mixed Hodge structure} (of {\em
shift} $s$) if $F$ and
$\overline{F}$
are $k+s$-opposed on $Gr^W_k(U)$. Recall that $\Uu [m]$ will be a mixed
Hodge complex (i.e. satisfies condition ${\bf MHC}$) 
if and only if $\Uu$ is a shifted mhs of shift $-m$.   More generally,
$\Uu [m]$ will be an $s$-shifted mhc if and only if $\Uu$ is 
an $s-m$-shifted mhs.

We call pre-namhs's which are formed out of pre-mhs's in the above way,
``Eilenberg-MacLane pre-namhs's''.
We need to know the ``Breen calculations'' for these guys, namely we
would like to calculate the $\pi _0$ (or other $\pi _i$) of the space of
morphisms from $K(\Uu , m)$ to $K(\Uu ', m')$. Recall that the space of
morphisms is denoted as
$$
PNAMHS _{1/} (K(\Uu , m), K(\Uu ', m')).
$$
A preliminary remark is that (because of the fact that $K(-,i)$
represents the delooping) we have
$$
\pi _iPNAMHS _{1/} (K(\Uu , m), K(\Uu ', m'))=
$$
$$
\pi _0 PNAMHS _{1/} (K(\Uu , m), K(\Uu ', m'-i)).
$$
Thus we only need to calculate the $\pi _0$.

\begin{theorem}
\mylabel{breen1}
Suppose $m,m' \geq 2$, and suppose that $m$ divides $m'$ with $m'=dm$.
Then for $m$ even we have
$$
\pi _0PNAMHS _{1/} (K(\Uu , m), K(\Uu ', m'))=Hom _{MHS}(Sym ^d(\Uu ),
\Uu ')
$$
whereas for $m$ odd we have
$$
\pi _0PNAMHS _{1/} (K(\Uu , m), K(\Uu ', m'))=Hom _{MHS}(\bigwedge ^d(\Uu ),
\Uu ').
$$
Here $Hom _{MHS}$ denotes the set of morphisms in the category of real
mixed Hodge structures.
\end{theorem}
{\em Proof:}
We write the proof  in the case $m$ even; for the case $m$ odd say the
same thing but with $Sym ^d$ replaced by $\bigwedge ^d$.

The space of morphisms in question is a homotopy fiber product of the
form
$$
A\times _CB
$$
where $A$ is the space of morphisms 
$$
Tot ^{F,W}(K(\Uu _{\cc} , m)) \rightarrow 
Tot ^{F,W}(K(\Uu ' _{\cc} , m'))
$$
relative to $\Aquot_{hod}\times \Aquot_{wt}$, where $B$ is the space of
morphisms
$$
Tot ^{W}(K(\Uu _{\rr} , m)) \rightarrow 
Tot ^{W}(K(\Uu ' _{\rr} , m'))
$$
relative to $\Aquot _{wt, \rr}$,
and where $C$ is the space of morphisms 
$$
Tot ^{W}(K(\Uu _{\cc} , m)) \rightarrow 
Tot ^{W}(K(\Uu ' _{\cc} , m'))
$$
relative to $\Aquot _{wt}$. We claim that these spaces of morphisms have
all components simply connected, and have $\pi _0$ equal to the spaces
of morphisms from $Sym ^d(\Uu )$ to $\Uu '$, respectively: for $A$, morphisms of
bifiltered vector spaces; for $B$, morphisms of real filtered vector
spaces; and for $C$,
morphisms of complex filtered vector spaces. 

Before proving the claim we recall the {\em relative Breen calculations}
in general: they say that if $\Ee$ and $\Ee '$ are vector 
bundles over a base $B$ then
$$
\pi _i \underline{Hom} \left(\frac{K(\Ee /B, m)}{B}, \frac{K(\Ee ' /B,
m')}{B} \right)
= H^{m'-i}(K(\Ee /B, m)/B ,\Ee ')
$$
$$
= \Ee ' \otimes _{\Oo _B}
(Sym \; \mbox{or} \; \bigwedge )^{\frac{m'-i}{m}}(\Ee ^{\ast}).
$$
Here the symmetric product is chosen if $m$ is even, the exterior product
if $m$ is odd; and when the exponent is fractional the answer is taken as zero.

Now getting back to the proof of the claim, note
first
that $Tot ^{W,F}(\Uu )$ is a vector bundle over $\Aquot_{hod}\times
\Aquot_{wt}$ and similarly for $\Uu '$; similarly 
$Tot ^W(\Uu _{\rr})$ is a real vector bundle over $\Aquot _{wt, \rr}$
and the same for $\Uu'_{\rr}$; and finally the previous phrase also holds
after complexification of everything.  The Breen calculations in the
relative case over $\Aquot _{hod}\times \Aquot _{wt}$, over $\Aquot
_{wt, \rr}$ or over $\Aquot _{wt}$, then say that the relative
$\underline{Hom}$ from $K(-,m)$ to $K(-', m')$ relative to the
appropriate base, has relative $\pi _0$ equal to the $Hom$ of vector
bundles from the symmetric power of the first vector bundle to the
second vector bundle; and the relative $\pi _1$ vanishes (this is
because of the hypothesis $m\geq 2$ so the terms in the symmetric
algebra are spaced out with at least one zero in between them); and the
higher $\pi _i$ are again vector bundles. Now the bases in question, 
either $\Aquot _{hod}\times \Aquot _{wt}$ or  $\Aquot
_{wt, \rr}$ or $\Aquot _{wt}$, have trivial higher cohomology with
coefficients in vector bundles. Thus only the relative $\pi _0$
contributes to the $\pi _0$ of the space of sections, and the trivial
relative $\pi _1$ gives trivial $\pi _1$ of all components of the space
of sections. Thus we get the simply-connectedness part of the claim, and
the
spaces of sections are the spaces of sections of the bundles given by
the Breen calculations, i.e. the vector bundle $Hom$'s from the
symmetric powers of the bundles for $\Uu$, to the bundles for $\Uu '$.
These spaces of sections are indeed the claimed spaces of morphisms from
$Sym ^d(\Uu )$ to $\Uu '$. This proves the claim.

Now we have identified the space of morphisms in question, as a fiber
product $A\times _CB$, and we have seen that $A$, $B$ and $C$ have
simply connected components. Therefore we have
$$
\pi _0(A\times _CB) = \pi _0(A) \times _{\pi _0(C)}\pi _0(B).
$$
In view of the above descriptions: of $\pi _0(A)$ as the space of
bifiltered morphisms from $Sym ^d (\Uu )$ to $\Uu '$; of 
$\pi _0(B)$ as the space of
real filtered morphisms from $Sym ^d (\Uu )$ to $\Uu '$; and of $\pi _0(C)$ as the space of
complex filtered morphisms from $Sym ^d (\Uu )$ to $\Uu '$, we obtain
that the fibered product is exactly the space of  morphisms of pre-mhs
from $Sym ^d (\Uu )$ to $\Uu '$. This proves the theorem.
\eop

The following theorem could also be stated as a ``conjecture'' since
we have only partially finished the proof; the remainder of the proof is left
as an exercise.

\begin{theorem}
\mylabel{breen2}
Suppose $m,m' \geq 2$, and suppose that $m'=dm-1$.
Then for $m$ even we have
$$
\pi _0PNAMHS _{1/} (K(\Uu , m), K(\Uu ', m'))=Ext^1 _{MHS}(Sym ^d(\Uu ),
\Uu ')
$$
whereas for $m$ odd we have
$$
\pi _0PNAMHS _{1/} (K(\Uu , m), K(\Uu ', m'))=Ext^1 _{MHS}(\bigwedge ^d(\Uu ),
\Uu ').
$$
Here $Ext^1 _{MHS}$ denotes the $Ext^1$ in the category of real
mixed Hodge structures. Outside of the present case and the
case
of the previous theorem, the answer is zero.
\end{theorem}
{\em Sketch of Proof:} Use the same proof as for Theorem \ref{breen1}.
The only place where things change is at the line where it is stated that the
relative $\pi _0$ is given by a certain formula and the relative $\pi _1$
vanishes. In the case treated by the last sentence, both the relative $\pi _0$
and the relative $\pi _1$ of that phrase vanish, so we get that the answer is
zero. In the case $m'=dm+1$ the relative $\pi _0$ vanishes and the relative
$\pi _1$ is given by the relative Breen calculations (it is the same as what occured
as relative $\pi _0$ in the proof of Theorem \ref{breen1}). The remainder of
the proof of Theorem \ref{breen1} may be followed through, bearing in mind that
the homotopy fiber product $A\times _CB$ is now a fiber product of
connected but not simply connected spaces. This gives a certain formula. We
leave as an exercise to the reader to identify the result given by this formula,
with the $Ext$ group in the category of real mixed Hodge structures 
(see Beilinson \cite{Beilinson}) as stated in
the theorem. 
\eop

\begin{center}
{\bf Shifting the weight filtration}
\end{center}

\mylabel{notationpage}

Something which will be useful in the next section is the following
notation.
Let $D\subset \Aquot$ denote the divisor defined by the origin.
If need be we denote the two occurences of this as $D_{hod}\subset
\Aquot _{hod}$ and $D_{wt}\subset \Aquot_{wt}$. Furthermore we also
denote by $D_{wt}$ the divisor 
$$
D_{wt}\times \Aquot_{hod}\subset \Aquot_{wt}\times \Aquot_{hod}.
$$
If $\Uu = (U,W)$ is a filtered vector bundle then we denote by 
$\Uu (-D)$ the filtered bundle whose total bundle of the weight
filtration is 
$$
Tot ^W(\Uu (-D)):= Tot^W(U) \otimes _{\Oo} \Oo _{\Aquot}(-D).
$$
This corresponds to renumbering or ``shifting'' the filtration $W$, namely
\mylabel{divisorpage}
$$
W_{m+1}(\Uu (-D)) = W_{m}(\Uu ).
$$
To see that this is the way the renumbering works, note that $u\in W_m(\Uu )$
if and only if $t^m u$ is a section of $T^W(\Uu )$ 
(see p. \pageref{filtdefpage}). We have the inclusion
of locally free sheaves over $\Aone$
$$
T^W(\Uu (-D)) \hookrightarrow T^W(\Uu ),
$$
and $t^mu$ is a section of $T^W(\Uu )$ if and only if $t^{m+1}u$ is
a section of $T^W(\Uu (-D))$, in other words 
$$
u\in W_m (\Uu ) \Leftrightarrow u\in W_{m+1}(\Uu (-D)).
$$

In the bifiltered situation we adopt a similar notation which gives rise
to the same notation for pre-mhs's: if $\Uu$ is a pre-mhs then we obtain
a new pre-mhs $\Uu (-D_{wt})$ where the weight filtration is shifted by
one.
Note that there is a canonical inclusion
$$
\Uu (-D_{wt})\rightarrow \Uu
$$
which induces the zero map on the associated-graded $Gr^W$, and in fact 
$\Uu (-D_{wt})$ can be characterized as the universal such object: 
any morphism of pre-mhs $\Uu '\rightarrow \Uu$ which induces zero on
$Gr^W$,
factors uniquely through the subobject $\Uu (-D_{wt})$.

\begin{center}
{\large \bf  Construction of a namhs $\Vv$ of homotopy type $S^2_{\cc}$}
\end{center}
\mylabel{aconstructionpage}

We would like to construct our $\Vv$ to have underlying homotopy type $V_B$ equal to
the complexified $2$-sphere $S^2_{\cc}$ (cf \cite{secondaryKS}). This means
that we want $\Vv$
to fit into a Postnikov fibration
sequence
$$
K(\Uu ', 3) \rightarrow \Vv \rightarrow K(\Uu , 2).
$$
The twisting class should be nontrivial (as soon as this class is nontrivial
the underlying homotopy type will be $S^2_{\cc}$).

In fact, we would like $\Vv$ to model the homotopy type of $\pp ^1$ (since
we know, {\em a priori}, that this should have a ``mixed Hodge structure'').

We have to choose the shifted mixed Hodge structures $\Uu $ and $\Uu '$
appropriately; we look at how to do this in a moment. First note that
the above sequence is obtained by specifying the {\em classifying
morphism}
$$
Q: K(\Uu , 2) \rightarrow K(\Uu ', 4).
$$
To be precise, $\Vv$ is the fiber of the  map $Q$, thus it is the
homotopy
fiber product
$$
\Vv = K(\Uu , 2) \times _{K(\Uu ', 4)}\ast .
$$
By the Breen-type calculations of the previous section, this is
equivalent to specifying a morphism of pre-mhs
$$
Sym ^2(\Uu ) \rightarrow \Uu '.
$$
Another thing to note is that we would like the resulting $\Vv$ to
support a linearization; this means that the associated-graded of the
morphism $Q$ for the weight filtration, should be zero. Another thing to
notice is that we would like the base $K(\Uu , 2)$
and the fiber  $K(\Uu ',3)$ to satisfy condition ${\bf MHC}$. Finally,
for our example we would like $\Vv$ to model the homotopy type of $\pp
^2$. Thus $\Uu $ and $\Uu '$ should be one-dimensional. 

In our case of rank one structures, the real structure and the Hodge
filtration determine a pure Hodge structure of some type $(p,p)$; on the
other hand, the weight filtration has only one nontrivial step, at
weight $w$. Thus the pre-mhs's $\Uu $ and $\Uu '$ are determined
respectively by specifying $p$ and $w$ (resp. $p'$ and $w'$). 

The condition that the morphism $Q$ should be zero on the
associated-graded,
is equivalent to saying that $Q$ factors through a morphism
$$
Q': Sym^2(\Uu ) \rightarrow \Uu '(-D_{wt})
$$
(see the end of the previous section p. \pageref{notationpage} 
for the notation $D_{wt}$).

We would like to get a $\Vv$ which models the homotopy type of $\pp ^1$.
During our original research on this
question,
we actually started out by looking at the
homotopy type of $\pp ^1$ and then worked backwards from there to deduce
what conditions needed to be required of a nonabelian mixed Hodge
structure.
In particular, one notices from the $\pp ^1$ example: the fact that one needs to
include the ``homotopy weight'' in the picture (although this was known
from another context of Hodge III \cite{hodge3}); and the condition that the Whitehead
product induce zero on the associated-graded of the weight filtration,
which suggested the general definition that the associated-graded of the
weight filtration should have a structure of a ``spectrum'' or
equivalently a perfect complex. 

To get back to the arithmetic at hand, we need to specify $p$, $w$, $p'$
and $w'$. Since we would like to have $\Uu$ be $\pi _2(\pp ^1)$ and $\Uu
'$ be $\pi _3(\pp ^1)$, we take $p=-1$ and $p' = -2$. This determines
the weights, namely $w=0$ and $w' = -1$. 

This works out correctly for the
condition ${\bf MHC}$: note that $LGr^W(\Vv )$ will be the product of $K(\Uu ,
2)$ and $K(\Uu ', 3)$. We have (according to the choices of $w$ and $w'$) 
$$
H^{-i}(LGr^W_k(\Vv )) =  \begin{array}{l}  
\Uu ,\;\; k=w=0, \; i=2;\\
\Uu ', \;\; k=w'=-1, \; i=3;\\
0 \; \mbox{otherwise}
\end{array}  .
$$
The purity condition ${\bf MHC}$ (cf p. \pageref{puritypage}) says that on 
$H^{-i}(LGr^W_k(\Vv ))$, the Hodge filtration
$F$ and its complex-conjugate $\overline{F}$ (taken with respect to the
real structure), are $m$-opposed with
$$
m = k-i.
$$
Note that according to our choice of $p$ and $p'$, we obtain $m=2p=-2$
for $\Uu$ (which works out as $-2=m=k-i=0-2$) and 
$m=2p'=-4$ for $\Uu '$ (which works out as $-4=m=k-i=-1-3$). Thus, with our
above choices of $p,w,p',w'$ we do indeed obtain the condition ${\bf MHC}$ for
$\Vv$.

\begin{theorem}
\mylabel{constr}
The extension $\Vv$ constructed as the fiber of the morphism $Q$ (with
$p$, $w$, $p'$ and $w,$ chosen as above), is a
nonabelian mixed  Hodge structure.
\end{theorem}
{\em Proof:}
Note that all of the total spaces of filtrations involved in $\Vv$ are
smooth, so the flatness and annihilator conditions are trivial. Our
choice of $p$, $w$, $p'$ and $w'$ gives that 
$K(\Uu , 2)$ and $K(\Uu ', 3)$ are both nonabelian mixed Hodge
structures, i.e. they satisfy the condition ${\bf MHC}$, as calculated
immediately above. Note that they
also satisfy the strictness condition ${\bf Str}$. The linearized
associated graded $LGr^W(\Vv )$ is a direct sum of those for $K(\Uu ,
2)$ and for $K(\Uu ', 3)$ (because the morphism $Q$ is zero on the
associated graded) so it satisfies condition ${\bf MHC}$.
\eop

{\em Remark:}
Note that $K(\Uu
', 4)$ does {\em not} satisfy the ${\bf MHC}$ condition (since it is translated by one
homotopical degree). In particular, we can't say that $Q$ is a morphism
of namhs. On the other hand, $K(\Uu '(-D_{wt}), 4)$ is a namhs, i.e. it
does satisfy the condition ${\bf MHC}$, so the morphism $Q'$ is a
morphism of namhs which makes it ``strict'' in a certain sense. This
phenomenon plays an important role in what goes on at the end of the
proof in the next section.

\begin{center}
{\large \bf  The namhs on cohomology of a smooth projective variety with
coefficients in $\Vv$}
\end{center}
\mylabel{thenamhspage}

Suppose now that $X$ is a smooth projective variety. 
We consider the situation of Conjecture \ref{nonabcoh}, for the particular namhs
$\Vv$ constructed above (Theorem \ref{constr}).

\begin{theorem}
\mylabel{computation}
If $X$ is a smooth projective variety and $\Vv$ is the nonabelian mixed  Hodge
structure constructed in Theorem \ref{constr}, then the linearized pre-namhs
$$
\Hh := \underline{Hom}(X_M, \Vv )
$$
constructed in Theorem \ref{basicconstruction}, 
is a nonabelian mixed Hodge structure.
\end{theorem}
 
The proof occupies the remainder of the paper.
Start by writing
$$
\Hh = \underline{Hom}(X_M, \Vv) =
$$
$$
\underline{Hom}(X_M, K(\Uu , 2) )\times _{\underline{Hom}(X_M, K(\Uu ', 4) )}
\underline{Hom}(X_M,  \ast ).
$$
We have to analyse this fiber product. Note
that
$\underline{Hom}(X_M,  \ast ) = \ast $. 

Put
$$
\Ee := \underline{Hom}(X_M, K(\Uu , 2))
$$
and 
$$
\Ff := \underline{Hom}(X_M, K(\Uu ', 4)).
$$
Also put
$$
\Ff ':= \underline{Hom}(X_M, K(\Uu '(-D_{wt}), 4)),
$$
where $D_{wt}$ was the divisor $\Aquot _{hod}\times [0]_{wt} \subset \Aquot
_{hod}\times
\Aquot_{wt}$. These are all pre-mhc's, and we have a morphism
$i:\Ff ' \rightarrow \Ff$ inducing zero on the associated-graded of the
weight filtration. We have constructed above a morphism (of underlying
pre-namhs's) $Q$ from $K(\Uu , 2)$ to $K(\Uu ', 4)$, which in turn yields
$$
Q_{\Hh}: DP\Ee \rightarrow DP\Ff
$$
factoring through $Q'_{\Hh}:DP\Ee \rightarrow DP\Ff '$. The fact that $i$ induces
zero on the $Gr^W$ implies that $Q_{\Hh}$ extends to a morphism of linearized
pre-namhs's. On the other hand the zero-section is a morphism (again, of
linearized pre-namhs's) 
$$
z: \ast \rightarrow DP\Ff .
$$
We have 
$$
\Hh = \Ee \times _{\Ff} \ast ,
$$
in other words $\Hh$ is the fiber of $Q_{\Hh}$ over the zero-section.

\begin{lemma}
\mylabel{theyremhcs}
In the above situation, $\Ee$ is a mixed Hodge complex and $\Ff$ is a
$1$-shifted mhc. 
\end{lemma}
{\em Proof:}
This is a trivial case of \cite{hodge3}, see also Theorem \ref{hodge3plus}. 
The only difficulty is that
$\Ee$ and $\Ff$ are truncations of the mixed Hodge complexes given by
\cite{hodge3}. Thus, {\em a priori} they are only truncated pre-mhc's.
However, in the present trivial case ($X$ being smooth projective), the
weight filtration is strict too (its spectral sequence degenerates at
$E_1$)
which means that the cohomology objects are bundles over $\Aquot_{wt}$;
thus truncating yields again a complex whose cohomology objects are
bundles, i.e. it is a perfect complex.

For the sign of the shift for $\Ff$, note that $K(\Uu ', 4)$ has
Hodge degree $h=-4$, weight $w= -1$, and homotopical degree $i=4$.
Now $h= (w-i) + 1$ which means that $K(\Uu ', 4)$ is $1$-shifted;
thus $\Ff$ is also $1$-shifted.
\eop

\begin{corollary}
The pre-namhs $\Hh$ satisfies the strictness condition ${\bf Str}$ as
well as the condition ${\bf MHC}$.
\end{corollary}
{\em Proof:}
This was pointed out in Corollary \ref{mhcatleast}; 
for the record we recall how this goes in this particular case. 
The connecting  morphism $Q$ induces zero on the associated-graded
$Gr^W$, so
$$
LGr^W(\Hh ) = LGr^W(\Ee ) \times \Omega LGr^W(\Ff ) .
$$
Here $\Omega$ is the loop-space, which on the level of complexes means
shifting to the right by one and truncating at $0$. As in the lemma, the
mixed Hodge complex $\Ff$ has a strict weight filtration in our easy
situation, so taking the truncation of the shift again yields a pre-mhc.
Recall that $\Ff$ was a $1$-shifted mhc; taking the
loop-space undoes this shift so $\Omega LGr^W(\Ff )$ is a mixed Hodge
complex.
Thus $LGr^W(\Hh )$ is a mixed Hodge complex (which includes the
strictness condition). 
\eop

\begin{center}
{\bf Some general theory}
\end{center}
In order to analyze the annihilator ideals in the fiber product $\Ee
\times _{\Ff}\ast$ we discuss this general type of situation.

\begin{lemma}
\mylabel{smoothnesscrit}
Suppose $f:X\rightarrow Y$ is a morphism of geometric $n$-stacks which
induces an isomorphism on the level of $\pi _0$, and is surjective on
$\pi _1$. Then $f$ is smooth.
\end{lemma}
{\em Proof:}
A morphism between geometric $n$-stacks is smooth if and only if it is
formally smooth (this remark was made e.g. in \cite{aaspects} \S 7.3), 
thus it suffices to prove that $f$ is formally
smooth.
Suppose $A\subset B$ are artinian local schemes.
The problem is to show that the map
$$
X(B) \rightarrow X(A) \times _{Y(A)}Y(B)
$$ 
is surjective on $\pi _0$.
Note that the component maps $X(B)\rightarrow Y(B)$ and 
$X(A)\rightarrow Y(A)$ are isomorphisms on $\pi _0$ and surjections on
$\pi _1$. These mean that the fibers of these morphisms are connected. 
We obtain from the statement for $A$ that the fibers of the map
$$
X(A) \times _{Y(A)}Y(B) \rightarrow Y(B)
$$
are connected, in particular this map induces an isomorphism on $\pi
_0$;
but then from the statement for $B$ we obtain that the map 
$$
\pi _0(X(B)) \rightarrow \pi _0(X(A) \times _{Y(A)}Y(B)) = \pi _0(Y(B))
$$
is an isomorphism. This yields the required lifting property.
\eop

\begin{corollary}
Suppose that $X$ is a geometric $n$-stack and suppose that $\pi _0(X)$
is represented by an algebraic space (i.e. it is a geometric $0$-stack).
Then
the morphism $X\rightarrow \pi _0(X)$ is smooth.
\end{corollary}
{\em Proof:}
Indeed, it clearly satisfies the hypotheses of the previous lemma.
\eop

The same argument gives furthermore

\begin{corollary}
Suppose that $X$ is a geometric $n$-stack and suppose that
there is $k\leq n$ such that  $\tau _{\leq k}(X)$ is a geometric $k$-stack.
Then
the morphism $X\rightarrow \tau _{\leq k}(X)$ is smooth.
\end{corollary}
{\em Proof:}
Again, it satisfies the hypotheses of the previous lemma.
\eop

The above smooth morphisms are useful for the following reason.

\begin{lemma}
Suppose $X\stackrel{g}{\rightarrow} Y\rightarrow
\Aone$ is a morphism of geometric $n$-stacks over $\Aone$.
If $g$ is smooth, then the annihilator ideals for $X$ are the pullbacks of the
annihilator ideals for $Y$.
\end{lemma}
{\em Proof:}
This follows from the definition of annihilator ideals, since a smooth chart
for $X$ gives by composing with $g$ a piece of a smooth chart for $Y$ (if $g$
wasn't surjective it might be necessary to add on some other pieces to the
chart). 
\eop

\begin{corollary}
Suppose $B$ is a scheme, and suppose $X$ and $Y$ are geometric 
$n$-stacks over $B$. Suppose that $U=\pi _0(X)$ and $V=\pi _0(Y)$ are
represented by schemes. Fix a section $b: B\rightarrow Y$ and let $F$ denote
the fiber of $X\rightarrow Y$ over $b$, in other words,
$$
F = X\times _YB.
$$
Let $G$ be the fiber of $U\rightarrow V$ over the projection $pb$ 
(here $p$ denotes the projection from $Y$ to $V$). Then the morphism
$$
F\rightarrow G
$$
is smooth.
\end{corollary}
{\em Proof:}
We can base change everything by the map $B\rightarrow V$, in particular
we can suppose that $V=B$. Note then that $G=U$. With this hypothesis, 
our section $b: B\rightarrow Y$
induces an isomorphism on $\pi _0$, so by Lemma \ref{smoothnesscrit} 
it is smooth. This implies that
$$
F=X\times _YB \rightarrow X
$$
is smooth. On the other hand, the morphism $X\rightarrow U$ is smooth, again
by Lemma \ref{smoothnesscrit}. Thus the composition $F\rightarrow U=G$ is 
smooth. 
\eop

We apply this to the case where the base is $B=\Aone$. We conclude 
that the 
annihilator ideals for the fiber $F$ are the pullbacks of the annihilator ideals
for the fiber $G$ of $\pi _0(X)\rightarrow \pi _0(Y)$ 
(under the hypothesis that these $\pi _0$ are schemes). The same remark holds
for $B=\Aone \times \Aone$ where we calculate the annihilator ideals
with respect to one of the variables. Similarly, flatness of the fiber with
respect to one of the variables is equivalent to flatness of the fiber with
respect to the fiber of the map on the $\pi _0$. 

To conclude, we find that in order to check the annihilator conditions 
and flatness conditions for the fiber of a morphism of linearized pre-namhs,
if the $\pi _0$ are schemes (say for the total space of the Hodge and weight
filtrations over $\Aone \times \Aone$) then it suffices to check these
conditions for the fiber of the induced morphism on $\pi _0$.

\begin{lemma}
Suppose $\Bb$ and $\Cc$ are linearized pre-namhs's, and suppose that 
$f:\Bb \rightarrow \Cc$ is a smooth morphism. Suppose that both $\Bb$
and $\Cc$ satisfy the ``mixed Hodge complex'' condition ${\bf
MHC}$ and strictness ${\bf Str}$.
Then $\Bb$ is a namhs if and only if $\Cc$ is a namhs. 
\end{lemma}
{\em Proof:}
The point is to say that since $f$ is smooth, 
the annihilator ideals of $\Bb$ are the pullbacks of those for $\Cc$
(and also flatness of the Hodge filtration is equivalent on $\Bb$ or
$\Cc$);
thus
the remaining conditions 
${\bf A1}$, ${\bf A2}$, ${\bf Fl}$ and ${\bf A3}$ for $\Bb$ and for
$\Cc$, are equivalent. 
\eop

\begin{center}
{\bf Continuation of the proof of Theorem \ref{computation}}
\end{center}

We analyze our fiber $\Hh$ using the general theory. Note first of all,

\begin{lemma}
In the above situation, the relative zeroth-homologies $H^0(\Ee )$, 
$H^0(\Ff )$ and 
$H^0(\Ff  ')$ are themselves pre-mhs's, in particular 
$$
Tot ^{W,F}H^0(\Ee )_{DR} = H^0(Tot^{W,F}(E_{DR})/\Aquot _{hod}\times
\Aquot _{wt} )
$$
is a vector bundle over $\Aquot _{hod}\times
\Aquot _{wt} $, and the same for $\Ff$ and $\Ff '$.
\end{lemma}
{\em Proof:}
This strictness comes from the fact that, exceptionally in this case, the
spectral sequences for the weight filtration 
degenerate at $E_1$ rather than $E_2$, as in 
Lemma \ref{theyremhcs}.
\eop

Note that 
$$
\pi _0(DP \Ee ) = H^0(\Ee ),\;\;\; \pi _0(DP \Ff ) = H^0(\Ff ).
$$
The morphism $Q$ is not linear so we represent it as a morphism 
between $\pi _0$.

\begin{corollary}
\mylabel{defK}
We get a morphism of linearized pre-namhs's which are sche\-mes
$$
Q_{\Hh ,0}: \pi _0(DP \Ee ) \rightarrow \pi _0(DP \Ff ).
$$
Let $\Kk$ be the fiber of this morphism
over the zero section of $\pi _0(DP \Ff )$. Then 
there is a smooth morphism of linearized pre-namhs's
$$
\Hh \rightarrow \Kk .
$$
\end{corollary}
{\em Proof:} From above.
\eop

\begin{corollary}
In our situation of Corollary \ref{defK}, the 
linearized pre-namhs $\Hh$ is a nonabelian
mixed Hodge structure if and only if $\Kk$ is.
\end{corollary}
\eop

Thus to finish the proof of Theorem \ref{computation} it suffices to show that $\Kk$ is a
namhs. 

We now analyze the situation a bit further. Recall that $\Uu$ is a
pure Hodge structure of Hodge type $(-1,-1)$ with a one-step weight
filtration
placed at weight $0$, and $\Uu '$ is a pure Hodge
structure of type $(-2,-2)$ with weight filtration placed at weight $-1$.

We have
$$
H^0(\Ee ) = H^2(X, \cc ) \otimes \cc ^{-1,-1}\{ 0\}
$$ 
and 
$$
H^0(\Ff ) = H^4(X, \cc ) \otimes \cc ^{-2,-2}\{ -1\} .
$$
The brackets represent the weight of the pre-mhs in question.
Note that $H^0 (\Ee )$ is a $0$-shifted mhs, pure of weight zero.
On the other hand, $H^0(\Ff )$ is a $1$-shifted mhs, 
where the Hodge filtration is pure of weight
zero but the weight filtration is shifted and located at $-1$.

The quadratic morphism $Q_{\Hh ,0}$ is just cup-product, composed with the
inclusion
$\cc \{ 0\} \rightarrow \cc \{ -1\}$. 

Recall that $\Ff ' = \Ff (-D_{wt})$. The morphism $Q_{\Hh , 0}$ factors through a
morphism
$$
Q'_{\Hh ,0}:  \pi _0(DP \Ee ) \rightarrow \pi _0(DP \Ff ')
$$
(this latter being just cup-product, without composition with the 
inclusion $\cc \{ 0\} \rightarrow \cc \{ -1\}$).
Let $\Kk '$ denote the zero-subscheme of $Q'_{\Hh , 0}$.

Note in the above terms that
$$
\pi _0(DP \Ff ')=H^0(\Ff ') = H^4(X, \cc ) \otimes \cc ^{-2,-2}\{ 0\} .
$$
It is a $0$-shifted mixed Hodge structure which is in fact pure of
weight zero (see the renumbering of page \pageref{divisorpage}).
Recall that the same was true of $H^0(\Ee )$.
In particular $Q'_{\Hh , 0}$ is a quadratic morphism of mixed Hodge
structures. We get back to $Q_{\Hh , 0}$ by multiplying $Q'_{\Hh , 0}$
by the coordinate $t$ defining the divisor $D_{wt}\subset \Aquot_{wt}$.

In the above particular case, the total spaces 
$$
Tot^{F,W}H^0(\Ee )\;\;\; \mbox{and} \;\;\;
Tot^{F,W}H^0(\Ff ') 
$$
are bundles over $\Aquot_{hod}\times \Aquot_{wt}$.
The weight filtrations are one-step and in degree zero, so these bundles
are actually pulled back from $\Aquot_{hod}$. Furthermore, 
the Hodge filtrations have splittings by the Hodge decomposition, and
the cup-product is compatible with the splitting. Thus on $\Aone _{hod}$
the bundles in question become trivial bundles, and the quadratic form
$Q'$ is also trivial i.e. it is a product. Therefore the pullback of
$\Kk '$ to $\Aone _{hod}\times \Aone_{wt}$ has a structure of product
$$
\Kk '\times _{\Aquot _{hod}\times \Aquot_{wt}}\Aone _{hod}\times
\Aone_{wt}
$$
$$
\cong \Aone _{hod}\times \Aone_{wt} \times C
$$
where $C$ is the quadratic cone in $H^2(X, \cc )$, kernel of the
cup-product map to $H^4(X, \cc )$. This product decomposition is
compatible with the inclusion into $Tot^{F,W}H^0(\Ee )$ (or rather the
restriction $T(\Ee )$ of this latter over $\Aone _{hod}\times \Aone_{wt}$).
Finally, $\Kk$ is the subscheme of $T(\Ee )$ defined by the
equation of $\Kk '$, multiplied by $t$. Thus $\Kk$ is the
scheme-theoretic union of $\Kk '$ and the inverse image of the divisor
$D_{wt}$. In particular it is trivializable in the $\Aone
_{hod}$-direction,
so $F$ is flat. Also, the structure as union of something flat, with the
inverse image of the reduced divisor $D_{wt}$, means that the
annihilator ideal condition ${\bf A1}$ is satisfied. 
The annihilator ideal $\ker (t)$ is just the ideal defining $\Kk '$,
and $Ann ^W(t; \Hh )$ is the restriction of this over
$D_{wt}=[0]_{wt}\times \Aquot_{hod}$.
Condition ${\bf
A2}$ is satisfied because $\Kk '$ is nonempty. 
Finally, the cup-product morphism is a morphism of mixed Hodge
structures, so the ideal defined by it is a sub-mixed Hodge structure in
the coordinate ring of $H^2(X, \cc )$. This is condition ${\bf A3}$.
\eop

It would seem to be a reasonable project to try to replicate the above types of
computations, for objects $\Vv$ supposed to represent the homotopy types of
projective spaces $\pp ^n$ or Grassmanian varieties. These types of explicit
computations would present an interest of their own, 
independant of any eventual general proof of Conjecture
\ref{nonabcoh}.

\end{document}